\documentclass{amsart}
\usepackage[utf8]{inputenc}

\usepackage{textcmds}
\usepackage{mathtools}
\usepackage{amssymb}
\usepackage{amsthm}
\usepackage{amsmath}
\usepackage{bm}
\usepackage{bbm}
\usepackage{enumerate}
\usepackage{cite}
\usepackage{xfrac}
\usepackage{xcolor}
\usepackage{tikz-cd}
\usepackage{lipsum}
\usepackage{adjustbox}
\usepackage{url}
\usepackage{hyperref}
\usepackage{multicol}

\usepackage[left=3cm,right=3cm,top=3cm,bottom=4.5cm]{geometry}

\theoremstyle{plain}
\newtheorem{theorem}{Theorem}[section]
\newtheorem{corollary}[theorem]{Corollary}
\newtheorem{lemma}[theorem]{Lemma}

\newtheorem{theoremintro}{Theorem}

\theoremstyle{definition}
\newtheorem{definition}[theorem]{Definition}

\theoremstyle{plain}

\theoremstyle{plain}

\theoremstyle{plain}

\theoremstyle{plain}
\newtheorem{fact}[theorem]{Fact}

\theoremstyle{remark}
\newtheorem{remark}[theorem]{Remark}

\theoremstyle{remark}
\newtheorem{notation}[theorem]{Notation}

\theoremstyle{remark}
\newtheorem{example}[theorem]{Example}

\theoremstyle{remark}
\newtheorem{convention}[theorem]{Convention}

\theoremstyle{plain}

\theoremstyle{plain}

\def\dotminussym#1{%
  \setbox0=\hbox{$-$}%
  \kern.5\wd0%
  \hbox to 0pt{\hss\hbox{$-$}\hss}%
  \raise.6\ht0\hbox to 0pt{\hss$.$\hss}%
  \kern.5\wd0%
}

\def\Ind#1#2{#1\setbox0=\hbox{$#1x$}\kern\wd0\hbox to 0pt{\hss$#1\mid$\hss}
\lower.9\ht0\hbox to 0pt{\hss$#1\smile$\hss}\kern\wd0}

\def\ind{\mathop{\mathpalette\Ind{}}}

\def\notind#1#2{#1\setbox0=\hbox{$#1x$}\kern\wd0
\hbox to 0pt{\mathchardef\nn=12854\hss$#1\nn$\kern1.4\wd0\hss}
\hbox to 0pt{\hss$#1\mid$\hss}\lower.9\ht0 \hbox to 0pt{\hss$#1\smile$\hss}\kern\wd0}

\title{Pseudofinite fields with additive and multiplicative character}

\author{Stefan Marian Ludwig}
\address{Albert-Ludwigs-Universität Freiburg,
Mathematisches Institut,
Abteilung für Mathematische Logik,
Ernst-Zermelo-Straße 1,
79104 Freiburg i. B., Germany}
\email{stefan.ludwig@mathematik.uni-freiburg.de}	
\thanks{SML has received funding from the European Union's Horizon 2020 research and innovation programme under the Marie Sk\l{}odowska-Curie grant agreement N\textsuperscript{\underline{o}} 945322.  \includegraphics[scale=0.025]{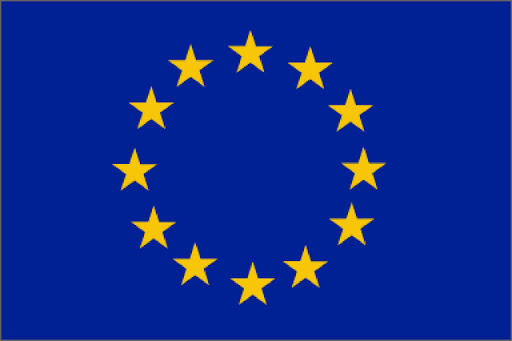}	Moreover, he was partially supported by GeoMod AAPG2019 (ANR-DFG), `Geometric and Combinatorial Configurations in Model Theory'.}



\begin{document}

\begin{abstract} 
We introduce the theory $\mathrm{PF}^{+,\times}$ of pseudofinite fields with generic additive and multiplicative character added as continuous logic predicates. Using the Weil bounds on character sums over finite fields as well as the Erd\H{o}s–Tur\`an–Koksma inequality we show that it is the asymptotic theory (in characteristic $0$) of finite fields with (sufficiently generic) additive and multiplicative character. Moreover, we establish quantifier elimination in a natural definitional expansion of the language and deduce that integration by the Chatzidakis-van den Dries-Macintyre counting measure is uniformly definable in the parameters. Finally, we show that $\mathrm{PF}^{+,\times}$ is a simple theory.
\end{abstract}

\maketitle

%\setcounter{tocdepth}{1}
%{
 % \hypersetup{linkcolor=black}
 % \tableofcontents
%}

% \renewcommand{\thesection}{\thechapter.\arabic{section}}
% \renewcommand{\thesubsection}{\thechapter.\arabic{section}.\arabic{subsection}}

\section{Introduction}\label{sectionintroductionpfplustimes}
Weil's work on the Riemann hypothesis for curves over finite fields and its consequences for character sums in finite fields \cite{weilexponentialsum}\cite{Weil1949NumbersOS} undoubtedly belong to the most influential number-theoretic works of the first half of the 20th century. The estimates for character sums were generalised by Bombieri \cite{bombieri1966exponential} and others in the following decades, for then to culminate in applications of Deligne's work on the Riemann hypothesis over finite fields \cite{Deligne1974}, \cite{Deligne1980}. 
Being unable to give a proper account, we will not go into detail about the numerous application of character sums and their estimates in analytic number theory, additive combinatorics and beyond and refer to presentations such as Chapter 11 of \cite{iwaniec2021analytic}, \cite{Kowalski_2011survey} or \cite{kowlaskiintroductionprobabilisticnumbertheory}.
In this article we will use a variant of Weil's bounds to consider the following model-theoretic question: \textit{What is the common theory of all finite fields together with an additive and multiplicative character?}\\
The model-theoretic study of finite fields goes back to the seminal work of Ax \cite{ax-elttheoryoffinitefields}. He determines the common theory of all finite fields in the language of rings. He showed that a field $F$ is an infinite model of the theory of finite fields (which we will call $\mathrm{PF}$) if and only if $F$ is a \textbf{pseudofinite field}, i.e., $F$ is perfect, PAC (every absolutely irreducible variety defined over $F$ has an $F$-rational point) and $\mathrm{Gal}(F)=\hat{\mathbb{Z}}$. Ax also gave a quantifier elimination result (up to algebraically bounded quantifiers) that was later used and refined by his student Kiefe in order to show a rationality result for zeta functions of definable sets \cite{kiefedefbsetszetafunction} in finite fields.
Another milestone in the study of pseudofinite fields was the 1992 paper \cite{ChatzidakisMacvdd} by Chatzidakis, van den Dries and Macintyre in which the authors establish a generalisation of the Lang-Weil bounds to definable sets over finite fields.
Pseudofinite fields were then used in many applications, often with the above result as a prominent ingredient, to different topics such as definable and algebraic groups over finite fields (e.g. \cite{Hrushovski1994GroupsDI} and
\cite{HrushovskiPillay1995}), motivic integration (see, among others, \cite{denefpasmotivicint}, \cite{nicaisemotInt} and more recently \cite{forey2025motivicfundamentallemma}) and additive combinatorics. See, for example, Tao's paper \cite{tao-expandingpolynomials} using a strengthening of Szemerédi's regularity lemma for definable sets over finite fields.\\
Returning to the context of character sums, in \cite{Kowalski2005ExponentialSO} Kowalski extends certain bounds on character sums over varieties over finite fields to sets definable (in $\mathcal{L}_{\mathrm{ring}})$ in finite fields. He requires the additive character to be non-trivial and discusses why the (purely) multiplicative case is more complicated. While for our purposes the multiplicative case is essential, Kowalski's work however indicates that character sums over finite fields are well suited for a model theoretic treatment. Further occurrences of character sums in model theory are in the work of Tran \cite{tran2019tamestructurescharactersums} who studies the classical first-order theory of the algebraic closure of a finite field $\bar{\mathbb{F}}_{p}$ together with a circular order stemming from the preimage of the circle under a multiplicative character $\chi:\bar{\mathbb{F}}_{q}\rightarrow S^{1}$. Moreover, Tomašić considers in \cite{exposumstomasic} character sums over pseudofinite fields and obtains estimates (Theorem 3.8 in \cite{exposumstomasic}) where the order of the additive and multiplicative characters are finite. In Theorem 6.1 in \cite{indmeasurepseudtomasic} he obtains definability (in the parameters) of summation over those sums.\\
Finally, in his paper on pseudofinite fields with an additive character \cite{Hrushovski2021AxsTW} Hrushovski gives a remarkable model-theoretic account of Kowalski's results. He generalises Ax's main result and determines the characteristic $0$-asymptotic theory of all finite fields with a non-trivial additive character. He works in the language $\mathcal{L}_{+}$ that consist of $\mathcal{L}_{\mathrm{ring}}$ together with a unary, complex-valued continuous logic predicate $\Psi$. The ambient theory is called $\mathrm{PF}^{+}$ and has quantifier elimination in a natural definitional expansion of the language. Furthermore, Hrushovski shows that $\mathrm{PF}^{+}$ is a simple theory. Another important result (Proposition 4.1 in \cite{Hrushovski2021AxsTW}) is the definability of the Chatzidakis-van den Dries-Macintyre counting measure in this context. More precisely, it is shown that the integral (non-standard sum) of an $\mathcal{L}_{+}$-definable predicate is uniformly definable in the ambient parameters.\\
In this article we will generalise those results to pseudofinite fields equipped with both additive and multiplicative character and thus allow for general character sums to be taken into account model-theoretically. Hrushovski's article \cite{Hrushovski2021AxsTW} will be the most prominent building block of this paper, as many of his ideas are also present in our work. While following the same structure, our axiomatisation and language for quantifier elimination will however be more intricate which is essentially due to the fact that it does not suffice to only ensure that the reducts consisting of the field together with the additive (resp. multiplicative) character are models of $\mathrm{PF}^{+}$ (resp. its multiplicative analogue $\mathrm{PF}^{\times}$). Moreover to prove that out theory $\mathrm{PF}^{+,\times}$ is indeed exactly the asymptotic theory of finite fields with additive and multiplicative character (both sufficiently generic) we have to use an approach different from Hrushovski's since his argument as given in \cite{Hrushovski2021AxsTW} does not adapt to the multiplicative case. To remedy this we use a classical result from additive combinatorics, the so-called Erdős–Turán–Koksma inequality. In order to prove uniform definability of integration with respect to the counting measure we show that Hrushovski's proof for $\mathrm{PF}^{+}$ can be adapted to this context despite of the fact that in Kowalski's work the bounds were not obtained for the multiplicative case.\\

\textbf{Presentation of main results.} When we say that a character is sufficiently generic in a pseudofinite field, we mean that it is not of finite order. Note that this is a very natural assumption given that a character of finite order is already interpretable in $\mathrm{PF}$. For example, for a multiplicative character of degree $d$ consider the cosets (in the multiplicative group) of the set defined by $x^{d}=1$.
We work in the language $\mathcal{L}_{+,\times}$, an expansion of the language $\mathcal{L}_{\mathrm{ring}}$ by two continuous logic predicates $\Psi$ and $\chi$ with range $S^{1}\cup \{0\}$ and extend the multiplicative character by $\chi(0)=0$. We give an axiomatisation (Definition \ref{defintionpfplustimes}) of the theory that we call $\mathrm{PF}^{+,\times}$ and obtain (using the Weil bounds) the following. 
\begin{theoremintro}[See Theorem \ref{theoremUPoffinitefieldsmodelofpfplustimes}]\label{theoremintroaxiomatisation}
      Let $\Psi_{q},\chi_{q}$ denote a nontrivial additive and a non-trivial multiplicative character on the finite field $\mathbb{F}_{q}$. Let $(F,\Psi,\chi)=\prod_{\mathcal{U}}(\mathbb{F}_{q},\Psi_{q},\chi_{q})$ be any characteristic $0$ ultraproduct of the structures $(\mathbb{F}_{q},\Psi_{q},\chi_{q})$ where $\chi$ is not of finite order. Then, $F\models \mathrm{PF}^{+,\times}$. 
\end{theoremintro}
Our next main result is quantifier elimination in a natural while technical definitional expansion of $\mathcal{L}_{+,\times}$. The language $\mathcal{L}_{\kappa,+,\times}^{\mathrm{sym}}$ only contains predicates defined with algebraically bounded quantifiers with the added predicates being certain (generalised) finite character sums (see Definition \ref{defintionlanguageforqepfplustimes}).
\begin{theoremintro}[See Theorem \ref{qepfplustimes}]\label{theoremintroqe}
The theory $\mathrm{PF}^{+,\times}$ has quantifier elimination in the language $\mathcal{L}_{\kappa,+,\times}^{\mathrm{sym}}$.
\end{theoremintro}
Our proof essentially has two parts. First, we show that when working in $\mathcal{L}_{\kappa,+,\times}^{\mathrm{sym}}$, we can extend every isomorphism of substructures to one between the corresponding (model-theoretic) algebraic closures (Lemma \ref{lemmapassingtoalgclosurepfplustimes}). Next, in Lemma \ref{lemmaembeddinglemmapfplustimes} we prove an embedding lemma which makes essential use of an equivalent characterisation of models of $\mathrm{PF}^{+,\times}$ (Lemma \ref{lemmapfplustimesequivofdefinition}). As for pure pseudofinite fields, Theorem \ref{theoremintroqe} implies model-completeness of $\mathrm{PF}^{+,\times}$ after adding constants for the coefficients of an irreducible polynomial of degree $n$ for every $n\in\mathbb{N}$ (Corollary \ref{corollarymodelcompleteness}).
Mimicking the corresponding statement for $\mathrm{PF}^{+}$ we do not only eliminate quantifiers but also connectives (Lemma \ref{lemmaalgebraofbasicpredicatespfplustimes}) and show that every quantifier-free formula can already be uniformly approximated by basic predicates (see Definition \ref{definitionbasicpredicatesPFplustimes}).\\
Our next result solves the initial problem of determining the limit theory of finite fields equipped with additive and multiplicative character. As explained above the interesting case is the one of sufficiently generic characters. Note that the additive character being of infinite order below already implies that the field is of characteristic $0$.
\begin{theoremintro}[See Theorem \ref{theorempfplustimesislimittheory}]
Let $T_{\Psi,\chi}$ be the common $\mathcal{L}_{+,\times}$-theory of all finite fields with additive character and multiplicative character.
The theory $\mathrm{PF}^{+,\times}$ is given by $T_{\Psi,\chi}$ together with a set of axioms stating that the multiplicative and additive character are of infinite order.   
\end{theoremintro}
In Theorem \ref{theoremintroaxiomatisation} we have already seen one direction, so the proof reduces to showing that for any completion $T$ of $\mathrm{PF}^{+,\times}$ there indeed is an ultraproduct as in Theorem \ref{theoremintroaxiomatisation} that is a model of $T$. We make use of the proof for pseudofinite fields that invokes the \v{C}ebotarev density theorem and are then left to find the \textit{right} characters for the finite fields in our ultraproduct. Here, we have to use a different strategy as in \cite{Hrushovski2021AxsTW} and deploy the Erdős–Turán–Koksma inequality (Fact \ref{facterdosturankoksmainequality}), essentially a quantitative version of Weyl's equidistribution theorem.\\
Finally, we consider summation of definable predicates over definable sets in finite fields, or equivalently, when working in the limit, integration with respect to the Chatzidakis-van den Dries-Macintyre counting measure. We obtain uniform definabilty of such expressions in the parameters. We state the finitary version of the theorem.
\begin{theoremintro}[See Theorem \ref{corollaryremovinglowerbounds}]
        Let $P(\bar{x},\bar{y})$ be an $\mathcal{L}_{+,\times}$-predicate and $B_{\bar{y}}$ an $\mathcal{L}_{\mathrm{ring}}$-definable set. For any finite field $\mathbb{F}_{q}$ where $q=p^{e}$ consider the function $f:\mathbb{F}_{q}^{|\bar{y}|}\rightarrow\mathbb{C}$ given by $f(\bar{y})=0$, if $B_{\bar{y}}=\emptyset$, and otherwise by\[f(\bar{y}):=\frac{\sum_{\bar{x}\in B_{\bar{y}}(\mathbb{F}_{q})}P(\bar{x},\bar{y})}{|B_{\bar{y}}(\mathbb{F}_{q})|}.\] 
Now for any $\epsilon>0$, there exists an $\mathcal{L}_{+,\times}$-predicate $\Tilde{P}(\bar{y})$  such that for any triple $(\mathbb{F}_{q},\Psi_{q},\chi_{q})$ (consisting of a finite field with additive and multiplicative character) and any $\bar{a}\in \mathbb{F}_{q}^{|\bar{y}|}$,
\[|\Tilde{P}(\bar{a})-f(\bar{a})|\leq\epsilon.\]
\end{theoremintro}
Our argument only uses the Weil bounds and does not rely on Deligne's work in contrast to Kowalski's estimates \cite{Kowalski2005ExponentialSO}. Also, in our work the case only using a multiplicative character is covered, as the difficulties Kowalski describes (e.g., Example 5 in \cite{Kowalski2005ExponentialSO}) are remedied by the fact that the character is sufficiently generic. Our strategy follows \cite{Hrushovski2021AxsTW} adding some of the details to the argument. Finally, using Theorem \ref{theoremintroqe} we can interpret the (non-standard) counting measure as a Keisler measure (Remark \ref{remarkKeislermeasure}).
In the final section we give some results for $\mathrm{PF}^{+,\times}$ that can be proved exactly in the same way as in \cite{Hrushovski2021AxsTW} such as decidablity (in a notion adapted to continuous logic), simplicity and the fact that $\mathrm{PF}^{+,\times}$ is a conservative expansion of $\mathrm{PF}$ in terms of definable sets. Moreover, we briefly discuss the purely multiplicative case (pseudofinite fields with only a multiplicative character) as it might be useful in applications to have results such as definable integration in a language that only contains the multiplicative character.

\subsection*{Structure of the article} 
After stating the necessary estimate on character sums in Section \ref{sectionpreliminariespfplustimes}, we axiomatise the theory $\mathrm{PF}^{+,\times}$ in Section \ref{sectionaxiomatisationpfplustimechapterfive} and show that every characteristic $0$ ultraproduct of finite fields with nontrivial additive and sufficiently generic multiplicative character is a model. In Section \ref{sectionqepfplustimeschapterfive} we obtain quantifier elimination in a natural definitional expansion of the language. Then, we show the other direction of the generalisation of Ax's theorem to this context in Section \ref{sectionlimittheoryplustimes}.
%that is, for every completion of $\mathrm{PF}^{+,\times}$ there is a model given by an ultraproduct as above. This uses a classical result from additive combinatorics, the Erdős–Turán–Koksma inequality. 
In Section \ref{sectiondefinableintegration} we show that integration with respect to the Chatzidakis-van den Dries-Macintyre counting measure is uniformly definable in the parameters to then conclude with some further remarks in Section \ref{sectionfurtheremarkspfplustimes}.

\subsection*{Acknowledgements} This article is based on Chapter 5 of my PhD thesis \cite{ludwig:tel-05236078} written at the École Normale Supérieure Paris. I would like to thank my supervisors Zoé Chatzidakis and Martin Hils for their invaluable support during my thesis.

\section{Preliminaries}\label{sectionpreliminariespfplustimes}
We briefly recall the preliminaries necessary to our work. We start by stating the main number-theoretic input, the \textit{Weil bounds} on character sums for curves over finite fields. Then, we briefly recall the notion of a pseudofinite field and explain how the continuous logic we use in this article corresponds to the one presented in the literature. Finally, we give an overview of the results on the theory $\mathrm{PF}^{+}$ in \cite{Hrushovski2021AxsTW} which we generalise in this article.

\subsection{Character sum estimates over affine curves}We start by recalling some basic facts about additive and multiplicative characters in finite fields.

\begin{definition}
    An additive character $\Psi_{q}$ on $\mathbb{F}_{q}$ is a group homomorphism $\Psi_{q}:(\mathbb{F}_{q},+)\rightarrow S^{1}$. Similarly, a multiplicative character $\chi_{q}$ on $\mathbb{F}_{q}$ is a group homomorphism $\chi_{q}:(\mathbb{F}_{q}^{\times},\cdot)\rightarrow S^{1}$. Whenever necessary we extend it by $\chi_{q}(0):=0$.
\end{definition}

\begin{definition}\label{definitionstandardcharacter}
Given a prime $p$ the \textit{standard character} on $\mathbb{F}_{p}$ is defined by the following map:
\[\Psi_{p}^{\mathrm{stan}}:n+p\mathbb{Z} \;\rightarrow \exp(2\pi i \frac{n}{p}).\]
For an arbitrary finite field $\mathbb{F}_{q}$ with $q=p^{e}$ we define the standard character via
\[\Psi_{q}^{\mathrm{stan}}(x):=\Psi_{p}^{\mathrm{stan}}(\mathrm{Tr}(x))\]
where $\mathrm{Tr}:\mathbb{F}_{q}\rightarrow \mathbb{F}_{p}$ is the trace map.
\end{definition}
\begin{fact}\label{remarkdescriptionalladdcharacters}
Any other additive character on $\mathbb{F}_{q}$ has the form $\Psi_{q}^{\mathrm{stan}}(cx)$ for a unique $c\in\mathbb{F}_{q}$.
\end{fact}

\begin{fact}
    For a finite field $\mathbb{F}_{q}$ and a generator $\gamma_{q}$ of its multiplicative group we denote by $\chi_{\gamma_{q}}$ the multiplicative character determined by sending $\gamma_{q}$ to $e^{2\pi i\frac{1}{q-1}}$. Then the multiplicative characters of $\mathbb{F}_{q}$ are all of the form $\chi_{\gamma_{q}}^{k}$ for some $1\leq k\leq q-1$.
 \end{fact}

In the following we will use an estimate of a general character sum over the rational points of a curve over a finite field. Similar bounds on exponential sums over finite fields date back to Weil \cite{weilexponentialsum} and were later generalised by Bombieri in \cite{bombieri1966exponential}. The estimate in the full generality as below was to our knowledge first presented in the paper \cite{perelmutterexponentialsum} (see its abstract for the statement we use below). Note that this estimate does not use Deligne's work on the Riemann hypothesis for varieties over finite fields \cite{Deligne1974}, \cite{Deligne1980}.

\begin{fact}\label{factperelmutterestimates}
Let $C$ be an absolutely irreducible smooth complete curve in $n$-many variables defined over a finite field $\mathbb{F}_{q}$ and $g,h$ rational functions\footnote{Whose value we set to be $0$ at its poles.} on $C$ that are defined over $\mathbb{F}_{q}$. Let $\Psi_{q}^{\mathrm{stan}}$ be the standard additive and $\chi_{q}$ a non-trivial multiplicative character on $\mathbb{F}_{q}$.\\
Assume that there are no $r,s\in\bar{\mathbb{F}}_{q}(C)$ such that in $\bar{\mathbb{F}}_{q}(C)$ we have both 
$g=s^{q}-s$ and $h=r^{d}$ where $d$ is the order of the multiplicative character $\chi_{q}$. 
Then, there is some natural number $k\in\mathbb{N}$ that only depends on $n$ as well as the degree of the polynomials that define $C$, $g$ and $h$ such that
\[\sum_{\bar{x}\in C(\mathbb{F}_{q})}\Psi_{q}^{\mathrm{stan}}(g(\bar{x}))\chi_{q}(h(\bar{x}))\leq k\cdot q^{\frac{1}{2}}.\]
    
\end{fact}

\subsection{Pseudofinite fields}
We briefly recall the definition of a pseudofinite field. Pseudofinite fields were first introduced and studied by Ax in \cite{ax-elttheoryoffinitefields}. In \cite{ChatzidakisMacvdd} the theory was further developed and the authors obtained a generalisation of the Lang-Weil bounds to definable sets as well as the definability of the (nonstandard) counting measure. To this day, the theory of pseudofinite fields has been actively studied in model theory and beyond, yielding many applications (see, for example, \cite{HrushovskiPillay1995} or \cite{tao-expandingpolynomials}).

\begin{definition}
    A field $F$ is called \textit{pseudofinite}, if it satisfies the following axioms.
    \begin{itemize}
        \item $F$ is perfect.
        \item $F$ is PAC, that is, every absolutely irreducible variety defined over $F$ has an $F$-rational point.
        \item $\mathrm{Gal}(F)=\hat{\mathbb{Z}}$.
    \end{itemize}
    The $\mathcal{L}_{\mathrm{ring}}$-theory that expresses these axioms will be denoted by $\mathrm{PF}$.
\end{definition}

\begin{fact}\label{theoremaxtheorempsf}(Theorem 9 in \cite{ax-elttheoryoffinitefields})
    A field $F$ is pseudofinite if and only if it is infinite and satisfies all $\mathcal{L}_{\mathrm{ring}}$-sentences true in all finite fields.
\end{fact}

\begin{fact}\label{kiefeQE}(See Theorem 2 in \cite{kiefedefbsetszetafunction}.) Every $\mathcal{L}_{\mathrm{ring}}$-formula $\phi(\bar{x})$ is equivalent modulo $\mathrm{PF}$ to a boolean combination of formulas of the form $\exists t\,g(\bar{x},t)=0$ where $g\in\mathbb{Z}[\bar{X},T]$.
    
\end{fact}

\subsection{Continuous logic}

Later we will work with a language in which predicates are allowed to take complex values. Let us briefly explain how this aligns with the usual presentation of real-valued continuous logic as given in \cite{mtfms}.
We can define a \textit{complex continuous logic language} exactly in the same way as for real-valued continuous logic, but associating to any predicate $P$ some compact complex box $C_{P}=\mathrm{Re}(C_{p})\times \mathrm{Im}(C_{p})\subseteq\mathbb{C}$ and working with continuous functions $g:C_{P_{1}}\times\dots\times C_{P_{n}}\rightarrow\mathbb{C}$ as connectives instead. We allow quantification either over the real or over the imaginary part of a complex quantifier-free predicate $P$ using $\inf$ and $\sup$. 
It is then straightforward to see that we can express this in real-valued continuous logic by associating to every ($n$-ary) predicate $P$ in the complex language $\mathcal{L}$ two ($n$-ary) predicates ${P}^{\mathrm{Re}},{P}^{\mathrm{Im}}$ with the corresponding ranges $\mathrm{Re}(C_{P})$ and $\mathrm{Im}(C_{P})$, respectively. One then inductively associates to every $\mathcal{L}$-formula $\phi$ two $\mathcal{L}_{r}$-formulas $\phi^{\mathrm{Re}}$ and $\phi^{\mathrm{Im}}$ such that for all $\bar{a}\in M^{k}$ where $\mathcal{M}$ is an ambient structure we have $\phi(\bar{a})=\phi^{\mathrm{Re}}(\bar{a})+i\cdot \phi^{\mathrm{Im}}(\bar{a})$. From now on, whenever we state that we work in $\mathcal{L}$ we may assume that we work in the corresponding real-valued language $\mathcal{L}_{r}$, and hence all notions and results from real-valued continuous logic can be used without further justification.

\subsection{Pseudofinite fields with additive character}
We briefly recall the main results of Hrushovski's work \cite{Hrushovski2021AxsTW} on the theory $\mathrm{PF}^{+}$ of pseudofinite fields with additive character in characteristic $0$ as the main goal of this article is to generalise them to the setting allowing for a multiplicative character as well.
In Section 2 of \cite{Hrushovski2021AxsTW} it is explained why working with classical discrete logic is problematic (potentially leads to an undecidable theory) when considering the common theory of finite fields with additive characters and thus the following approach using continuous logic is used. 
We work in the language $\mathcal{L}_{+}:=\mathcal{L}_{\mathrm{ring}}\cup\{\Psi\}$ where $\Psi$ is a $S^{1}$-valued continuous logic predicate. Below, a rational hyperplane over $F$ in $\mathbb{A}^{n}$ is a variety defined by an equation of the form $\sum_{1\leq i\leq n}z_{i}X_{i}=b$ where $b\in F$ and $z_{i}\in\mathbb{Z}$ for $1\leq i\leq n$ and $z_{i}\neq 0$ for some $1\leq i\leq n$. We say that it has height $\leq m$, if $|z_{i}|\leq m$ for all $1\leq i\leq n$.

\begin{notation}
    We write $\Psi^{(n)}(\bar{x}):=(\Psi(\bar{x}),\dots,\Psi(\bar{x}))$ (and later similarly for the multiplicative character $\chi^{(n)}$). Moreover, we use $S^{1}$ for the unit circle and $\mathbb{T}^{n}\cong S^{1}\times\cdots\times S^{1}$ for the $n$-dimensional complex torus.
\end{notation}

\begin{definition}(Definition 3.4 in \cite{Hrushovski2021AxsTW})
An $\mathcal{L}_{+}$-structure $(F,\Psi)$ is a model of $\mathrm{PF}^{+}$ if the following axioms hold:
\begin{enumerate}[(1)]
\item $\mathrm{Gal}(F)=\hat{\mathbb{Z}}$ and $\mathrm{char}(F)=0$.
\item $\Psi$ is an additive character on $F$.
\item Let $n,m\in\mathbb{N}$. Let $h\in\mathbb{Q}[z_{1},z_{1}^{-1},\dots,z_{n},z_{n}^{-1}]$ be a finite Fourier series (Laurent polynomial) with degrees $\leq m$, real-valued on $\mathbb{T}^{n}$, with no constant term. For any absolutely irreducible curve $C$ over $F$ with $C\subset\mathbb{A}^{n}$, not contained in any rational hyperplane over $F$ of height at most $m$,
\[\sup\{h(\Psi^{(n)}(\bar{x}))\;:\;\bar{x}\in C(F)\}\geq 0.\]
\end{enumerate}

\end{definition}
The first important result (Proposition 3.9 in \cite{Hrushovski2021AxsTW}) on $\mathrm{PF}^{+}$ which we will generalise in Section \ref{sectionqepfplustimeschapterfive} to $\mathrm{PF}^{+,\times}$ is quantifier elimination in a natural definitional expansion of $\mathcal{L}_{+}$. We refrain at this point from defining the language and instead state the following immediate corollary. 

\begin{corollary}(Consequence of Proposition 3.9 in \cite{Hrushovski2021AxsTW})
$\mathrm{PF}^{+}$ has quantifier reduction up to algebraically bounded quantifiers:
Let $(F_{1},\Psi_{1}),(F_{2},\Psi_{2})$ be two models of $\mathrm{PF}^{+}$
 and $E$ a common relatively algebraically closed subfield on which $\Psi_{1}$ and $\Psi_{2}$ coincide, then
 \[(F_{1},\Psi_{1})\equiv_{E} (F_{2},\Psi_{2}).
 \]
 % In particular, the completions of $\mathrm{PF}^{+,\times}$ are determined by the isomorphism types of $F\cap\bar{\mathbb{Q}}$ together with the corresponding restrictions of $\Psi$ and $\chi$ for $(F,\Psi,\chi)\models\mathrm{PF^{+,\times}}$.
\end{corollary}
The next statement on the asymptotic theory of the structures $(\mathbb{F}_{q},\Psi_{q})$ will be generalised in Section \ref{sectionlimittheoryplustimes}. Our proof however is different from the one in \cite{Hrushovski2021AxsTW} as the proof of a key lemma (Lemma 3.17 in \cite{Hrushovski2021AxsTW}) does not generalise to our context.
\begin{theorem}(Proposition 3.16 in \cite{Hrushovski2021AxsTW})
The theory $\mathrm{PF}^{+}$ is the common theory of the structures $(\mathbb{F}_{q},\Psi_{q})$ together with a set of axioms stating that the characteristic is $0$.
\end{theorem}
Finally, we state the following result on the uniform definability of integration with respect to the counting measure which we will generalise in Section \ref{sectiondefinableintegration} to $\mathrm{PF}^{+,\times}$. It can be considered as a generalisation of the definability of the Chatzidakis-van den Dries-Macintyre counting measure \cite{ChatzidakisMacvdd} to this context. We state the theorem close to its presentation in \cite{Hrushovski2021AxsTW}. A finite version, similar to the one given in the introduction, can be deduced using compactness.

\begin{theorem}(Proposition 4.1 in \cite{Hrushovski2021AxsTW})
  Let $P(\bar{x},\bar{y})$ be an $\mathcal{L}_{+}$-predicate and $B_{\bar{y}}$ an $\mathcal{L}_{\mathrm{ring}}$-definable set. For any finite field $\mathbb{F}_{q}$ where $q=p^{e}$ consider the function $f_{q}:\mathbb{F}_{q}^{|\bar{y}|}\rightarrow\mathbb{C}$ given by $f(\bar{y})=0$, if $B_{\bar{y}}=\emptyset$, and otherwise by\[f_{q}(\bar{y}):=\frac{\sum_{\bar{x}\in B_{\bar{y}}(\mathbb{F}_{q})}P(\bar{x},\bar{y})}{|B_{\bar{y}}(\mathbb{F}_{q})|}.\] 
Now, in an ultraproduct $\prod_{\mathcal{U}}(\mathbb{F}_{q},\Psi_{q})\models \mathrm{PF}^{+}$ the function $f:=\mathrm{lim}_{\mathcal{U}}f_{q}$ is definable.
\end{theorem}

\section{Axiomatisation}\label{sectionaxiomatisationpfplustimechapterfive}

\begin{notation}

Let $\mathbb{G}_{m},\mathbb{G}_{a}$ denote the multiplicative and additive group (of a field) respectively. We say that a morphism $f:\mathbb{G}_{a}^{l}\rightarrow\mathbb{G}_{a}^{k}$ is an \textit{integral linear transformation}, if $\bar{y}=f(\bar{x})$ is given as follows. For $1\leq j\leq k$ we have $y_{j}=a_{1,j}x_{1}+\dots+a_{l,j}x_{l}$ with the $a_{i,j}$ being integers. Similarly, we say that $f:\mathbb{G}_{m}^{l}\rightarrow\mathbb{G}_{m}^{k}$ is an \textit{integral multiplicative transformation}, if $\bar{y}=f(\bar{x})$ is given as follows. For $1\leq j\leq k$ we have $y_{j}=\prod_{1\leq i\leq l}x_{i}^{a_{i,j}}$ with the $a_{i,j}$ being integers. 
    
\end{notation}
%a morphism $g:\mathbb{G}_{m}^{l}\rightarrow\mathbb{G}_{m}^{k}$
As already in \cite{Hrushovski2021AxsTW} (see Section 2 in the same paper for an explanation of why this approach is chosen) for the purely additive case we consider the characters as predicates in continuous logic.
\begin{definition}\label{definitionlanguageplustimes}
    We define as before the language $\mathcal{L}_{+}$ to be the ring language equipped with a unary continuous logic predicate $\Psi$ taking values in $S^{1}\subset\mathbb{C}$. Similarly, we define $\mathcal{L}_{\times}$ to be the ring language equipped with a unary continuous logic predicate $\chi$ taking values in $S^{1}\cup\{0\}\subset\mathbb{C}$. Finally, we define $\mathcal{L}_{+,\times}$ to be the ring language equipped with both $\Psi,\chi$ as before.
\end{definition}
In light of Definition \ref{definitionlanguageplustimes} it would certainly seem natural simply to work inside $\mathbb{G}_{m}$ when dealing with a multiplicative character $\chi$. However, as we will later work with an additive and multiplicative character at the same time, we will introduce the following slightly clumsy notation to make sure that there is no confusion about the number of variables involved in the logical expressions.
 
\begin{notation}
    Given an algebraic set $W\subseteq\mathbb{A}^{n}$, we set $W^{\prime}:=W\backslash Z$ where $Z=Z_{1}\cup\dots\cup Z_{n}$ for $Z_{i}$ defined by the equation $X_{i}=0$ where $\bar{X}=(X_{1},\dots,X_{n})$. If $W\subseteq Z$ we say that $W$ is \textit{zero-degenerate}.
    We will often identify the multiplicative group $\mathbb{G}_{m}^{n}$ with the affine open $(\mathbb{A}^{n})^{\prime}$. This convention will be of importance when we have to keep track of the variables in which we work.
\end{notation}
\begin{convention}
     Sometimes we will consider an integral multiplicative transformation $f:\mathbb{G}_{m}^{l}\rightarrow\mathbb{G}_{m}^{k}$ as a map $f:\mathbb{A}^{l}\rightarrow\mathbb{A}^{k}$. This will always implicitly mean that $f$ is constant $0$ outside of $(\mathbb{A}^{l})^{\prime}$. 
\end{convention}

%\subsection{The multiplicative case}

% for $C^{\prime}(F):=\{\bar{a}\in C(F)\,|\,a_{i}\neq 0\;\text{for all}\; 1\leq i\leq n\}$
%\subsection{The mixed case}

\begin{convention}
    In the following a Laurent-monomial will always be a \textit{primitive} Laurent-monomial. i.e., it is normalised to coefficient $1$. Thus, for us a constant Laurent-monomial will be the constant $1$. A non-trivial Laurent-monomial is then a non-constant one.
%    i.e., one that contains at least one variable in which it has non-zero degree.
\end{convention}

\begin{definition}\label{definitionhyperplanescosetsassociatedtopolynomial}
    Let $h\in\mathbb{Q}[Y_{1},Y_{1}^{-1},\dots,Y_{n},Y_{n}^{-1},Z_{1},Z_{1}^{-1},\dots,Z_{n},Z_{n}^{-1}]$ be a finite Laurent polynomial with degree $\leq m$, real valued on $\mathbb{T}^{2n}=S^{1}\times\dots\times S^{1}$ with no constant term. Write $h=\sum c_{i}h_{i}$ where $c_{i}\in\mathbb{Q}$ and the $h_{i}$ are Laurent monomials and factorise each $h_{i}$ in one monomial $h_{i}^{+}$ in the variables $Y_{1},Y_{1}^{-1},\dots,Y_{n},Y_{n}^{-1}$ and one monomial $h_{i}^{\times}$ in the variables $Z_{1},Z_{1}^{-1},\dots,Z_{n},Z_{n}^{-1}$. As $h$ has no constant term, at least one of the factors $h_{i}^{+}$ or $h_{i}^{\times}$ is non-trivial.\\ 
    Let $F$ be a field. To $h_{i}^{+}(Y_{1},Y_{1}^{-1},\dots,Y_{n},Y_{n}^{-1})=\prod_{1\leq j\leq n}Y_{j}^{s_{i,j}}$ where $s_{i,j}\in\mathbb{Z}$ for $1\leq j\leq n$, we associate the family of rational hyperplanes $D_{i}^{+}(f)$ defined in the variables $x_{1},\dots,x_{n}$ by $\sum_{1\leq j\leq n}s_{i,j}x_{j}=f$ where $f\in F$. 
    Similarly, to $h_{i}^{\times}(Z_{1},Z_{1}^{-1},\dots,Z_{n},Z_{n}^{-1})=\prod_{1\leq j\leq n}Z_{j}^{t_{i,j}}$ where $t_{i,j}\in\mathbb{Z}$ for $1\leq j\leq n$, we associate the family of multiplicative cosets $D_{i}^{\times}(f)$ defined in the variables $x_{1},\dots,x_{n}$ by $\prod_{1\leq j\leq n}x_{j}^{t_{i,j}}=f$ where $f\in F^{\times}$. 
    % We say that $D_{i}^{+}(f)$ (respectively $D_{i}^{\times}(f)$) is non-trivial, if $h_{i}^{+}$ (respectively $h_{i}^{\times}$) is non-trivial.
\end{definition}

\begin{definition}\label{defintionpfplustimes}
    We define $F$ to be a model of the $\mathcal{L}_{+,\times}$-theory $\mathrm{PF}^{+,\times}$, if it satisfies the following axiom scheme:
    \begin{enumerate}[(1)]
        \item $F$ is a field of characteristic $0$  with Galois group $\mathrm{Gal}(F)=\hat{\mathbb{Z}}$.
        \item $\chi: (F^{\times},\cdot)\rightarrow S^{1}$ is a group homomorphism and $\chi(0)=0$.
        \item $\Psi:(F,+)\rightarrow S^{1}$ is a group homomorphism.
        \item The following holds for any $m,m^{\prime}\in\mathbb{N}$:\\
        Let $h\in\mathbb{Q}[Y_{1},Y_{1}^{-1},\dots,Y_{n},Y_{n}^{-1},Z_{1},Z_{1}^{-1},\dots,Z_{n},Z_{n}^{-1}]$ be a finite Laurent polynomial with degree $\leq m$, real-valued on $\mathbb{T}^{2n}=S^{1}\times\dots\times S^{1}$ with no constant term. Write $h=\sum_{1\leq i\leq e}c_{i}h_{i}$ and use the notation of Definition \ref{definitionhyperplanescosetsassociatedtopolynomial}. 
        Further assume that $C\subset\mathbb{A}^{n}$ is a non zero-degenerate, absolutely irreducible curve defined over $F$ of degree $\leq m^{\prime}$, such that for all $1\leq i\leq e$ at least one of the following holds:
        \begin{itemize}
            \item[(+)]  $h_{i}^{+}$ is non-trivial and there is no $f\in F$ such that $C$ is contained in $D_{i}^{+}(f)$, or
             \item[($\times$)] $h_{i}^{\times}$ is non-trivial and there is no $f\in F^{\times}$ such that $C^{\prime}$ is contained in $D_{i}^{\times}(f)$.
        \end{itemize}
        Then, the following holds:
        \[\sup\{h(\Psi^{(n)}(\bar{x}),\chi^{(n)}(\bar{x}))\,|\,\bar{x}\in C^{\prime}(F)\}\geq 0.\]
    \end{enumerate}
\end{definition}

\begin{remark}
    Note that in particular any model $F$ of $\mathrm{PF}^{+,\times}$ is a pseudofinite field. To see this one has only to check that  $F$ is PAC, that is, every absolutely irreducible variety that is define over $F$ has a $F$-rational point. Using Bertini this reduces to curves (see Section 11.2 \cite{fried2008field} for the argument) and then it follows immediately using axiom (4).
\end{remark}

% for $C^{\prime}(F):=\{\bar{a}\in C(F)\,|\,a_{i}\neq 0\;\text{for all}\; 1\leq i\leq n\}$
While the above axiom scheme appears rather technical we will now give an equivalent characterisation in Definition \ref{definitionporpertyequivtoaxiompfplustimes} that will be easier to handle later on. This equivalence (Lemma \ref{lemmafirstdirectionequivaxiom} and \ref{lemmapfplustimesequivofdefinition}) generalises Lemma 3.5 in \cite{Hrushovski2021AxsTW}. Our proof however turns out to be more involved as additive and multiplicative character have to be taken care of at the same time. We start by recalling a result from \cite{Hrushovski2021AxsTW} that will prove useful later on.

\begin{fact}\label{factlaurentpolynomialapproximation}(Lemma 3.5 in \cite{Hrushovski2021AxsTW}, first part of the proof.)
    For any $n\in\mathbb{N}$ and open subset $U\subseteq\mathbb{T}^{n}$, we find a Laurent polynomial with rational coefficients (without constant term) $h$ in $n$-variables that takes real values on $\mathbb{T}^{n}$ such that $h(\bar{a})\geq 0$ implies $\bar{a}\in U$ for $\bar{a}\in\mathbb{T}^{n}$.
\end{fact}

\begin{definition}\label{definitionporpertyequivtoaxiompfplustimes}
    We say that a field $F$ has the property $(\star)$, if the following holds:\\
    Given any non zero-degenerate, absolutely irreducible curve $C\subset \mathbb{A}^{n}$ defined over $F$, $\alpha:\mathbb{A}^{n}\rightarrow\mathbb{A}^{k} $ an integral linear transformation and $\beta:\mathbb{A}^{n}\rightarrow\mathbb{A}^{l}$ an integral multiplicative transformation such that $C_{1}:=\alpha(C^{\prime})$ is not contained in any rational hyperplane over $F$ and $C_{2}:=\beta(C^{\prime})$ is not contained in any multiplicative coset over $F$. Further, let $\Delta(F)$ denote the diagonal belonging to $\alpha,\beta,C^{\prime}$, i.e., $(\bar{x},\bar{y})\in\Delta(F)$ if and only if $(\bar{x},\bar{y})=(\alpha(\bar{z}),\beta(\bar{z}))$ for some $\bar{z}\in C^{\prime}(F)$.\\
    Then, we have that
$(\Psi^{k}\times\chi^{l})(\Delta(F))$ is dense in $\mathbb{T}^{k+l}$.
\end{definition}

\begin{lemma}\label{lemmafirstdirectionequivaxiom}
A field $F$ of characteristic $0$ with $\mathrm{Gal}(F)=\hat{\mathbb{Z}}$ together with an additive character $\Psi$ and a multiplicative character $\chi$ satisfies axiom (4) of the theory $\mathrm{PF}^{+,\times}$ (Definition \ref{defintionpfplustimes}) if it has property $(\star)$ (Definition \ref{definitionporpertyequivtoaxiompfplustimes}).
\end{lemma}

\begin{proof}
 Let $\mathcal{U}\models \mathrm{PF}$ be a monster model containing $F$ as a subfield. Let $\bar{a}=(a_{1},\dots,a_{n})\in\mathcal{U}^{n}$ be a generic point of $C$. The (multiplicative) subgroup $G_{\bar{a}}^{\times}$ generated by $\{a_{1},\dots,a_{n}\}$ in $\mathcal{U}^{\times}/F^{\times}$ is torsion-free (as $F$ is relatively algebraically closed in $\mathcal{U}$) and thus a free abelian group with a free generating set given by $b_{1},\dots,b_{l}$ and $b_{i}=a_{1}^{z_{1,i}}\cdots a_{n}^{z_{n,i}}$ where $(z_{1,i},\dots,z_{n,i})\in\mathbb{Z}^{n}$ for all $1\leq i\leq l$. Thus, we may define an integral multiplicative transformation $\beta:\mathbb{A}^{n}\rightarrow \mathbb{A}^{l}$ by setting $(\beta(\bar{x}))_{i}:=\prod_{1\leq s\leq n}x^{z_{s,i}}$ for $1\leq i\leq l$ and denote by $C_{2}$ the image of $C^{\prime}$ under $\beta$. Then $C_{2}$ is not contained in any multiplicative coset over $F$. Further, as by construction $\{b_{1},\dots,b_{l}\}$ generates $G_{\bar{a}}^{\times}$, we have for all $1\leq i\leq n$ that $a_{i}=\Tilde{\beta_{i}}(\bar{b})\cdot f_{i}$ where $\Tilde{\beta}:=(\Tilde{\beta}_{1},\dots,\Tilde{\beta}_{n})$ is an integral multiplicative transformation $\Tilde{\beta}:\mathbb{A}^{l}\rightarrow \mathbb{A}^{n}$ and $(f_{1},\dots,f_{n})=:\bar{f}\in (F^{\times})^{n}$. Define the map $\gamma:\mathbb{A}^{l}\rightarrow\mathbb{A}^{n}$ by $(\gamma(\bar{x}))_{i}=\Tilde{\beta_{i}}(\bar{x})\cdot f_{i}$ for all $1\leq i\leq n$.\\ 
    Now $\gamma\circ\beta$ yields the identity on $C^{\prime}$. Similarly, working additively (as it was already done in \cite{Hrushovski2021AxsTW}) we find an isomorphism $\alpha: C\rightarrow C_{1}\subset \mathbb{A}^{k}$ given by an integral linear transformation (for the additive case we do not have to pass to $C^{\prime}$ even) such that $C_{1}$ is not contained in any rational hyperplane over $F$. Let $\Delta(F)$ denote the diagonal of $C^{\prime},\alpha,\beta$ as in Definition \ref{definitionporpertyequivtoaxiompfplustimes}. By property $(\star)$ we then obtain that $(\Psi^{k}\times\chi^{l})(\Delta(F))$ is dense in $\mathbb{T}^{k+l}$.\\
    Without loss of generality we assume that $h_{i}^{\times}$ was non-trivial and that $C^{\prime}$ was not contained in $D_{i}^{\times}(f)$ for any $f\in F$. The additive case, that is, when $h_{i}^{+}$ is non-trivial and $C^{\prime}$ is not contained in $D_{i}^{+}(f)$ for any $f\in F$ goes the same way and is already closer to the argument given in Lemma 3.5 in \cite{Hrushovski2021AxsTW}.\\
    For each monomial $h_{i}^{\times}(Z_{1},\dots Z_{n})$ we can consider the formal expression $h_{i}^{\times}\circ\Tilde{\beta}:=H_{i}^{\times}$ as a monomial in $l$-variables. $H_{i}^{\times}$ is itself again non-trivial ($\neq 1$) since otherwise $h_{i}^{\times}$ would be constant on the image of $\gamma$ which contains $C^{\prime}$ and thus $C^{\prime}$ would be contained in the rational coset defined by $h_{i}^{\times}(Z_{1},\dots Z_{n})=h_{i}^{\times}(\bar{f})$ contradicting the assumption on $C$. Consequently, \begin{equation}\label{equationformcommutingwithcharacter}h_{i}^{\times}\circ\chi^{(n)}\circ\gamma(\bar{x})=\chi\circ h_{i}^{\times}\circ\gamma(\bar{x})=\chi(H_{i}^{\times}(\bar{x}))\cdot\chi(H_{i}^{\times}(\bar{f}))=H_{i}^{\times}(\chi^{(l)}(\bar{x}))\cdot \delta_{i} 
    \end{equation}
    where $\delta_{i}=H_{i}^{\times}(\chi^{(l)}(\bar{f}))\in S^{1}$. In the first expression of equation (\ref{equationformcommutingwithcharacter}) $h_{i}$ is evaluated in $\mathbb{C}^{n}$ and in the second on $F^{n}$. Similarly $H_{i}^{\times}$ is evaluated in $F^{l}$ in the third and in $\mathbb{C}^{n}$ in the fourth expression.\\
    Now, we can consider the last expression as a multiple of a Laurent-monomial over $\mathbb{C}$ and set $g_{i}^{\times}:=\delta_{i}\cdot H_{i}^{\times}\in \mathbb{C}[S_{1},S_{1}^{-1},\dots,S_{l},S_{l}^{-1}]$.
    Since $g_{i}^{\times}$ is non-trivial we have \begin{equation}\label{eqationintegrationgitimes}\int_{\mathbb{T}^{l}}g_{i}^{\times}=0\end{equation} where we consider integration with respect to the Haar measure on $\mathbb{T}^{l}$.
    Finally, we will put everything together.
    Let the monomial $g_{i}^{+}\in \mathbb{C}[T_{1},T_{1}^{-1},\dots,T_{k},T_{k}^{-1}]$ be obtained from $h_{i}^{+}$ with the same procedure. Note, that at this point $g_{i}^{+}$ might be trivial, because we only worked with the assumption that $C^{\prime}$ was not contained in any $D_{i}^{\times}(f)$.\\
    Nevertheless, if we consider $g_{i}\in \mathbb{C}[T_{1},T_{1}^{-1},\dots,T_{k},T_{k}^{-1},S_{1},S_{1}^{-1},\dots,S_{l},S_{l}^{-1}]$ to be given by the product $g_{i}^{\times}(S_{1},S_{1}^{-1},\dots,S_{l},S_{l}^{-1})\cdot g_{i}^{+}(T_{1},T_{1}^{-1},\dots,T_{k},T_{k}^{-1})$ we obtain from equation (\ref{eqationintegrationgitimes}) that \begin{equation}\label{equationintegrationgi}\int_{\mathbb{T}^{k+l}}g_{i}=\int_{\mathbb{T}^{k+l}} g_{i}^{+}\cdot g_{i}^{\times}= \left(\int_{\mathbb{T}^{k}}g_{i}^{+}\right)\cdot\left(\int_{\mathbb{T}^{l}}g_{i}^{\times}\right)=0.\end{equation}
    We return to consider the monomial $h_{i}$ that we factorised as $h_{i}=h_{i}^{\times}\cdot h_{i}^{+}$. We obtain the following for any $\bar{x}\in C^{\prime}(F)$ from equation (\ref{equationformcommutingwithcharacter}) and its additive analogue. (Recall that $\beta$ is the inverse of $\gamma$.)
    \begin{equation}\label{equationproductofhiandgithesame}h_{i}^{\times}(\chi^{(n)}(\bar{x}))\cdot h_{i}^{+}(\Psi^{n}(\bar{x}))=g_{i}^{\times}(\chi^{(l)}(\beta(\bar{x})))\cdot g_{i}^{+}(\Psi^{(k)}(\alpha(\bar{x}))) \end{equation}
    Recall that we had $h=\sum c_{i}h_{i}$. Define $g:=\sum c_{i}g_{i}$. From equation (\ref{equationproductofhiandgithesame}) we get for every $\bar{x}\in C$ that
    \begin{equation}\label{equationhequalsg}
        h(\Psi^{(n)}(\bar{x}),\chi^{(n)}(\bar{x}))=g(\Psi^{(k)}(\alpha(\bar{x})),\chi^{(l)}(\beta(\bar{x}))).
    \end{equation}
As $h$ was assumed to be real-valued on $\mathbb{T}^{2n}$, it follows from equation (\ref{equationhequalsg}) that $g$ is real-valued on $\Psi^{(k)}\times\chi^{(l)}(\Delta)\subseteq \mathbb{T}^{k+l}$. But now, we use that we obtained from property $(\star)$ that $\Psi^{(k)}\times\chi^{(l)}(\Delta)$ is dense in $\mathbb{T}^{k+l}$. It follows, that $g$ is real-valued on $\mathbb{T}^{k+l}$. Further, we deduce from equation (\ref{equationintegrationgi}) that $\int_{\mathbb{T}^{k+l}}g=0$ and thus $\sup_{\mathbb{T}^{k+l}}g\geq 0$. Combining this with the fact that $\Psi^{(k)}\times\chi^{(l)}(\Delta(F))$ is dense in $\mathbb{T}^{k+l}$ we finally can conclude that
\begin{equation}
    \sup_{\bar{x}\in C^{\prime}(F)}h(\Psi^{(n)}(\bar{x}),\chi^{(n)}(\bar{x}))=\sup_{\bar{x}\in C^{\prime}(F)}g(\Psi^{(k)}(\alpha(\bar{x})),\chi^{(l)}(\beta(\bar{x})))\geq 0
\end{equation}
which completes the proof of the first direction.
\end{proof}

\begin{lemma}\label{lemmapfplustimesequivofdefinition}
    If a field $F$ of characteristic $0$ with additive character $\Psi$ and multiplicative character $\chi$ satisfies axiom (4) of the theory $\mathrm{PF}^{+,\times}$, then it has property $(\star)$. In particular, the two notions are equivalent.
\end{lemma}
\begin{proof}
Let $C,\alpha=(\alpha_{1},\dots,\alpha_{k}),\beta=(\beta_{1},\dots,\beta_{l}),C_{1},C_{2}$ be given as in property $(\star)$. We want to show that axiom (4) of $\mathrm{PF}^{+,\times}$ implies that $(\Psi^{k}\times\chi^{l})(\Delta(F))$ is dense in $\mathbb{T}^{k+l}$. We consider $\alpha_{j}$ for $1\leq j\leq k$ and $\beta_{i}$ for $1\leq i\leq l$ as Laurent monomials in $n$ variables $X_{1},X_{1}^{-1},\dots,X_{n},X_{n}^{-1}$. That is, if $\alpha_{j}(\bar{x})=\sum_{1\leq s\leq n} z_{j,s}x_{s}$ we associate to it the monomial $\prod_{1\leq s\leq n} X_{s}^{z_{j,s}}$ and similarly for the $\beta_{i}$. Let $U\subseteq\mathbb{T}^{k+l}$ be an open subset.
By Fact \ref{factlaurentpolynomialapproximation} we obtain a real-valued (on $\mathbb{T}^{k+l}$) Laurent polynomial $g(T_{1},T_{1}^{-1},\dots,T_{k},T_{k}^{-1},S_{1},S_{1}^{-1},\dots,S_{l},S_{l}^{-1})$ with rational coefficients and no constant term such that $g(\bar{t})\geq 0$ implies $\bar{t}\in U$ for $\bar{t}\in\mathbb{T}^{k+l}$. Now we replace in $g$ the variables $T_{j}$ for $1\leq j\leq k$ by the monomials $\alpha_{j}(Y_{1},Y_{1}^{-1},\dots,Y_{n},Y_{n}^{-1})$ and the variables $S_{i}$ for $1\leq i\leq l$ by the monomials $\beta_{i}(Z_{1},Z_{1}^{-1},\dots,Z_{n},Z_{n}^{-1})$.\\
We name the polynomial from $\mathbb{Q}[Y_{1},Y_{1}^{-1},\dots,Y_{n},Y_{n}^{-1},Z_{1},Z_{1}^{-1},\dots,Z_{n},Z_{n}^{-1}]$ that we obtain like this by $h$. Then $h$ remains real valued on $\mathbb{T}^{2n}$. Next, we show that $h$ has no constant term and $C$ satisfies the condition in axiom (4) of $\mathrm{PF}^{+,\times}$ with respect to $h$, i.e., for any monomial of $h$ $(+)$ or $(\times)$ (here even \textit{and} instead of \textit{or} in case both $h_{i}^{+}$ and $h_{i}^{\times}$ are non-trivial) is fullfilled:\\ 
Assume that not, then $g_{i}^{+}$ (resp. $g_{i}^{\times}$) would be constant on the image of $\alpha$ (resp. $\beta$). As the later contains $C_{1}$ (resp. $C_{2}^{\prime})$, it follows that $C_{1}$ (resp. $C_{2}^{\prime}$) would be contained in the rational hyperplane (resp. rational coset) associated to $g_{i}^{+}$ (resp. $g_{i}^{\times}$).
%$C_{1}$ or $C_{2}$ would be contained in a rational hyperplane over $F$ (rational coset over $F$, respectively) as $g_{i}^{+}$ (resp. $g_{i}^{\times}$) would be constant on the image of $\alpha$ (resp. $\beta$) which contains $C_{1}^{\prime}$ (resp. $C_{2}^{\prime})$ (and thus $C_{1}$ (resp. $C_{2}^{\prime}$) would be contained in the rational hyperplane (resp. rational coset) associated to $g_{i}^{+}$). Note that the latter is almost the same argument as in the first direction, only working with the inverses of the maps we worked with before.
Again, with the same argument as in the first direction we obtain that
    \[  h(\Psi^{(n)}(\bar{x}),\chi^{(n)}(\bar{x}))=g(\Psi^{(k)}(\alpha(\bar{x})),\chi^{(l)}(\beta(\bar{x}))).\]
And now by axiom (4) of $\mathrm{PF}^{+,\times}$ applied to $C$ and $h$, we obtain that     
\[
    \sup_{\bar{x}\in C^{\prime}(F)}g(\Psi^{(k)}(\alpha(\bar{x})),\chi^{(l)}(\beta(\bar{x})))=\sup_{\bar{x}\in C^{\prime}(F)}h(\Psi^{(n)}(\bar{x}),\chi^{(n)}(\bar{x}))\geq 0
\]and thus $(\Psi^{k}\times\chi^{l})(\Delta(F))\cap U$ is non-empty. As $U$ was chosen arbitrarily, it follows that $(\Psi^{k}\times\chi^{l})(\Delta(F))$ is dense in $\mathbb{T}^{k+l}$.  
\end{proof}

\begin{theorem}\label{theoremUPoffinitefieldsmodelofpfplustimes}
    Let $\Psi_{q},\chi_{q}$ denote a non-trivial additive and a non-trivial multiplicative character on the finite field $\mathbb{F}_{q}$. Let $(F,\Psi,\chi)=\prod_{\mathcal{U}}(\mathbb{F}_{q},\Psi_{q},\chi_{q})$ be any characteristic $0$ ultraproduct of the structures $(\mathbb{F}_{q},\Psi_{q},\chi_{q})$ where $\chi$ is not of finite order. Then, $F\models \mathrm{PF}^{+,\times}$.
\end{theorem}
\begin{proof}
    We only have to show that axiom (4) of Definition \ref{defintionpfplustimes} holds asymptotically (with respect to the ultrafilter $\mathcal{U}$) for the $\mathbb{F}_{q}$. When we write \textit{for almost all} in the following, this will always be meant with respect to $\mathcal{U}$.\\
    Let $h\in\mathbb{Q}[Y_{1},Y_{1}^{-1},\dots,Y_{n},Y_{n}^{-1},Z_{1},Z_{1}^{-1},\dots,Z_{n},Z_{n}^{-1}]$ be a finite Laurent polynomial with degree $\leq m$, real valued on $\mathbb{T}^{2n}$ with no constant term. Let $C$ be an absolutely irreducible curve over $F$ satisfying the assumption from axiom (4) with respect to $h$. We write $h=\sum_{1\leq i\leq k} c_{i}h_{i}$ and consider the monomials $h_{i}$ for $1\leq i\leq k$ of $h$ separately. Write $h_{i}$ as a product of Laurent monomials $h_{i}^{+}$ and $h_{i}^{\times}$ where $h_{i}^{+}$ is in the variables $Y_{1},Y_{1}^{-1},\dots,Y_{n},Y_{n}^{-1}$ and $h_{i}^{\times}$ is in the variables $Z_{1},Z_{1}^{-1},\dots,Z_{n},Z_{n}^{-1}$. We define $g_{i}$ to be the integral linear transformation corresponding to $h_{i}^{+}$, i.e., such that $\Psi(g_{i}(\bar{x}))=h_{i}^{+}(\Psi^{n}(\bar{x}))$ for $\bar{x}\in F^{n}$.
    Correspondingly, we take $f_{i}:=h_{i}^{\times}$ and observe that $\chi(f_{i}(\bar{x}))=h_{i}^{\times}(\chi^{n}(\bar{x}))$ for $\bar{x}\in F^{n}$.\\
We will write $C_{q}$ for the curve over $\mathbb{F}_{q}$ that corresponds to $C$. As it will then be clear from the context over which field we are working, we will still write $g_{i}$ and $f_{i}$ for their $\mathbb{F}_{q}$-counterparts.
Now, for any $1\leq i\leq k$, we have that $f_{i}$ or $g_{i}$ is not constant on $C^{\prime}(F)$ and hence for almost all $q$ they are not both constant on $C_{q}^{\prime}(\mathbb{F}_{q})$. 
Our aim is to apply Fact \ref{factperelmutterestimates}. We can assume that $\Psi_{q}=\Psi_{q}^{\mathrm{stan}}$ for almost all $q$ as otherwise we can always reduce to this case in the below argument by multiplying with some $a\in\mathbb{F}_{q}$ using Fact \ref{remarkdescriptionalladdcharacters}.
We therefore have to show that for almost all $q$ for any $1\leq i\leq k$ we cannot find $s_{i},r_{i}\in\bar{\mathbb{F}}_{q}(C_{q})$ such that both $g_{i}=s_{i}^{q}-s_{i}$ and $f_{i}=r_{i}^{d_{q}}$ where $d_{q}$ is the order of $\chi_{q}$.\\
First, assume that $g_{i}$ is not constant on $C_{q}$. For almost all $q$ the curve $C_{q}$ is irreducible over $\bar{\mathbb{F}}_{q}$. If $g_{i}(C_{q})$ were finite, then $g_{i}$ would already have to be constant on $C_{q}$ due to irreducibility. It follows that the subfield $L_{i}$ of $\bar{\mathbb{F}}_{q}(C_{q})$ given by $L_{i}:=\bar{\mathbb{F}}_{q}(g_{i})$ fulfills $\text{tr.deg}(L_{i}/\bar{\mathbb{F}}_{q})=1$. We thus have that the degree of the field extension $\bar{\mathbb{F}}_{q}(C_{q})/L_{i}$ is uniformly bounded by some $N\in\mathbb{N}$. Assume that $s_{i}\in\overline{\mathbb{F}}_{q}(C_{q})$ exists such that $g_{i}=s_{i}^{q}-s_{i}$. Consider the field extension $\bar{\mathbb{F}}_{q}(s_{i})/L_{i}$. From $g_{i}=s_{i}^{q}-s_{i}$ we can deduce that $[\bar{\mathbb{F}}_{q}(s_{i}):L_{i}]=q$ since $g_{i}$ is transcendental over $\bar{\mathbb{F}}_{q}$. For $q>N$, this yields a contradiction. 
   Next, we assume that $f_{i}$ is not constant on $C_{q}$. As before, we work over $\bar{\mathbb{F}}_{q}$ for $q$ big and define $L_{i}$ as $\bar{\mathbb{F}}_{q}(f_{i})$ and assume that there is $r_{i}\in \bar{\mathbb{F}}_{q}(C_{q})$ such that $r_{i}^{d_{q}}=f_{i}$. Again we obtain that $[\bar{\mathbb{F}}_{q}(C_{q}):L_{i}]$ is uniformly bounded by some $N\in\mathbb{N}$ contradicting that $\bar{\mathbb{F}}_{q}(r_{i})/L_{i}$ is of degree $d_{q}$ whenever $d_{q}>N$ (which is true for almost all $q$, as $\chi$ is not of finite order).
We finally want to apply Fact \ref{factperelmutterestimates}. Note therefore, as the size of the sets of singular points of $C_{q}$ is bounded\footnote{More precisely, as Fact \ref{factperelmutterestimates} is stated for projective curves, we should say that if we consider the map on rational points given by the restriction of the identification of $C_{q}$ with an open subset of its smooth completion $\Tilde{C_{q}}$, then $|\Tilde{C}_{q}(F_{q})\backslash C_{q}(F_{q})|$ is bounded uniformly in $q$.} by a constant depending on the polynomials defining it, we obtain from Fact \ref{factperelmutterestimates} for sufficiently large $q$ that
    \[\left| \sum_{\bar{x}\in C_{q}^{\prime}(\mathbb{F}_{q})}\Psi_{q}(g_{i}(\bar{x}))\chi_{q}(f_{i}(\bar{x}))\right| \leq t_{i}q^{\frac{1}{2}}\] for some constant $t_{i}>0$. Recall that $c_{i}$ was the coefficient of the monomial $h_{i}$ in $h$. Then, it follows for $s:=\max_{1\leq i\leq k}t_{i}\cdot|c_{i}|$ that
    \[\left|\sum_{\bar{x}\in C_{q}^{\prime}(\mathbb{F}_{q})}h(\Psi_{q}^{n}(\bar{x}))g(\chi_{q}^{n}(\bar{x}))\right|\leq \sum_{1\leq i\leq k}\left|\sum_{\bar{x}\in C_{q}^{\prime}(\mathbb{F}_{q})}\Psi_{q}(g_{i}(\bar{x}))\chi_{q}(f_{i}(\bar{x}))\right|\leq sq^{\frac{1}{2}}.\]
    From that we obtain as a consequence that
    \[\sup\{h(\Psi^{n}(\bar{x}),\chi^{n}(\bar{x}))\,|\,\bar{x}\in C_{q}^{\prime}(\mathbb{F}_{q})\}\geq -sq^{\frac{1}{2}}|C_{q}^{\prime}(\mathbb{F}_{q})|^{-1}.\]
    Finally, by the Lang-Weil bounds \cite{LangWeilBounds}, we have $q^{\frac{1}{2}}|C_{q}^{\prime}(\mathbb{F}_{q})|^{-1}\rightarrow 0$ for $q\rightarrow\infty$ which yields the result.
\end{proof}

\section{Quantifier elimination}\label{sectionqepfplustimeschapterfive}
We recall a definition and fact from \cite{Hrushovski2021AxsTW} which describe definable functions that we can add to $\mathcal{L}_{\mathrm{ring}}$ to obtain quantifier elimination for $\mathrm{PF}$ as well as that substructures are definably closed in the enriched language. 

\begin{definition}\label{definitiondefclosedfunctionsforPF}
We work in $F\models\mathrm{PF}$. Let $P(\bar{X},Y)$ and $Q(\bar{X},Y)$ be integral polynomials. We define $\kappa_{P,Q}(\bar{a})=b$, if there exists at least one $d\in F$ such that $P(\bar{a},d)=0$ and for any such $d$ we have $Q(\bar{a},d)=b$. If no such $b$ exists or if there is no root of $P(\bar{a},X)$ in $F$, we set $\kappa_{P,Q}(\bar{a})=0$.
\end{definition}

%The following result is stated in 3.2 of \cite{Hrushovski2021AxsTW}. We did not find a proof in the literature.

\begin{fact}\label{factqepseudofinitefieldsfunctional}(3.2 of \cite{Hrushovski2021AxsTW}, see Fact 2.1.16 in \cite{ludwig:tel-05236078} for a proof)
    Let $\mathcal{L}_{\mathrm{ring},\kappa}$ be the language $\mathcal{L}_{\mathrm{ring}}$ together with the functions $\kappa_{P,Q}$ from Definition \ref{definitiondefclosedfunctionsforPF}. The theory $\mathrm{PF}$ (together with the interpretations of the $\kappa_{P,Q}$ as above) has elimination of quantifiers, and substructures are definably closed. 
\end{fact}

We now describe the predicate symbols we add to obtain quantifier elimination.

\begin{definition}\label{defintionlanguageforqepfplustimes}
    Let $\bar{a}=(a_{1},\dots,a_{n})$ be an $n$-tuple in $F\models \mathrm{PF}^{+,\times}$. Let $\bar{z}=(z_{1},\dots,z_{m})$, $\bar{y}=(y_{1},\dots, y_{n})$ and $\phi(\bar{z},\bar{y})$ be an algebraic formula in the parameters $\bar{z}$, i.e., $\phi(\bar{z},\bar{a})$ is finite for all $\bar{a}$ in all models of $\mathrm{PF}^{+,\times}$. We write $X_{\bar{a}}^{\phi}$ for the set defined by $\phi(\bar{z},\bar{a})$. Further, let $g:\mathbb{A}^{n}\rightarrow\mathbb{A}^{1}$ be an integral linear transformation and $h:\mathbb{A}^{n}\rightarrow\mathbb{A}^{1}$ an integral multiplicative transformation. %(both can be assumed without constant term). 
    We define the terms $\Theta_{\mathrm{sym}}^{\phi,g,h}(\bar{y})$ as follows:
    \[\Theta_{\mathrm{sym}}^{\phi,g,h}(\bar{a}):=\sum_{\bar{z}\in X_{\bar{a}}^{\phi}}\Psi(g(\bar{z}))\chi(h(\bar{z})).\]
    We define $\mathcal{L}_{\kappa,+,\times}^{\mathrm{sym}}$ to be the language consisting of $\mathcal{L}_{+,\times}$ together with n-ary predicate symbols $\Theta_{\mathrm{sym}}^{\phi,g,h}(\bar{y})$ for all triples $(\phi,g,h)$ as above and symbols for the $\kappa_{P,Q}$ as in Definition \ref{definitiondefclosedfunctionsforPF}. The $\mathcal{L}_{\kappa,+,\times}^{\mathrm{sym}}$-theory $\mathrm{PF}^{+,\times}$ extends the $\mathcal{L}_{+,\times}$-theory $\mathrm{PF}^{+,\times}$ by the interpretations of the new symbols as described above.
\end{definition}

The results in this section are adaptations of the corresponding results for the theory $\mathrm{PF}^{+}$ in \cite{Hrushovski2021AxsTW}. The following Lemmas \ref{lemmaoncoefficientpolynomialoverC}, \ref{lemmaGaloisextensionpfplustimes}, \ref{lemmapassingtoalgclosurepfplustimes} and \ref{lemmaembeddinglemmapfplustimes} follow a similar structure as Lemma 3.8 and Proposition 3.9 in \cite{Hrushovski2021AxsTW}. Our presentation of Lemma \ref{lemmaembeddinglemmapfplustimes} is closer to the one of quantifier elimination for $\mathrm{PF}$ as presented in lecture notes by Chatzidakis \cite{chatzidakislecturenotespsf}. We start with a technical lemma that was already implicitly present in the proof of Lemma 3.8 in \cite{Hrushovski2021AxsTW}.

\begin{lemma}\label{lemmaoncoefficientpolynomialoverC}
Let $G$ be a finite group that acts on a finite set $B$ and $\zeta: B\rightarrow\mathbb{C}^{\times}$ a map. Consider the $G$-action on the ring $\mathbb{C}[T_{b}]_{b\in B}$ that fixes $\mathbb{C}$ and satisfies $g.T_{b}:=T_{g.b}$. Let $f=\sum_{b\in B}\zeta(b)T_{b}$ and $P_{f}=\prod_{g\in G}g.f$. In this case the monomials of $P_{f}$ all have coefficients that are sums of terms of the form $\frac{1}{K}\sum_{g\in G}\zeta(g.b_{1})\cdots \zeta(g.b_{k})$ for some $b_{1},\dots,b_{k}\in B$ and $k=|G|, K\in\mathbb{N}$.    
\end{lemma}

\begin{proof}
We fix a monomial $Q=\alpha\cdot\prod_{b\in\beta}T_{b}$ where $\beta$ is an unordered tuple of size $n=|G|$ consisting of elements from $B$ and $\alpha$ is the coefficient in $P_{f}$. Let $R$ be the set of bijections\footnote{By a bijection $f:G\rightarrow\beta$ between a set $G$ and a multiset $\beta$ we mean by abuse of notation a function $f:G\rightarrow\{b|b\in\beta\}$ from $G$ to the underlying set of $\beta$ such that for all $b$ that occur in $\beta$ we have that $|f^{-1}(b)|$ equals the number of occurrences of $b$ in $\beta$.} $\gamma:G\rightarrow\beta$ and write $b_{\gamma,h}:=\gamma(h)$ as well as $\Tilde{b}_{\gamma,h}:=h^{-1}(b_{\gamma,h})$ for $h\in G$. We then have\[\alpha=\sum_{\gamma\in R}\prod_{h\in G}\zeta(h^{-1}(b_{\gamma,h}))=\sum_{\gamma\in R}\prod_{h\in G}\zeta(\Tilde{b}_{\gamma,h})\]
    The group $G$ acts on $R$ by $g.\gamma(h):=\gamma(hg^{-1})$ for $\gamma\in R$ and $g,h\in G$. Let $\bigcup_{1\leq j\leq k}G\cdot\tau_{j}$ be a partition of $R$ into $G$-orbits (with representatives $\tau_{j}$) and corresponding stabilisers $G_{\tau_{j}}$. It then follows that
    \[\alpha=\sum_{1\leq j\leq k}\sum_{\gamma\in G\cdot\tau_{j}}\prod_{h\in G}\zeta(h^{-1}(b_{\gamma,h}))=\sum_{1\leq j\leq k}\frac{1}{|G_{\tau_{j}}|}\sum_{g\in G}\prod_{h\in G}\zeta(h^{-1}(b_{g.\tau_{j},h})).\]
    Next, we use $b_{g.\tau_{j},h}=b_{g.\tau_{j},hg}=b_{\tau_{j},h}$ (the second equality holds since $g.\tau_{j}(hg)=\tau_{j}(hgg^{-1})$) to obtain  
    % will use the following equality
    % \[\prod_{h\in G}\zeta(h^{-1}(b_{g.\tau_{j},h}))=\prod_{h\in G}\zeta(g^{-1}h^{-1}(b_{g.\tau_{j},hg}))=\prod_{h\in G}\zeta(g^{-1}h^{-1}(b_{\tau_{j},h})).\]
    % The first equality is trivial and the second holds since $g.\tau_{j}(hg)=\tau_{j}(hgg^{-1})$. Combining this with the above, we obtain 
    \[\sum_{g\in G}\prod_{h\in G}\zeta(h^{-1}(b_{g.\tau_{j},h}))=\sum_{g\in G}\prod_{h\in G}\zeta(g^{-1}(h^{-1}(b_{\tau_{j},h})))=\sum_{g\in G}\prod_{h\in G}\zeta(g^{-1}(\Tilde{b}_{\tau_{j},h})).\]
\end{proof}

\begin{lemma}\label{lemmaGaloisextensionpfplustimes}
Let $B$ be a finite Galois extension of a field $A$ where $A\subseteq B\subseteq F\models \mathrm{PF}$. Moreover, assume that $B$ carries two additive characters $\Psi$ and $\Psi^{\prime}$ and multiplicative characters $\chi$ and $\chi^{\prime}$ that agree on $A$ such that moreover the terms $\Theta_{\mathrm{sym}}^{\phi,g,h}(\bar{y})$ are the same for all $\bar{a}\in A^{|\bar{a}|}$ whether calculated with $\Psi$ and $\chi$ or with $\Psi^{\prime}$ and $\chi^{\prime}$.
Then there is always some $\tau\in \mathrm{Gal}(B/A)$ such that $\Psi(b)=\Psi^{\prime}(\tau(b))$ and $\chi(b)=\chi^{\prime}(\tau(b))$ for all $b\in B$.
\end{lemma}
\begin{proof} 
As $G:=\mathrm{Gal}(B/A)$ is finite, it suffices to find for any finite set $B_{0}$ of $B$ some element $\tau\in G$ such that $\Theta(b)=\chi(\tau(b))$ for all $b\in B_{0}$. We can enlarge $B_{0}$ to be $G$-invariant and so that $B=A(B_{0})$ as fields. We further assume that $1\in B_{0}$.
        We define variables $S_{b},T_{b}$ for any $b\in B_{0}$ and define a $G$-action on $\mathbb{C}[S_{b},T_{b}]_{b\in B_{0}}$ that fixes $\mathbb{C}$ by setting $g.S_{b}=S_{g(b)}$ and $g.T_{b}=T_{g(b)}$ for $g\in G$. We consider the linear polynomials $f=\sum_{b\in B_{0}}(\Psi(b)S_{b}+\chi(b)T_{b})$ and $f^{\prime}=\sum_{b\in B_{0}}(\Psi^{\prime}(b)S_{b}+\chi^{\prime}(b)T_{b})$. We work with the polynomials $P_{f}=\prod_{g\in G}g.f$ and $P_{f^{\prime}}=\prod_{g\in G}g.f^{\prime}$. As $\mathbb{C}[S_{b},T_{b}]_{b\in B_{0}}$ is a unique factorisation domain it suffices to show that $P_{f}=P_{f^{\prime}}$ to obtain that $f$ divides $P_{f^{\prime}}$ and thus $f=g.f^{\prime}$ for some $g\in G$. Now we consider the monomials of $P_{f}$. From Lemma \ref{lemmaoncoefficientpolynomialoverC} we deduce that their coefficients are given as sums of terms of the form $\frac{1}{K}\sum_{g\in G}\Psi(g(b_{j_{1}}))\cdots \Psi(g(b_{j_{l}}))\chi(g(b_{i_{1}}))\cdots \chi(g(b_{i_{k}}))$ where $l+k=|G|$ and $b_{j_{1}},\dots,b_{j_{l}},b_{i_{1}},\dots,b_{i_{k}}\in B_{0}$ (and similarly for $P_{f^{\prime}}$) as well as $K\in\mathbb{N}$.\\
    The set $X:=\{(g.b_{j_{1}},\dots,g.b_{j_{l}},g.b_{i_{1}},\dots,g.b_{i_{k}})\,|\,g\in G\}$ is clearly finite and $G$-invariant hence $Aut(F/A)$-invariant and thus definable over $A$ in $\mathcal{L}_{\mathrm{ring},\kappa}$ (and consequently $\mathcal{L}_{\mathrm{ring}}$-definable) by a formula $\phi(\bar{x},\bar{a})$ for $\bar{a}$ a tuple in $A$. Consider $\alpha:\mathbb{A}^{l+k}\rightarrow \mathbb{A}^{1}$ defined by $\alpha(\bar{z}):=z_{1}+\dots +z_{l}$ and $\beta:\mathbb{A}^{l+k}\rightarrow \mathbb{A}^{1}$ defined by $\beta(\bar{z}):=\prod_{1\leq i\leq k}z_{l+i}$. Then we obtain
    \[\sum_{g\in G}\Psi(g(b_{j_{1}}))\cdots \Psi(g(b_{j_{l}}))\chi(g(b_{i_{1}}))\cdots \chi(g(b_{i_{k}}))=\sum_{\bar{x}\in X_{\bar{a}}^{\phi}}\Psi(\alpha(\bar{x}))\chi(\beta(\bar{x}))\]
    and the latter is by definition $\Theta_{\mathrm{sym}}^{\phi,\alpha,\beta}(\bar{a})$. Since this works for all monomials occuring in $P_{f}$ we obtain that $P_{f}=P_{f^{\prime}}$ which completes the proof.
\end{proof}

\begin{lemma}\label{lemmapassingtoalgclosurepfplustimes}

Let $A$ and $A^{\prime}$ be isomorphic substructures of $F\models \mathrm{PF}^{+,\times}$ in the language $\mathcal{L}_{\kappa,+,\times}^{\mathrm{sym}}$. Then the isomorphism $\alpha$ between $A$ and $A^{\prime}$ extends to an isomorphism between the relative algebraic closures of $A$ and $A^{\prime}$.
    
\end{lemma}

\begin{proof}
Since we have quantifier elimination for the theory $\mathrm{PF}$ in the language $\mathcal{L}_{\mathrm{ring},\kappa}$ by Fact \ref{factqepseudofinitefieldsfunctional}, we obtain a field isomorphism that extends $\alpha$ to the relative algebraic closures of $A$ and $A^{\prime}$. Also note that an $\mathcal{L}_{+,\times}$-structure isomorphism between algebraically closed structures is already an $\mathcal{L}_{\kappa,+,\times}^{\mathrm{sym}}$-isomorphism. We can thus reduce to the situation of working with $A$ and its relative algebraic closure $\Tilde{A}$ carrying two additive characters $\Psi, \Psi^{\prime}$ and two multiplicative characters $\chi,\chi^{\prime}$, that extend the one on $A$ and we have to find an $\mathcal{L}_{\mathrm{ring},\kappa}$-structure automorphism $\tau$ of $\Tilde{A}$ over $A$ such that for all $c\in\Tilde{A}$ we have $\Psi(c)=\Psi^{\prime}(\tau(c))$ and $\chi(c)=\chi^{\prime}(\tau(c))$.\\
By compactness we can reduce to finitely generated algebraic extensions of $A$. However, to apply Lemma \ref{lemmaGaloisextensionpfplustimes} we have to make sure that we can assume these extensions to be Galois. But this follows as by Fact \ref{factqepseudofinitefieldsfunctional} $\mathrm{dcl}_{\mathrm{ring},\kappa}(A)=A$, so in particular $\mathrm{dcl}_{\mathrm{ring}}(A)=A$, and thus $\mathrm{Gal}(\Tilde{A}/A)=\mathrm{Aut}_{\mathrm{ring}}(\Tilde{A}/A)$.
\end{proof}

% We now turn to prove an embedding lemma in strong analogy to Lemma \ref{lemmaembeddinglemma}.

\begin{lemma}\label{lemmaembeddinglemmapfplustimes}
Let $F$ be an $\aleph_{1}$-saturated model of $\mathrm{PF}^{+,\times}$ with Galois generator $\sigma\in \mathrm{Gal}(F^{\mathrm{alg}}/F)$, $K\subseteq F$ a countable substructure, relatively algebraically closed as a field in $F$. Let $E\supseteq K$ be an extension of $\mathcal{L}_{+}$-structures such that the following hold.
\begin{itemize}
    \item $E$ is a field of transcendence degree $1$ over $K$.
    \item $K$ is relatively algebraically closed in $E$.
    \item $\Psi:(E,+)\rightarrow S^{1}$ is a group homomorphism.
    \item $\chi:(E^{\times},\cdot)\rightarrow S^{1}$ is a group homomorphism (and $\chi(0)=0$).
    \item There is a topological generator $\sigma^{\prime}$ of $\mathrm{Gal}(E^{\mathrm{alg}}/E)$ such that $\sigma^{\prime}\mapsto\sigma$ induces an isomorphism $\mathrm{Gal}(E^{\mathrm{alg}}/E)\rightarrow \mathrm{Gal}(F^{\mathrm{alg}}/F)$.
\end{itemize}
Then, we find a field embedding $\tau$ of $E^{\mathrm{alg}}$ over $K^{\mathrm{alg}}$ into $F^{\mathrm{alg}}$ such that $\tau(\sigma^{\prime}(a))=\sigma(\tau(a))$ for all $a\in E^{\mathrm{alg}}$, $\tau(E)$ is relatively algebraically closed in $F$ and $\tau$ preserves $\Psi$ and $\chi$.

% Let $F$ be an $\aleph_{1}$-saturated model of $\mathrm{PF}^{+,\times}$, $K\subseteq F$ a relatively algebraically closed countable substructure and $E\supseteq K$ be a relatively algebraically closed substructure of some $\aleph_{1}$-saturated model $F^{\prime}\models \mathrm{PF}^{+,\times}$ such that $E$ is of transcendence degree $1$ over $K$ and $E$ is countable. Further, let $\sigma^{\prime}$ be a generator of $\mathrm{Gal}(E^{\mathrm{alg}}/E)$ and $\sigma$ of $ \mathrm{Gal}(F^{\mathrm{alg}}/F)$. We assume that $\sigma^{\prime}\rightarrow \sigma$ induces an isomorphism $\mathrm{Gal}(E^{\mathrm{alg}}/E)\rightarrow \mathrm{Gal}(F^{\mathrm{alg}}/F)$.\\
% Then, we find a field embedding $\tau$ of $E^{\mathrm{alg}}$ over $K^{\mathrm{alg}}$ into $F^{\mathrm{alg}}$ such that $\tau(\sigma^{\prime}(a))=\sigma(\tau(a))$ for all $a\in E^{\mathrm{alg}}$, such that $\tau(E)$ is relatively algebraically closed in $F$ and such that $\tau$ preserves $\Psi$ and $\chi$.

\end{lemma}
\begin{proof}
We work in a big algebraically closed field which we assume to contain all the fields mentioned above and replace $E$ by some $K$-isomorphic copy, in order to assume that $E$ is linearly disjoint from $F$ over $K$. We take $\Tilde{\sigma}$ as an extension of $(\sigma^{\prime},\sigma)\in \mathrm{Gal}(E^{\mathrm{alg}}/E)\times_{\mathrm{Gal}(K^{\mathrm{alg}}/K)}\mathrm{Gal}(F^{\mathrm{alg}}/F) \cong \mathrm{Gal}(E^{\mathrm{alg}}F^{\mathrm{alg}}/EF)$ to $(EF)^{\mathrm{alg}}$ and consider the field $M:=\mathrm{Fix}(\Tilde{\sigma})$.
We choose some extension of $\Psi$ to a group homomorphism $(M,+)\rightarrow (\mathbb{T},\cdot)$ and similarly for $\chi$ and $(M^{\times},\cdot)$ both arbitrarily chosen outside of the additive (resp. multiplicative) span of $E$ and $F$. Thus, we extend $M$ to an $\mathcal{L}_{+,\times}$-structure.\\
The fields $E$ and $F$ are by construction relatively algebraically closed in $M$ and moreover $\sigma^{\prime}\mapsto\Tilde{\sigma}\mapsto\sigma$ induces isomorphisms of the corresponding Galois groups. Hence we obtain that $M^{\mathrm{alg}}=ME^{\mathrm{alg}}=MF^{\mathrm{alg}}$ and the latter is the ring generated by $M$ and $F^{\mathrm{alg}}$.
As $E$ is countable, there are some countable $M_{0}\subset M$ as well as $F_{0}\subset F$ such that $F_{0}^{\mathrm{alg}}[M_{0}]$ contains $E^{\mathrm{alg}}$. We may assume that $F_{0}$ is relatively algebraically closed in $F$, simply by passing to its relative algebraic closure. It follows that $F_{0}(M_{0})$ is an extension of $F_{0}$ of transcendence degree $1$ with $F_{0}$ being relatively algebraically closed in $F_{0}(M_{0})$ (as $F$ already was in $F(M_{0})$).\\
We will now show that it then suffices to find an $F_{0}$-morphism $F_{0}(M_{0})\rightarrow F$ which preserves the values of $\Psi$ and $\chi$. (Note that it then has to be injective as $\chi(x)=0\iff x=0$). Assume for the moment that such a morphism exists. Then, it extends to an $F_{0}^{\mathrm{alg}}$-morphism $\tau:F_{0}^{\mathrm{alg}}[M_{0}]\rightarrow F^{\mathrm{alg}}$ and this is indeed the $\tau$ we are after because of the following. For any $a\in E^{\mathrm{alg}}$ we can write $a=\sum m_{i}b_{i}$ for some $m_{i}\in M_{0}$ and $b_{i}\in F_{0}^{\mathrm{alg}}$ and then
\[\tau(\sigma^{\prime}(a))=\tau(\Tilde{\sigma}(a))=\sum \tau(m_{i})\Tilde{\sigma}(b_{i})=\sum \tau(m_{i})\sigma(b_{i})=\sigma(\tau(a)).\]
Moreover, $a$ is fixed by $\sigma^{\prime}$ if and only if $\tau(a)$ is fixed by $\sigma$. Thus, $\tau(E)$ is relatively algebraically closed in $F_{0}$. We now construct this $F_{0}$-morphism $F_{0}(M_{0})\rightarrow F$.\\
By compactness it suffices to consider finite sets $\{m_{1},\dots,m_{n}\}\subseteq M_{0}$. Let $C\subset \mathbb{A}^{n}$ be the variety with generic point $(m_{1},\dots,m_{n})$, then we may assume that $C$ is a curve. Thus, we have to show that we find some $\bar{c}\in C(F)\backslash C(F_{0})$ such that moreover $\Psi(m_{i})=\Psi(c_{i})$ and $\chi(m_{i})=\chi(c_{i})$ for all $1\leq i\leq n$. Here we can assume that $C$ is not zero-degenerate by passing to a suitable projection if necessary.\\
By compactness it suffices to show for any finite subset $\{\bar{\gamma}_{1},\dots,\bar{\gamma}_{t}\}$ of $C(F_{0})$ that we can find some $\bar{c}$ as above in $C(F)\backslash \{\bar{\gamma}_{1},\dots,\bar{\gamma}_{t}\}$. Therefore, we work in $\mathbb{A}^{2n}$ and let the curve $\Tilde{C}$ to be defined by the polynomials that define $C$ together with all the polynomials of the form $X_{n+i}\gamma_{k,i}=1$ where $\gamma_{k,i}$ is the $i$-th coordinate of $\gamma_{k}$. (Note that there is no harm in assuming $\gamma_{k,i}\neq 0$ as our original curve was not zero-degenerate.) We then choose an (absolutely) irreducible component of $\Tilde{C}$ which we again name $C$ and have thus reduced the problem to find $\bar{c}\in C(F)$ such that $\Psi(m_{i})=\Psi(c_{i})$ and $\chi(m_{i})=\chi(c_{i})$ for all $1\leq i\leq n$ and $\bar{m}$ a generic point of $C$.\\
Let $\{d_{1},\dots,d_{k}\}$ be a free generating set of the additive group generated by $\{m_{1},\dots,m_{n}\}$ over $F_{0}$ and $\alpha:\mathbb{A}^{n}\rightarrow\mathbb{A}^{k}$ be the integral linear transformation such that $\alpha(\bar{m})=\bar{d}$. Further, let $\{e_{1},\dots,e_{l}\}$ be a free generating set of the multiplicative group generated by $\{m_{1},\dots,m_{n}\}$ over $F_{0}^{\times}$ and $\beta:\mathbb{A}^{n}\rightarrow\mathbb{A}^{l}$ be the integral multiplicative transformation such that $\beta(\bar{m})=\bar{e}$.\\
Now we use Lemma \ref{lemmapfplustimesequivofdefinition} and obtain by property $(\star)$ from Definition \ref{definitionporpertyequivtoaxiompfplustimes} applied on the datum given by $C$, $\alpha$ and $\beta$ that $(\Psi^{k}\times\chi^{l})(\Delta(F))$ is dense in $\mathbb{T}^{k+l}$ where $\Delta$ is the diagonal defined as in Definition \ref{definitionporpertyequivtoaxiompfplustimes}. But this means that we find by saturation $(f_{1},\dots,f_{n})\in C(F)$ with $\Psi^{k}(\alpha(\bar{f}))=\Psi^{k}(\alpha(\bar{m}))$ and $\chi^{l}(\beta(\bar{f}))=\chi^{l}(\beta(\bar{m}))$. But then it already follows that $\Psi^{n}(\bar{f})=\Psi^{n}(\bar{m})$ and $\chi^{n}(\bar{f})=\chi^{n}(\bar{m})$ as the $m_{i}$ for $1\leq i\leq n$ lie in the additive group generated by $d_{1},\dots,d_{k}$ and in the multiplicative group generated by $e_{1},\dots,e_{l}$ which completes the proof.
\end{proof}

We need the following two basic facts that already appear in \cite{Hrushovski2021AxsTW} in order to complete the proof of quantifier elimination.

\begin{fact}\label{facttowerfullsubfieldsyieldsfullsubfield}(e.g., see proof of Lemma 3.10 in \cite{Hrushovski2021AxsTW})
    Let $F$ be a pseudofinite field and $K\subset F$ a subfield, such that the restriction map $\mathrm{Gal}(F^{\mathrm{alg}}/F)\rightarrow \mathrm{Gal}(K^{\mathrm{alg}}/K)$ is an isomorphism. If $K\subset E\subset F$ is a tower of fields and $E$ is relatively algebraically closed in $F$, then the restriction map $\mathrm{Gal}(F^{\mathrm{alg}}/F)\rightarrow \mathrm{Gal}(E^{\mathrm{alg}}/E)$ is also an isomorphism.
\end{fact}

\begin{fact}(Lemma 3.11 in \cite{Hrushovski2021AxsTW}.)\label{factexistenceoffullextensionoverQ}
Let $T$ be a completion of $\mathrm{PF}$. Then there exists $F\models T$ and a subfield $F_{1}$ of $F$ of transcendence degree $1$ over $\mathbb{Q}$ such that the restriction $\mathrm{Gal}(F^{\mathrm{alg}}/F)\rightarrow \mathrm{Gal}(F_{1}^{\mathrm{alg}}/F_{1})$ is an isomorphism.
    
\end{fact}

\begin{theorem}\label{qepfplustimes}
    The theory $\mathrm{PF}^{+,\times}$ has quantifier elimination in the language $\mathcal{L}_{\kappa,+,\times}^{\mathrm{sym}}$.
\end{theorem}

\begin{proof}
As $\mathcal{L}_{\kappa,+,\times}^{\mathrm{sym}}$ is countable, it suffices to show that we can extend every isomorphism between countable substructures $K_{1}$ and $K_{2}$ to a back-and-forth system between $\aleph_{1}$-saturated models $F_{1}$ and $F_{2}$. 
Using Lemma \ref{lemmapassingtoalgclosurepfplustimes}, we can assume that $F_{1}$ and $F_{2}$ share a common countable substructure $K$ which is a relatively algebraically closed subfield of $F_{1}$ and $F_{2}$. It then suffices to show that every $E_{1}\subseteq F_{1}$, which is relatively algebraically closed of transcendence degree $1$ over $K$, can be embedded into $F_{2}$ so that the image is relatively algebraically closed again.\\
We now simply choose for all our fields at hand topological generators of their Galois group. Given those choices, we obtain a map $\mathrm{Gal}(F_{2}^{\mathrm{alg}}/F_{2})\rightarrow \mathrm{Gal}(E_{1}^{\mathrm{alg}}/E_{1})$ obtained by identifying the Galois generators. If this map is an isomorphism, we can apply Lemma \ref{lemmaembeddinglemmapfplustimes} and are done.
Now to obtain this we can apply Fact $\ref{facttowerfullsubfieldsyieldsfullsubfield}$ from which it follows that we only have to make sure that $\mathrm{Gal}(F_{1}^{\mathrm{alg}}/F_{1})\rightarrow \mathrm{Gal}(K^{\mathrm{alg}}/K)$ and  $\mathrm{Gal}(F_{2}^{\mathrm{alg}}/F_{2})\rightarrow \mathrm{Gal}(K^{\mathrm{alg}}/K)$ are both isomorphisms. So, it suffices to show that there is always an $\mathcal{L}_{+,\times}$-isomorphism fixing $K$ between two relatively algebraically closed extensions $K_{1}/K$,  $K_{2}/K$ in $F_{1}$,$F_{2}$, respectively, such that $\mathrm{Gal}(F_{1}^{\mathrm{alg}}/F_{1})\rightarrow \mathrm{Gal}(K_{1}^{\mathrm{alg}}/K_{1})$ and  $\mathrm{Gal}(F_{2}^{\mathrm{alg}}/F_{2})\rightarrow \mathrm{Gal}(K_{2}^{\mathrm{alg}}/K_{2})$ are isomorphisms.\\
From Fact \ref{factexistenceoffullextensionoverQ} and saturation it follows that there is such an extension $H_{1}$ of transcendence degree $1$ over $\mathbb{Q}$ in $F_{1}$. Let $\Tilde{K}_{1}$ be the relative algebraic closure of $H_{1}K$ in $F_{1}$. It has transcendence degree at most $1$ over $K$ and $\mathrm{Gal}(F_{1}^{\mathrm{alg}}/F_{1})\rightarrow \mathrm{Gal}(\Tilde{K}_{1}^{\mathrm{alg}}/\Tilde{K}_{1})$ is an isomorphism. Let $b_{i}$ be an $i$-tuple in $\Tilde{K}_{1}$ consisting of the coefficients of an irreducible polynomial of degree $i$ over $\Tilde{K}_{1}$ for any $i\in\mathbb{N}$. Further, let $K_{1}$ be the relative algebraic closure of $\cup_{i\in \mathbb{N}}Kb_{i}$. Now $K_{1}$ is countable and $\mathrm{Gal}(F_{1}^{\mathrm{alg}}/F_{1})\rightarrow \mathrm{Gal}(K_{1}^{\mathrm{alg}}/K_{1})$ is still an isomorphism. We reproduce this construction obtain $K_{2}\subseteq F_{2}$.\\
Finally, let $M$ be an amalgam of $K_{1}$ and $K_{2}$ over $K$ constructed in the same way as in in the proof of Lemma \ref{lemmaembeddinglemmapfplustimes}. We can then, by Lemma \ref{lemmaembeddinglemmapfplustimes}, embed $M$ into $F_{1}$ over $K_{1}$, and similarly, $M$ into $F_{2}$ over $K_{2}$ such that (in both cases) the image $K_{1}^{\prime}$ (resp. $K_{2}^{\prime}$) is relatively algebraically closed. But now $K_{1}^{\prime}$ and $K_{2}^{\prime}$ are $K$-isomorphic, which completes the proof.
\end{proof}

\begin{corollary}\label{corollaryeltequivpfplustimes}

Let $(F_{1},\Psi_{1},\chi_{1}),(F_{2},\Psi_{2},\chi_{2})$ be two models of $\mathrm{PF}^{+,\times}$
 and $E$ a common relatively algebraically closed subfield on which $\Psi_{1}$ and $\Psi_{2}$ as well as $\chi_{1}$ and $\chi_{2}$ coincide, then
 \[(F_{1},\Psi_{1},\chi_{1})\equiv_{E} (F_{2},\Psi_{2},\chi_{2}).
 \]
  In particular, the completions of $\mathrm{PF}^{+,\times}$ are determined by the isomorphism types of $F\cap\bar{\mathbb{Q}}$ together with the corresponding restrictions of $\Psi$ and $\chi$ for $(F,\Psi,\chi)\models\mathrm{PF^{+,\times}}$.
 \end{corollary}
\begin{proof}
    Since $E$ is algebraically closed, it carries the same $\mathcal{L}_{\kappa,+,\times}^{\mathrm{sym}}$-structure if computed from $(\Psi_{1},\chi_{1})$ or $(\Psi_{2},\chi_{2})$.
\end{proof}

 \begin{corollary}\label{corollarymodelcompleteness}
   Let $\mathcal{L}_{+,\times,\bar{c}}$ consist of $\mathcal{L}_{+,\times}$ with additional constant tuples $\bar{c}_{i}$ of length $i$ for every $i\in\mathbb{N}_{\geq 2}$.  The $\mathcal{L}_{+,\times,\bar{c}}$-theory $\mathrm{PF}^{+,\times}$ obtained by adding sentences stating that the polynomials with parameter tuple $\bar{c}_{i}$ are irreducible is model complete.
 \end{corollary}
 \begin{proof}
     Direct from Corollary \ref{corollaryeltequivpfplustimes}.
 \end{proof}

\begin{corollary}\label{corollarymodelcompletion}
The theory $\mathrm{PF}^{+,\times}$ is the model completion of the theory of structures $(F_{0},\Psi_{0},\chi_{0},\bar{c})$ where $F_{0}$ is a field of
characteristic $0$  with $\mathrm{Gal}(F_{0})=\hat{\mathbb{Z}}$ and $\Psi_{0}$ (resp. $\chi_{0}$) are additive (resp. multiplicative) characters and $\bar{c}$ a tuple of constants interpreted as in Corollary \ref{corollarymodelcompleteness}.
\end{corollary}
\begin{proof}
We can extend every such structure to model of $\mathrm{PF}^{+,\times}$ by a chain-construction using Lemma \ref{lemmapfplustimesequivofdefinition}. Then the results already follows from Corollary \ref{corollarymodelcompleteness}. 
\end{proof}

\textbf{Basic predicates.} We will now show that we can uniformly approximate all quantifier-free $\mathcal{L}_{\kappa,+,\times}^{\mathrm{sym}}$-definable predicates (and thus all definable predicates) by what we call basic predicates. This adapts the corresponding result for $\mathrm{PF}^{+}$ in \cite{Hrushovski2021AxsTW} to our setting.

\begin{definition}\label{definitionbasicpredicatesPFplustimes}
    We define a \textit{basic predicate} in $n$-variables $\bar{x}=(x_{1},\dots,x_{n})$ as follows.
    \begin{enumerate}[(1)]
        
 \item A predicate of the form \footnote{Again, if we work in a completion of $\mathrm{PF}^{+,\times}$ with $\chi(\mathbb{Q})$ dense in $S^{1}$ we do not need a sum but it suffices to work with predicates of the form $\lambda\cdot\Theta_{\mathrm{sym}}^{\phi,g,h}(\kappa_{1}(\bar{x}),\dots\kappa_{m}(\bar{x}))$.}
 \[\lambda\cdot\Theta_{\mathrm{sym}}^{\phi,g,h}(\kappa_{1}(\bar{x}),\dots\kappa_{m}(\bar{x}))+i\cdot \lambda^{\prime}\cdot\Theta_{\mathrm{sym}}^{\phi^{\prime},g^{\prime},h^{\prime}}(\kappa^{\prime}_{1}(\bar{x}),\dots\kappa_{m^{\prime}}^{\prime}(\bar{x}))\]
    where the $\kappa_{1}(\bar{x}),\dots\kappa_{m}(\bar{x}),\kappa_{1}^{\prime}(\bar{x}),\dots\kappa_{m^{\prime}}^{\prime}(\bar{x})$ are $\mathcal{L}_{\mathrm{ring},\kappa}$-terms in the variables $\bar{x}=(x_{1},\dots,x_{n})$, $\lambda,\lambda^{\prime}\in\mathbb{Q}$ and $\Theta_{\mathrm{sym}}^{\phi,g,h}$ is as in Definition \ref{defintionlanguageforqepfplustimes}.
    \item A predicate of the form
    \[\sum_{1\leq i\leq k}\beta_{i}(\bar{x})\cdot\mathbbm{1}_{C_{i}}\] where the $\beta_{i}$ are as in case (1) and $C_{1},\dots,C_{k}$ is a disjoint partition of $F^{|\bar{x}|}$ into $\mathcal{L}_{\mathrm{ring}}$-definable sets.
\end{enumerate}
\end{definition}

% \begin{definition}\label{definitionbasicpredicatesPFplustimes}
%     We define a \textit{basic predicate} in $n$-variables $\bar{x}=(x_{1},\dots,x_{n})$ to be one of the following two: 
%     \begin{enumerate}[(1)]
%         \item An $\mathcal{L}_{\mathrm{ring},\kappa}$-predicate $\phi(\bar{x})_{(c_{1},c_{2})}$ where $\phi(\bar{a})=c_{1}\in\mathbb{Q}(i)^{\mathrm{alg}}\subset\mathbb{C}$ if $\phi(\bar{a})$ holds for the tuple $\bar{a}$ in $F\models \mathrm{PF}$ and otherwise $\phi(\bar{a})=c_{2}\in\mathbb{Q}(i)^{\mathrm{alg}}\subset\mathbb{C}$.
%         \item A predicate of the form  
%     \[\lambda\cdot\Theta_{\mathrm{sym}}^{\phi,g,h}(\kappa_{1}(\bar{x}),\dots\kappa_{m}(\bar{x}))\]
%     where the $\kappa_{1}(\bar{x}),\dots\kappa_{m}(\bar{x})$ are $\mathcal{L}_{\mathrm{ring},\kappa}$-terms for $\bar{x}=(x_{1},\dots,x_{n})$ and $\lambda\in\mathbb{Q}(i)^{\mathrm{alg}}\subset\mathbb{C}$.
%     \item A predicate of the form
%     \[\sum_{1\leq i\leq k}\beta_{i}(\bar{x})\cdot\mathbbm{1}_{C_{i}}\] where the $\beta_{i}$ are as in case (1) or (2) and $C_{1},\dots,C_{k}$ is a disjoint partition of $F\models \mathrm{PF}$ into $\mathcal{L}_{\mathrm{ring}}$-definable sets.
%     \end{enumerate}
% \end{definition}

\begin{lemma}\label{lemmaalgebraofbasicpredicatespfplustimes}
    Every $\mathcal{L}_{\kappa,+,\times}^{\mathrm{sym}}$-definable predicate in $n$-variables can be uniformly approximated by basic predicates in $n$-variables.
\end{lemma}

\begin{proof}
 Let $S_{n}$ denote the space of $n$-types equipped with the logic topology (Definition 8.4 of \cite{mtfms}). The space of definable predicates in $n$ variables can be identified with $C(S_{n},\mathbb{C})$, the $C^{*}$-algebra of complex-valued continuous functions with respect to the topology of uniform convergence on $S_{n}$ (see Proposition 8.10 in \cite{mtfms}). By the complex Stone-Weierstrass-theorem it then suffices to show that the uniform limits of basic formulas form a unital $C^{*}$-algebra that separates points in $S_{n}$. Seperation of points, i.e., that for every two distinct complete $n$-types $p_{1},p_{2}$ there is a basic predicate that takes different values for the realisations of $p_{1}$ and $p_{2}$, already follows from quantifier elimination of $\mathrm{PF}^{+,\times}$ in $\mathcal{L}_{\kappa,+,\times}^{\mathrm{sym}}$. To show that it is a unital $C^{*}$-algebra, it will suffice to show that the basic formulas are closed under (pointwise) addition, multiplication and complex conjugation as well as scalar multiplication. As unions, intersections and complements of $\mathcal{L}_{\mathrm{ring}}$-definable sets are again $\mathcal{L}_{\mathrm{ring}}$-definable the only thing to take care of are sums, products and complex conjugates of basic predicates as in case (1) of Definition \ref{definitionbasicpredicatesPFplustimes}.\\
We start by showing closure under multiplication. Consider two basic predicates $\Theta_{\mathrm{sym}}^{\phi,g,h}(\bar{x})$ and $\Theta_{\mathrm{sym}}^{\Tilde{\phi},\Tilde{g},\Tilde{h}}(\bar{x})$ as in Definition \ref{definitionbasicpredicatesPFplustimes}. Let $\gamma(\bar{z},\bar{x})$ be the $\mathcal{L}_{\mathrm{ring}}$-formula defining the cartesian product of the definable sets defined by $\phi$ and $\Tilde{\phi}$ with parameters $\bar{x}$, i.e., $\gamma(\bar{z},\bar{x})$ defines $X_{\bar{x}}^{\phi}\times X_{\bar{x}}^{\Tilde{\phi}}=:X_{\bar{x}}^{\gamma}$. Write $\bar{z}_{1}$ for the subtuple of $\bar{z}$ that corresponds to the variables of $X_{\bar{x}}^{\phi}$ and correspondingly write $\bar{z}_{2}$ for the one corresponding to $X_{\bar{x}}^{\Tilde{\phi}}$. Let $\alpha:\mathbb{A}^{|\bar{z}|}\rightarrow\mathbb{A}^{1}$ be the rational function defined by $\alpha(\bar{z}):=g(\bar{z}_{1})+\Tilde{g}(\bar{z}_{2})$ and  $\beta:\mathbb{A}^{|\bar{z}|}\rightarrow\mathbb{A}^{1}$ be defined by $\beta(\bar{z}):=h(\bar{z}_{1})\cdot\Tilde{h}(\bar{z}_{2})$. It then follows that
\[\Theta_{\mathrm{sym}}^{\phi,g,h}(\bar{x})\cdot\Theta_{\mathrm{sym}}^{\Tilde{\phi},\Tilde{g},\Tilde{h}}(\bar{x})=\Theta_{\mathrm{sym}}^{\gamma,\alpha,\beta}(\bar{x}) \] which resolves the case of closure under products.\\
Now, $\bar{q}_{1},\bar{q}_{2}$ will denote tuples of elements of $\mathbb{Q}$ and $r_{1},r_{2},s_{1},s_{2}\in\mathbb{Q}$ which we will specify later on. Let $\gamma(\bar{z},v,w,\bar{x})$ define (in the variables $(\bar{z},v,w)$) the following set (which depends on the choice of $\bar{q}_{1},\bar{q}_{2},r_{1},r_{2},s_{1},s_{2}$): \[X_{\bar{x}}^{\gamma}:=\left(X_{\bar{x}}^{\phi}\times\{(\bar{q}_{2},r_{2},s_{2})\}\right)\cup\left(\{\bar{q}_{1}\}\times X_{\bar{x}}^{\Tilde{\phi}}\times \{(r_{1},s_{1})\}\right).\] 
We set $\alpha:\mathbb{A}^{|\bar{z}|+2}\rightarrow\mathbb{A}^{1}$ to be $\alpha(\bar{z},v,w):=g(\bar{z}_{1})+\Tilde{g}(\bar{z}_{2})+v$ and $\beta:\mathbb{A}^{|\bar{z}|+2}\rightarrow\mathbb{A}^{1}$ to be $\beta(\bar{z},v,w):=h(\bar{z}_{1})\cdot\Tilde{h}(\bar{z}_{2})\cdot w$. Then, if we choose $\bar{q}_{1},\bar{q}_{2},r_{1},r_{2},s_{1},s_{2}$ such that $\Tilde{g}(\bar{q}_{2})=-r_{2}$ and $g(\bar{q}_{1})=-r_{1}$ as well as $\Tilde{h}(\bar{q}_{2})=s_{2}^{-1}$, $h(\bar{q}_{1})=s_{1}^{-1}$ and $s_{1}\neq s_{2}$, then we obtain that
\[\Theta_{\mathrm{sym}}^{\phi,g,h}(\bar{x})+\Theta_{\mathrm{sym}}^{\Tilde{\phi},\Tilde{g},\Tilde{h}}(\bar{x})=\Theta_{\mathrm{sym}}^{\gamma,\alpha,\beta}(\bar{x}). \]
It is now be straightforward to see that closure under addition follows for all basic predicates using that the allowed coefficients are in $\mathbb{Q}$. Finally, to obtain the complex conjugate of $\Theta_{\mathrm{sym}}^{\phi,g,h}(\bar{x})$ it suffices to replace $g$ by $-g$ and $h$ by $h^{-1}$.
% Now if the above terms come with scalar multipliers we can argue as in the proof of Lemma \ref{lemmabasicformulasaredensepftimes} once we use that for a (discrete) first-order formula $\phi(\bar{x},\bar{y})$ that defines a finite set (in $\bar{x}$) for all $\bar{y}$ we can partition the parameter space $F^{|\bar{y}|}$ into definable sets such that for all $\bar{y}$ in a given set of the partition the definable sets (in $\bar{x}$) have the same size.\footnote{This then allows us to do definably a case distinction for the definable functions $\kappa_{i}(\bar{x})$.}
\end{proof}

\section{Limit theory}\label{sectionlimittheoryplustimes}
In this section we will consider the limit theory of prime fields with some non-trivial additive character and some (sufficiently generic) multiplicative character. We will prove a generalisation of Ax's result in this context. For the purely additive case this was done in \cite{Hrushovski2021AxsTW}. In our case however the method used in \cite{Hrushovski2021AxsTW} does not work and instead we employ a classical result from additive combinatorics, the Erdős–Turán–Koksma inequality. Before stating it, we recall a basic fact.

\begin{fact}(Lemma 3 in \cite{IwasawaKummerextensions})\label{factiwasawakummerextension}
 For an algebraic number field $L$ let $\mathrm{Tor}(L)$ denote its group of roots of unity. Then $L^{\times}$ is isomorphic to the direct sum of $\mathrm{Tor}(L)$ and $\mathbb{Z}^{(\omega)}$, a free abelian group on countably many generators.
\end{fact}

The difficulty for $\mathrm{PF}^{+,\times}$ is that there is no sensible multiplicative analogue of Lemma 3.18 in \cite{Hrushovski2021AxsTW} and we thus have to use different techniques to prove our result about the limit theory. In \cite{Hrushovski2021AxsTW} it was already mentioned that one could possible evoke an effective version of Weyl's equidistribution criterion. We will follow this approach in our case. Before stating Lemma \ref{lemmafinitefieldmultiplicativecorrection} we recall the necessary standard result from additive combinatorics that we will use. It is referred to as the \textit{Erdős–Turán–Koksma inequality} in the literature and was proven by Erdős and Turán
in the one-dimensional case in \cite{erdosturan1} and \cite{erdosturan2} and generalised by Koksma to higher dimensions in \cite{Koksma1950SomeTO}.

\begin{notation}
    Let $\bar{x}=\bar{x}_{1},\bar{x}_{2},\dots$ be a sequence of tuples in $[0,1)^{d}$. For a set $B\subseteq [0,1)^{d}$ we write $|B\cap X_{n}|$ for the number $|\{1\leq i\leq n\,|\,\bar{x}_{i}\in B\}|$.
    The discrepancy $\mathrm{Dis}(X_{n})$ of the finite sequence
    $X_{n}=\bar{x}_{1},\dots,\bar{x}_{n}$ is defined as follows.
    \[\mathrm{Dis}(X_{n}):=\sup_{B\in J}\left|\frac{|B\cap X_{n}|}{n}-\mu(B)\right|\]where $\mu$ is the Lebesgue-measure and $J$ the set of $d$-dimensional boxes defined as $[a_{1},b_{1})\times\cdots\times [a_{d},b_{d})$ where $0\leq a_{i}\leq b_{i}\leq1$ for all $1\leq i\leq d$.
\end{notation}
We present the Erdős–Turán–Koksma inequality as on page 116 in \cite{kuipers1974uniform}. 
\begin{fact}\label{facterdosturankoksmainequality}
    We work in the above notation. Let $H\in\mathbb{N}$, $\bar{h}\in\mathbb{Z}^{d}$, $z(\bar{h}):=\prod_{1\leq j\leq d}\max\{1,|h_{j}|\}$ and $|\bar{h}|_{\infty}:=\max_{1\leq j\leq d}|h_{j}|$. Then, there exists some constant $C_{d}>0$ only depending on $d$ such that
    \[\mathrm{Dis}(X_{n})\leq C_{d}\left(\frac{1}{H}+\sum_{0<|\bar{h}|_{\infty}\leq H}\frac{1}{z(\bar{h})}\left|\frac{1}{n}\sum_{1\leq i\leq n}e^{2\pi i \langle\bar{h},\bar{x}_{i}\rangle}\right|\right).\]
\end{fact}

We now apply the Erdős–Turán–Koksma inequality to obtain the following Lemma that will be the key ingredient in our later proofs.

\begin{lemma}\label{lemmafinitefieldmultiplicativecorrection}
Let $(\bar{\gamma}_{m})_{m\in\mathbb{N}}$ be a sequence of $d$-tuples $\bar{\gamma}_{m}=(\gamma_{m,1},\dots,\gamma_{m,d})\in \mathbb{T}^{d}$ in the $d$-Torus. Further, let $U$ be a non-empty open subset of $\mathbb{T}^{d}$ and let $R,K,f\in\mathbb{N}$ and $1\leq f\leq R$. Assume that for any $H\in\mathbb{N}$ there is an infinite subsequence $(\bar{\gamma}_{m})_{m\in N(H)}$ such that for all $\bar{\alpha}\in\mathbb{Z}^{d}$ with $0<|\bar{\alpha}|_{\infty}\leq H$ and all $m\in N(H)$ the inequality $\prod_{i=1}^{d}\gamma_{m,i}^{\alpha_{i}}\neq1$ holds.\\
Then, there is an infinite subsequence $(\bar{\gamma}_{m})_{m\in M}$ (of our initial sequence) such that for all $m\in M$ there is $l_{m}\in \mathbb{N}$ such that
\begin{itemize}
    \item $\bar{\gamma}_{m}^{l_{m}}:=(\gamma_{m,1}^{l_{m}},\dots,\gamma_{m,d}^{l_{m}})\in U$,
    \item $l_{m}\mod R=f$,
    \item For $\beta_{m,i}:=\gamma_{m,i}^{l_{m}}\in S^{1}$ we have $\max_{1\leq i\leq d}\mathrm{ord}(\beta_{m,i})\geq K$ where $\mathrm{ord}$ denotes the multiplicative order of an element in $S^{1}$.
\end{itemize}
% {\color{red}
% Let $(p_{m})_{m\in\mathbb{N}}$ be a sequence of strictly increasing prime numbers and $\chi_{p_{m}}$ multiplicative characters on $\mathbb{F}_{p_{m}}$. Let $\bar{c}_{m}=(c_{m_{1}},\dots,c_{m_{d}})\in\mathbb{F}_{p_{m}}^{d}$ be a sequence of $d$-tuples. Further, let $U$ be a non-empty open subset of $\mathbb{T}^{d}$ and let $R,K,f\in\mathbb{N}$ and $1\leq f\leq R$.\\
% Assume that for any $H\in\mathbb{N}$ there is an infinite subsequence $(p_{m})_{m\in M(H)}$ such that for all $\bar{\alpha}\in\mathbb{Z}^{d}$ with $0<|\bar{\alpha}|_{\infty}\leq H$ and all $m\in\mathbb{N}$ the inequality $\prod_{i=1}^{d}\chi_{p_{m}}(c_{m_{i}})^{\alpha_{i}}\neq1$ holds.\\
% Then, there is an infinite subsequence $(p_{m})_{m\in M}$ (of our initial sequence) such that for all $m\in M$ there is $l_{m}\in M(H)$ such that
% \begin{itemize}
%     \item $\chi^{(d)}_{p_{m}}(\bar{c}_{m})^{l_{m}}\in U$,
%     \item $l_{m}\mod R=f$,
%     \item For $\beta_{m,i}:=\chi_{p_{m}}(c_{m_{i}})^{l_{m}}\in S^{1}$ we have $\max_{1\leq i\leq d}\mathrm{ord}(\beta_{m,i})\geq K$ where $\mathrm{ord}$ denotes the multiplicative order of an element in $S^{1}$.
% \end{itemize}
% }
\end{lemma}
\begin{proof}
We define for every $m\in\mathbb{N}$ a sequence $\bar{x}^{m}=\bar{x}_{1}^{m},\bar{x}^{m}_{2},\dots$ of tuples in $[0,1)^{d}$ as follows. We take $x^{m}_{1,j}$ for $1\leq j\leq d$ to be the element from $[0,1)$ such that $\gamma_{m,j}=\exp(2\pi i x^{m}_{1,j})$. Next, for $k\geq 2$ and $1\leq j\leq d$ we set $x^{m}_{k,j}:=k\cdot x^{m}_{1,j} \mod \mathbb{Z}$. Here we identify $[0,1)$ with $\mathbb{R}/\mathbb{Z}$ and to lighten the notation we will not always write $\text{mod}\;\mathbb{Z}$ when it is clear from the context. For example below by $\langle\bar{\alpha},\bar{x}_{i}^{m}\rangle$ we will always mean $\langle\bar{\alpha},\bar{x}_{i}^{m}\rangle\mod\mathbb{Z}$.
For any nonzero integer vector $(\alpha_{1},\dots,\alpha_{d})$ we have 
 \[\prod_{i=1}^{d}\gamma_{m,i}^{\alpha_{i}}=1\;\iff\;\sum_{1\leq i\leq d}\alpha_{i}x_{1,i}^{m}=0\mod\mathbb{Z} \;\iff\;e^{2\pi i \langle\bar{\alpha},\bar{x}_{1}^{m}\rangle}=1.\]
 Further, note that $ \langle\bar{\alpha},\bar{x}^{m}_{k}\rangle= k \langle\bar{\alpha},\bar{x}^{m}_{1}\rangle$ holds for any $k\in\mathbb{N}$. It follows that $\frac{1}{n}\sum_{1\leq j\leq n}e^{2\pi i \langle\bar{\alpha},\bar{x}_{j}^{m}\rangle}\rightarrow 0$ for $n\rightarrow\infty$ if and only if $\prod_{i=1}^{d}\gamma_{m,i}^{\alpha_{i}}\neq1$. 
 %(As $m$ is fixed the sequence will be eventually periodic (of length $\leq p_{m}-1$) and the sub-sum corresponding to every such period is $0$.) {\color{red}Otherwise use Weyl's criterion.}
 Now let $H\in\mathbb{N}$ be given. We pass to the subsequence $(\bar{\gamma}_{m})_{m\in N(H)}$ as given by the assumption. Then, for all $m\in N(H)$, $\prod_{i=1}^{d}\gamma_{m,i}^{h_{i}}\neq 1$ holds for all $\bar{h}\in\mathbb{Z}^{d}$ with $|\bar{h}|_{\infty}\leq H$. It now follows from Fact \ref{facterdosturankoksmainequality} for $X_{n}^{m}:=\bar{x}_{1}^{m},\dots,\bar{x}_{n}^{m}$ that $\mathrm{Dis}(X_{n}^{m})\leq \frac{C_{d}}{H}+C_{d}\cdot\epsilon(n)$ for all $n\in\mathbb{N}$ where $\epsilon(n)\rightarrow 0$ for $n\rightarrow \infty$.\\
 Finally, w.l.o.g. we can assume that $U\subseteq \mathbb{T}^{d}$ is the image of an open box in $[0,1)^{d}$ under the map $x\rightarrow \exp(2\pi i x)$. Let $\lambda$ be the pushforward of the Lebesgue-measure on $[0,1)^{d}$. It then follows that, if $\mathrm{Dis}(X_{n}^{m})<\lambda(U)$, then there is some $1\leq l\leq n$ such that \[(e^{2\pi i x^{m}_{l,1}},\dots,e^{2\pi i x^{m}_{l,d}})=(\gamma_{m,1}^{l},\dots,\gamma_{m,d}^{l})\in U\] 
Thus, if $H$ is chosen big enough, that is $\frac{C_{d}}{H}<\lambda(U)$, then for all $m\in N(H)$ we indeed find $l_{m}\in\mathbb{N}$ such that $\bar{\gamma}^{l_{m}}\in U$.\\
For the second point first note that the above argumentation works similarly when defining $x^{m}_{k,j}:=(k\cdot R+f)\cdot x^{m}_{1,j} \mod \mathbb{Z}$. We then have $ \langle\bar{\alpha},\bar{x}^{m}_{k}\rangle= (k\cdot R+f) \langle\bar{\alpha},\bar{x}^{m}_{1}\rangle$ for any $k\in\mathbb{N}$.
Thus, it follows that \[\frac{1}{n}\sum_{1\leq j\leq n}e^{2\pi i \langle\bar{\alpha},\bar{x}_{j}\rangle}=e^{2\pi if\langle\bar{\alpha},\bar{x}_{1}^{m}\rangle}\cdot\frac{1}{n}\sum_{1\leq j\leq n}e^{2\pi ij\cdot (kR \langle\bar{\alpha},\bar{x}_{1}^{m}\rangle)}\]
As before, it now follows that for all $m\in N(H\cdot R)$ we can find $l_{m}$ such that $\bar{\gamma}^{l_{m}}\in U$ and $l_{m}$ is a multiple of $(k\cdot R+f)$.
For the third bullet point, it suffices to apply the above to some open subset $U^{\prime}\subseteq U$ such that at least the projection to one coordinate of $U^{\prime}$ does not contain any element of order $\leq K$.
\end{proof}

\begin{notation}\label{notationisomtypepftimes}Let $(F,\Psi,\chi)\models \mathrm{PF}^{+,\times}$ and write $\mathbb{Q}^{\mathrm{alg}}$ as an increasing union of finite Galois-extensions $(L_{n})_{n\in\mathbb{N}}$ of $\mathbb{Q}$. For each $n$ denote by $(F_{n},\Psi,\chi)$ the induced structure on $F_{n}:=F\cap L_{n}$. Further $\beta_{n}$ will denote a generator of $F_{n}$ over $\mathbb{Q}$. Let $\phi_{n}(x)$ be an $\mathcal{L}_{\mathrm{ring}}$-formula that describes the $\mathcal{L}_{\mathrm{ring}}$-isomorphism type of $F_{n}$ in the following sense: For any $F^{\prime}\models \mathrm{PF}_{0}$ we have that
      $F^{\prime}\models \exists \beta_{n}\phi_{n}(\beta_{n})$ if and only if $F_{n}\cong_{\mathcal{L}_{\mathrm{ring}}} L_{n}\cap F^{\prime}$.\footnote{See (6.18) in \cite{chatzidakislecturenotespsf} for an explicit construction.}\\
      Now we fix such a $\beta_{n}$ as well as a tuple $\bar{\alpha}=(\alpha_{1},\dots,\alpha_{e})\in F_{n}$ that (multiplicatively) generates a free abelian subgroup of the multiplicative group $F_{n}^{\times}$ isomorphic to $\mathbb{Z}^{e}$. Let $\zeta_{1},\dots,\zeta_{k}$ be an enumeration of all roots of unity contained in $F_{n}$ and $\iota_{1},\dots,\iota_{v}$ a $\mathbb{Q}$-basis of $F_{n}$. We write the $\alpha_{i}$ as $\alpha_{i}=h_{i}(\beta_{n})$, where $h_{i}(X)\in\mathbb{Q}[X]$. Similarly write $\zeta_{j}=g_{j}(\beta_{n})$ and $\iota_{s}=f_{s}(\beta_{n})$ for $g_{j}(X),f_{s}(X)\in\mathbb{Q}[X]$.
\end{notation}

\begin{definition}\label{definitionthetasentencepfplustimes}
 We work in $(F,\Psi,\chi)\models \mathrm{PF}^{+,\times}$ and consider for some fixed $n\in\mathbb{N}$ a triple $I=(\bar{\alpha},\bar{\iota},l)$ as well as $\bar{\zeta}$ where $l\in\mathbb{N}$ and $\bar{\alpha},\bar{\iota},\bar{\zeta}$ are as in Notation \ref{notationisomtypepftimes}.
 % For some given $n,l\in\mathbb{N}$ and $\bar{\alpha},\bar{\zeta},\bar{\iota}$ we define an $\mathcal{L}_{+,\times}$-sentence $\Theta_{n,l}^{\bar{\alpha}}$ in the following. 
Let $r_{i}=\chi(\alpha_{i})$ for $1\leq i\leq e$, $t_{j}=\chi(\zeta_{j})$ for $1\leq j\leq k$ and $u_{s}=\Psi(\iota_{s})$ for $1\leq s\leq v$
 as determined by the $\mathcal{L}_{+,\times}$-theory of $(F_{n},\Psi,\chi)$. Next, we define
\[\phi^{\times}_{n,I}(\beta_{n}):=
\bigwedge_{1\leq i\leq e}\left|\chi\left(h_{i}(\beta_{n})\right)-r_{i}\right|\leq\frac{1}{l}\;\land\;\bigwedge_{1\leq j\leq k}\left|\chi\left(g_{j}(\beta_{n})\right)-t_{j}\right|\leq\frac{1}{l}
\]
and
\[\phi^{+}_{n,I}(\beta_{n})=\bigwedge_{1\leq i\leq v}\left|\Psi\left(f_{i}(\beta_{n})\right)-u_{i}\right|\leq\frac{1}{l}.\]
Finally, we define the sentence $\Theta_{n,I}$ as follows and set $\Theta_{n}=\cup_{I}\Theta_{n,I}$.   
\[
\Theta_{n,I}:=\exists\beta_{n}\left(\phi_{n}(\beta_{n})\;\land\;\phi^{+}_{n,l}(\beta_{n})\;\land\;\phi^{\times}_{n,l}(\beta_{n})\right).
\]
%       Now we define the set of $\mathcal{L}_{+,\sigma}$-sentences $\Theta_{n}$ as $\Theta_{n}:=\{\Theta_{n,s}\,|\,s\in\mathbb{N}\}$.
\end{definition}

\begin{remark}
    There is no problem in working with an existential quantifier instead of an infimum in the sentences $\Theta_{n,I}$. Indeed, modulo $\mathrm{PF}^{+,\times}$, the sentences $\Theta_{n,I}$ are expressible in $\mathcal{L}_{+,\times}$ since the quantification is over a finite (strongly) definable set and thus infimum and existential quantifier coincide.
\end{remark}

\begin{lemma}\label{lemmasentencesdescribingisomtypepfplustimes}

    The set of sentences $\Theta_{n}$ describes the $\mathcal{L}_{+,\times}$-isomorphism type of $(F_{n},\Psi,\chi)$.

\end{lemma}

\begin{proof}
Assume that $(F_{n},\Psi,\chi)$ and $(F_{n}^{\prime},\Psi^{\prime},\chi^{\prime})$ both satisfy some given $\Theta_{n}$.
    By the piegeonhole-principle we find some $\beta_{n}$ (as in Notation \ref{notationisomtypepftimes}) that witnesses all the $\Theta_{n,I}$ and similary we find such a $\beta_{n}^{\prime}$ in $F^{\prime}$. Then, taking into account Fact \ref{factiwasawakummerextension}, it is direct to see that $\beta_{n}\rightarrow\beta_{n}^{\prime}$ induces an $\mathcal{L}_{+,\times}$-isomorphism between $(F_{n},\Psi,\chi)$ and $(F_{n}^{\prime},\Psi^{\prime},\chi^{\prime})$.
\end{proof}

\begin{lemma}\label{lemmaonesentencesufficientPFplustimes}
    Let $\Theta_{n_{1},I_{1}},\dots,\Theta_{n_{k},I_{k}}$ be as in Definition \ref{definitionthetasentencepfplustimes}. We can find some $\Theta_{n,I}$ that implies all the $\Theta_{n_{1},I_{1}},\dots,\Theta_{n_{k},I_{k}}$.
    
\end{lemma}
\begin{proof}
    Let $F_{n}$ be a Galois extension of $\mathbb{Q}$ that contains all the $F_{n_{1}},\dots,F_{n_{k}}$. Next, choose $\bar{\alpha}$ and $\bar{\iota}$ such that all the elements of all $\bar{\alpha}_{i}$ and $\bar{\iota}_{i}$ for $1\leq i\leq k$ are contained in the additive, resp. multiplicative, subgroup generated by $\bar{\alpha}$, resp. $\bar{\iota}$. Finally, given this data it suffices to choose $l$ sufficiently big to obtain $\Theta_{n,I}$ as wanted.
\end{proof}

\begin{theorem}\label{theorempfplustimesislimittheory}
 Let $T_{\Psi,\chi}$ be the common $\mathcal{L}_{+,\times}$-theory of all prime fields with an additive character and multiplicative character.
The theory $\mathrm{PF}^{+,\times}$ is given by $T_{\Psi,\chi}$ together with a set of axioms stating that the order of the characters is infinite (which implies that the characteristic of the field is $0$). 
\end{theorem}
\begin{proof}
By Theorem \ref{theoremUPoffinitefieldsmodelofpfplustimes} any ultraproduct $\prod_{\mathcal{U}}(\mathbb{F}_{q},\Psi_{q},\chi_{q})$ of characteristic $0$ with additive and multiplicative character of infinite order is a model of $\mathrm{PF}^{+,\times}$. It remains to show that for any completion $T$ of $\mathrm{PF}^{+,\times}$, we find such an ultraproduct that is a model of $T$. From Corollary \ref{corollaryeltequivpfplustimes} it follows, that it suffices to find for some given model $(F,\Psi,\chi)\models\mathrm{PF}^{+,\times}$ an ultraproduct $(F^{\prime},\Psi^{\prime},\chi^{\prime})=\prod_{\mathcal{U}}(\mathbb{F}_{q},\Psi_{q},\chi_{q})$ as above such that $\bar{\mathbb{Q}}\cap F$ and $\bar{\mathbb{Q}}\cap F^{\prime}$ are $\mathcal{L}_{+,\times}$-isomorphic. By Lemma \ref{lemmasentencesdescribingisomtypepfplustimes} the $\mathcal{L}_{+,\times}$-isomorphism type of $F\cap\bar{\mathbb{Q}}$ is determined by the set of sentences $\cup_{n\in\mathbb{N}}\Theta_{n}$. Using compactness and Lemma \ref{lemmaonesentencesufficientPFplustimes}, we thus have to show that given some $\Theta_{n,I}$, there is always an infinite strictly increasing sequence of primes $(p_{m})_{m\in\mathbb{N}}$ with additive and multiplicative characters $\Psi_{p_{m}}$ and $\chi_{p_{m}}$ of unbounded order such that $(\mathbb{F}_{p_{m}},\Psi_{p_{m}},\chi_{p_{m}})\models \Theta_{n,I}$.\\
We work with the notation (corresponding to $\Theta_{n,I}$) of Definition \ref{definitionthetasentencepfplustimes}. We apply the \v{C}ebotarev Density Theorem following Proposition 7 of \cite{ax-elttheoryoffinitefields} to obtain a strictly increasing sequence of primes $(p_{m})_{m\in\mathbb{N}}$ such that the $\mathcal{L}_{\mathrm{ring}}$-sentence describing the $\mathcal{L}_{\mathrm{ring}}$-isomorphism type of $F_{n}$ holds in $\mathbb{F}_{p_{m}}$ in the sense of \cite{ax-elttheoryoffinitefields}.\\
We will now repeatedly refine this sequence $(p_{m})_{m\in\mathbb{N}}$ always keeping an infinite subsequence. Let $R\in\mathbb{N}$ be the least common multiple of the finitely many $k\in\mathbb{N}$ such that $F_{n}$ contains all $k$-th roots of unity. As this is witnessed by the $\mathcal{L}_{\mathrm{ring}}$-isomorphism type of $F_{n}$, we can assume by passing to a subsequence that $R|(p_{m}-1)$ holds for all $m\in\mathbb{N}$. Let $\lambda_{R}$ be a generator of the group of $R$-th roots of unity (written as a rational polynomial in $\beta_{n}$) and $s_{R}:=\chi(\lambda_{R})$ in $(F,\Psi,\chi)$. Now let $\Tilde{\chi}_{p_{m}}$ be a multiplicative character of order $p_{m}-1$ on $\mathbb{F}_{p_{m}}$. By passing to a subsequence we can assume that there is $1\leq f\leq R$ such that for all $m\in\mathbb{N}$ the equality $\Tilde{\chi}_{p_{m}}(\lambda_{R})^{f}=s_{R}$ holds.\\
Recall that $\Tilde{\chi}_{p_{m}}$ is of order $p_{m}-1$ and $\{h_{1}(\beta_{n}),\dots,h_{e}(\beta_{n})\}$ generates (multiplicatively) a free abelian subgroup isomorphic to $\mathbb{Z}^{e}$ in $F_{n}$. From this it follows that for any $H\in\mathbb{N}$ we can pass to a subsequence such that then for all $p_{m}$ (from the subsequence) and for all non-zero integer vector $(\rho_{1},\dots,\rho_{e})$ with $|\bar{\rho}|_{\infty}\leq H$ we have that \[\prod_{1\leq i\leq e}\Tilde{\chi}_{p_{m}}(h_{i}(\beta_{n}))^{\rho_{i}}\neq1.\]
We apply Lemma \ref{lemmafinitefieldmultiplicativecorrection} to the sequence of $e$-tuples $(\Tilde{\chi}_{p_{m}}(h_{1}(\beta_{n})),\dots,\Tilde{\chi}_{p_{m}}(h_{e}(\beta_{n})))_{m\in\mathbb{N}}$ and obtain an infinite subsequence (which we again denote by) $(p_{m})_{m\in\mathbb{N}}$ such that for all $m\in\mathbb{N}$ there is some $r_{m}\in\mathbb{N}$ such that $|\Tilde{\chi}_{p_{m}}(h_{i}(\beta_{n}))^{r_{m}}-t_{i}|\leq \frac{1}{l}$ for all $1\leq i\leq e$ and $r_{m}\mod R=f$. Consequently, for all $m\in\mathbb{N}$ we have $\Tilde{\chi}_{p_{m}}(\lambda_{R})^{r_{m}}=s_{R}$. We then set $\chi_{p_{m}}:=(\Tilde{\chi}_{p_{m}})^{r_{m}}$. Note that using the third bullet point in Lemma \ref{lemmafinitefieldmultiplicativecorrection} it follows that for any $K\in\mathbb{N}$ we can find the $\chi_{p_{m}}$ as above such that moreover $\mathrm{ord}(\chi_{p_{m}})\geq K$.\\
    We continue to work in the sequence obtained thus far and let $\Tilde{\Psi}_{p_{m}}$ be some non-trivial (and thus of order $p_{m}$) additive character on $\mathbb{F}_{p_{m}}$. We use the exact same argument as in the multiplicative case but this time applied to the sequence of tuples $(\Tilde{\Psi}_{p_{m}}(f_{1}(\beta_{n})),\dots,\Tilde{\Psi}_{p_{m}}(f_{v}(\beta_{n})))_{m\in\mathbb{N}}$ to obtain an infinite subsequence $(p_{m})_{m\in\mathbb{N}}$ as well as $r_{m}\in\mathbb{N}$ such that $\left|\Tilde{\Psi}_{p_{m}}\left(f_{i}(\beta_{n})\right)^{r_{m}}-u_{i}\right|\leq\frac{1}{l}$ for all $1\leq i\leq v$ and $m\in\mathbb{N}$. We set $\Psi_{p_{m}}:=\Tilde{\Psi}_{p_{m}}^{r_{m}}$ and conclude that $(\mathbb{F}_{p_{m}},\Psi_{p_{m}},\chi_{p_{m}})\models \Theta_{n,I}$ for all $m\in\mathbb{N}$ which finishes the proof.\end{proof}
The multiplicative case differs from the additive one in the sense that there is no generic choice of a \textit{standard multiplicative character} since there is no generic choice of a generator of the multiplicative group. Or, in other words the choice of $\Psi_{p}^{\mathrm{stan}}(x)=\exp(2\pi i\frac{x}{p})$ as standard additive character (as in \cite{Hrushovski2021AxsTW}) is only \textit{standard} after accepting $1$ as a standard choice of the generator of the additive group. However, we can still  consider the question of the limit theory of prime and finite fields with standard additive character and some sufficiently generic multiplicative character.\\
As described in Section 5 in \cite{Hrushovski2021AxsTW} in the case of prime fields with standard additive character axioms must be added. The conjectural axiomatisation in \cite{Hrushovski2021AxsTW} is shown to relate to a number-theoretic conjecture due to Duke, Friedlander and Iwaniec (see Section 5.5 in \cite{Hrushovski2021AxsTW}). We will not make any progress in this direction but outline that, as in \cite{Hrushovski2021AxsTW}, when one passes to finite fields a version of Theorem \ref{theorempfplustimesislimittheory} can be obtained. Recall that the standard additive character $\Psi_{q}^{\mathrm{stan}}$ on a finite field $\mathbb{F}_{q}$ is defined as $\Psi_{q}^{\mathrm{stan}}:=\Psi_{p}^{\mathrm{stan}}(\mathrm{Tr}(x))$ for $\mathrm{Tr}:\mathbb{F}_{q}\rightarrow\mathbb{F}_{p}$ the trace. The proof is the same as before once we apply the result from \cite{Hrushovski2021AxsTW}.

\begin{theorem}
 Let $T_{\Psi,\chi}^{\mathrm{stan}}$ be the common $\mathcal{L}_{+,\times}$-theory of all finite fields with the standard additive character and some multiplicative character.
The theory $\mathrm{PF}^{+,\times}$ is given by $T_{\Psi,\chi}^{\mathrm{stan}}$ together with a set of axioms stating that the multiplicative character is of infinite order and that the characteristic of the field is $0$.   
\end{theorem}

\begin{proof}
We work in the same setup as in the proof of Theorem \ref{theorempfplustimesislimittheory} and have to find a sequence $q_{m}$ with $(\mathbb{F}_{q_{m}},\Psi_{q_{m}}^{\mathrm{stan}},\chi_{p_{m}})\models \Theta_{n,I}$. We find a sequence $(\mathbb{F}_{q_{m}},\Psi_{q_{m}}^{\mathrm{stan}})$ corresponding to the $\mathcal{L}_{+}$-isomorphism type of $(F_{n},\Psi)$ by Proposition 3.16 (and its proof) in \cite{Hrushovski2021AxsTW}. Then we use the exact same proof as for Theorem \ref{theorempfplustimesislimittheory} to find a subsequence as well as multiplicative characters $\chi_{q_{m}}$ such that $(\mathbb{F}_{q_{m}},\Psi_{q_{m}}^{\mathrm{stan}},\chi_{p_{m}})\models \Theta_{n,I}$.\end{proof}

\begin{remark}
Note that the proof of Theorem \ref{theorempfplustimesislimittheory} yields something slightly stronger. For any $\sigma\in\mathrm{Gal}(\bar{\mathbb{Q}}/\mathbb{Q})$ and $\Psi,\chi$ additive and multiplicative characters on $L=\mathrm{Fix}(\sigma)$ there is some $(F,\Psi,\chi)=\prod_{\mathcal{U}}(\mathbb{F}_{q},\Psi_{q},\chi_{q})$ such that $L$ and $F\cap\bar{\mathbb{Q}}$ (together with the respective characters) are $\mathcal{L}_{+,\times}$-isomorphic. In particular this gives a description of all the completions $T$ of $\mathrm{PF}^{+,\times}$.
\end{remark}

\section{Definable integration}\label{sectiondefinableintegration}

\subsection{Facts around the Chatzidakis-van den Dries-Macintyre counting measure}

We start by recalling the main result of \cite{ChatzidakisMacvdd}.
\begin{notation}
    In the following let $\phi(\bar{x},\bar{y})$ be an $\mathcal{L}_{\mathrm{ring}}$-formula. For some field $K$ and $\bar{a}\in K^{|\bar{y}|}$ we write $\phi(K,\bar{a})$ for the corresponding definable set in $K^{|\bar{x}|}$. If not specified otherwise, $F$ will be a pseudofinite field of characteristic $0$.
\end{notation}

\begin{fact}(Theorem 3.7 in \cite{ChatzidakisMacvdd})\label{factsizeofdebsets}
    For any $\mathcal{L}_{\mathrm{ring}}$-formula $\phi(\bar{x},\bar{y})$ there is a positive constant $C$ and a finite set $D_{\phi}$ of tuples $(d_{1},\mu_{1}),\dots,(d_{k},\mu_{k})\in\{0,\dots,|\bar{x}|\}\times\mathbb{Q}_{>0}\cup\{0,0\}$ as well as $\mathcal{L}_{\mathrm{ring}}$-formulas $\psi_{1}(\bar{y}),\dots,\psi_{k}(\bar{y})$ such that for all finite fields $\mathbb{F}_{q}$ and $\bar{a}\in \mathbb{F}_{q}^{|\bar{y}|}$ there is $1\leq i\leq k$ such that $\bar{a}\in\psi_{i}(\mathbb{F}_{q})$ and  
    \[||\phi(\mathbb{F}_{q},\bar{a})|-\mu_{i} q^{d_{i}}|\leq Cq^{d_{i}-\frac{1}{2}}\]
    holds if and only if $\bar{a}\in\psi_{i}(\mathbb{F}_{q})$. Consequently, for any pseudofinite $F$ and $\phi(\bar{x},\bar{y})$ the formulas $\psi_{1}(\bar{y}),\dots,\psi_{k}(\bar{y})$ yield a partition of $F^{|\bar{y}|}$.
\end{fact}

\begin{notation}
   For $\bar{a}\in F^{|\bar{y}|}$ we denote by $(d_{\phi}(\bar{a}),\mu_{\phi}(\bar{a}))$ the tuple corresponding to the unique $1\leq l\leq k$ such that $F \models \psi_{l}(\bar{a})$. We call $d_{\phi}(\bar{a})$ the \textit{dimension} of the definable set $\phi(\bar{x},\bar{a})$ and $\mu_{\phi}(\bar{a})$ its \textit{multiplicity}.
\end{notation}

We will now collect some intermediate results from the proof of Fact \ref{factsizeofdebsets} and a slight variation of them, similar to what is presented in Theorem 7 in \cite{Kowalski2005ExponentialSO}.

\begin{notation}
    $W_{\bar{y}}$ will denote an algebraic set defined over $\mathbb{Z}$ in the variables $(\bar{x},\bar{y})$ where we will think of the $\bar{y}$ as placeholders for parameters. I.e., when we write $W_{\bar{a}}$ for a tuple $\bar{a}$ of the length of $\bar{y}$ in some field $F$, then this denotes the corresponding algebraic set in the variables $\bar{x}$ with parameters $\bar{a}$. Similarly, when we say that something holds uniformly in $\bar{y}$ for $W_{\bar{y}}$ we mean that it holds for the corresponding family of algebraic sets in the variables $\bar{x}$.\\
    In the same manner we treat definable sets $B_{\bar{y}}$.\\
    Next, for any prime power $q$ we fix $\beta_{q}$ to be a function $\beta_{q}:\mathbb{F}_{q}^{|\bar{x}|}\rightarrow \mathbb{C}$ such that the $\beta_{q}$ are uniformly bounded (in $q$) and $\beta:F\rightarrow \mathbb{C}$ will be the corresponding ultralimit where $F$ is a characteristic $0$ ultraproduct of finite fields. If the field in which we work is clear from the context, we might leave out the subscript and write $\beta$ for $\beta_{q}$. By $\pi$ we will denote the projection to the first $|\bar{x}|$-many variables whenever this is well defined.
\end{notation}

\begin{fact}\label{factkowlaskireductiontoalgsets}(See Theorem 7 in \cite{Kowalski2005ExponentialSO} for a different but closely related calculation.)
For any $\phi(\bar{x},\bar{y})$ there exist $e,l\in\mathbb{N}$ as well as $\mathcal{L}_{\mathrm{ring}}$-formulas $\psi_{1}(\bar{x},\bar{y},\bar{z}),\dots,\psi_{l}(\bar{x},\bar{y},\bar{z})$ and algebraic sets $W^{i,j}_{\bar{y}\bar{z}}$ (in the variables $\bar{v}=(\bar{x},\bar{u})$) for $1\leq i\leq l$ and $1\leq j\leq e$ and $k_{0}\in\mathbb{N}$ such that the following holds.
For any $\mathbb{F}_{q}$ with $q\geq k_{0}$ there exists $\bar{b}\in\mathbb{F}_{q}^{|\bar{z}|}$ such that for any $\bar{a}\in\mathbb{F}_{q}^{|\bar{y}|}$ we have:
\begin{itemize}
	\item $\phi(\mathbb{F}_{q},\bar{a})$ is given by the disjoint union $\psi_{1}(\mathbb{F}_{q},\bar{a},\bar{b})\,\dot{\cup}\cdots\dot{\cup}\,\psi_{l}(\mathbb{F}_{q},\bar{a},\bar{b})$.
    \item For any $1\leq i\leq l$ we have for $R^{i}_{q}:=|\psi_{i}(\mathbb{F}_{q},\bar{a},\bar{b})|$, $S^{i,j}_{q}:=|W^{i,j}_{\bar{a}\bar{b}}(\mathbb{F}_{q})|$ and $\pi(S)^{i,j}_{q}:=|\pi(W^{i,j}_{\bar{a}\bar{b}}(\mathbb{F}_{q}))|$ that
    \[\frac{1}{R^{i}_{q}}\sum_{\bar{x}\in\psi_{i}(\mathbb{F}_{q},\bar{a},\bar{b})}\beta(\bar{x})=\sum_{1\leq j\leq e}\frac{\pi(S)^{i,j}_{q}}{R^{i}_{q}}\left(\frac{1}{S^{i,j}_{q}}\sum_{\bar{v}\in W^{i,j}_{\bar{a}\bar{b}}(\mathbb{F}_{q})}\beta(\pi(\bar{v}))\right)\]
%	\item For any $1\leq i\leq l$ we have	\[\sum_{\bar{x}\in\psi_{i}(\mathbb{F}_{q},\bar{a},\bar{b})}\beta(\bar{x})=\sum_{1\leq j\leq e}\frac{(-1)^{j+1}}{j!}\sum_{\bar{v}\in W^{i,j}_{\bar{a}\bar{b}}(\mathbb{F}_{q})}\beta(\pi(\bar{v})).\]
\end{itemize}

Moreover, if $\phi(\bar{x},\bar{a})$ has dimension $1$ for all $\bar{a}\in\mathbb{F}_{q}^{|\bar{y}|}$, then the $W^{i,j}_{\bar{y}\bar{z}}$ can all be assumed to be algebraic curves (not necessarily irreducible).\\
%The above statements work similarly if in the assumptions we let $\bar{y}$ vary over some $\mathcal{L}_{\mathrm{ring}}$-definable set instead of the whole field.
\end{fact}

\begin{proof}
The discussion on page 123 and 124 in \cite{ChatzidakisMacvdd} (or see the proof of Theorem 7 in \cite{Kowalski2005ExponentialSO}) yields that for every $1\leq i\leq l$ we obtain a partition
    \[\psi_{i}(\mathbb{F}_{q},\bar{a},\bar{b})=\dot{\bigcup}_{1\leq j\leq e}\pi\left(W^{i,j}_{\bar{a}\bar{b}}(\mathbb{F}_{q})\right)\]and all (non-empty) fibers of the map $\pi:W^{i,j}_{\bar{a}\bar{b}}(\mathbb{F}_{q})\rightarrow\psi_{i}(\mathbb{F}_{q},\bar{a},\bar{b})$ are of size $j$ and thus $j\cdot\pi(S)^{i,j}_{q}=S^{i,j}_{q}$.
\end{proof}

We will need the following statement that was given in \cite{ChatzidakisMacvdd} under the name \textit{decomposition-intersection procedure}. See \cite{vandendriesSchmidt} for some of the underlying results.

\begin{fact}(Theorem 1.15 in \cite{ChatzidakisMacvdd}.)\label{factdecompintersect}
For any $W_{\bar{y}}$ we find an algebraic set $\Tilde{W}_{\bar{d}(\bar{y})}=V_{\bar{d}(\bar{y})}^{1}\cup\cdots\cup V_{\bar{d}(\bar{y})}^{t}$ and $k_{0}\in\mathbb{N}$ such that the following holds. For any field $K$ of characteristic $0$ or $\geq k_{0}$ and for any $\bar{a}\in K^{n}$ we have that
    \begin{itemize}
        \item $W_{\bar{a}}(K)=\Tilde{W}_{\bar{d}(\bar{a})}(K)$
 %       =V_{\bar{d}(\bar{a})}^{1}(F)\dot{\cup}\cdots\dot{\cup} V_{\bar{d}(\bar{a})}^{t}(F)$.
        \item For any $1\leq i\leq t$ the algebraic set $V_{\bar{d}(\bar{a})}^{i}$ is absolutely irreducible or $V_{\bar{d}(\bar{a})}^{i}(K)=\emptyset$.
        \item The parameters $\bar{d}(\bar{a})$ are in $acl(\bar{a})\cap K$ uniformly (in $\bar{a}$).\footnote{There exist polynomials $h_{1}(Z,\bar{Y}),\dots, h_{l}(Z,\bar{Y})\in\mathbb{Z}[Z,\bar{Y}]$ such that for every $\bar{a}\in K^{n}$ and every $1\leq j\leq l$ the element $d_{j}(\bar{a})$ (the $j$-th coordinate of the tuple $\bar{d}(\bar{a})$) is given by a root of the polynomial $h_j(Z,\bar{a})$.}
    \end{itemize}
      
\end{fact}

The above Fact \ref{factsizeofdebsets} gives rise to a measure on definable subsets of a definable set in a pseudofinite field as follows.

\begin{definition}\label{definitionmeasureonPF}
    For $B_{\bar{y}}$ defined by $\phi(\bar{x},\bar{y})$ we denote by $\mu_{B_{\bar{y}}}$ the family of maps defined as follows. For $\bar{a}\in F^{|\bar{y}|}$ and a definable subset $D_{\bar{a}}\subseteq B_{\bar{a}}$ defined by the formula $\Tilde{\phi}(\bar{x},\bar{y})$ we set
    \[\mu_{B_{\bar{a}}}(D_{\bar{a}}):=\begin{cases}
       \frac{\mu_{\Tilde{\phi}}(\bar{a})}{\mu_{\phi}(\bar{a})}\;\;\;\;,\text{if}\;d_{\Tilde{\phi}}(\bar{a})=d_{\phi}(\bar{a})\;\text{and}\;\mu_{\phi}(\bar{a})\neq 0.\\
       0\;\;\;\;\;\;\;\;\;,\text{otherwise}.
    \end{cases}\]
\end{definition}

\begin{fact}(4.10 in \cite{ChatzidakisMacvdd})\label{factzoemasureiscountingmeasure}
    For every $B_{\bar{y}},F$ and every $\bar{a}\in F^{|\bar{y}|}$ the map $\mu_{B_{\bar{a}}}$ yields a finitely additive probability measure on the boolean algebra of $\bar{a}$-definable subsets of $B_{\bar{a}}$. Moreover, if $F=\prod_{\mathcal{U}}\mathbb{F}_{q}$ is an ultraproduct of finite fields, then for any definable $D_{\bar{a}}$ we have that $\mu_{B_{\bar{a}}}(D_{\bar{a}})$ is given as the ultralimit (along $\mathcal{U}$) of the quantity $\frac{|D_{\bar{a}_{q}}(\mathbb{F}_{q})|}{|B_{\bar{a}_{q}}(\mathbb{F}_{q})|}$ for any representative $\bar{a}=(\bar{a}_{q})/\mathcal{U}$.
    %   Moreover $\mu_{B_{\bar{a}}}$ satisfies Fubini in the following sense: For any $\bar{a}$-definable set $H$ and function $f:B_{\bar{a}}\rightarrow H$ that is onto, $\bar{a}$-definable
\end{fact}

\subsection{Uniform definable integration of definable predicates}

We now return to work in the language $\mathcal{L}_{+,\times}$. Our interest will be the uniform definability of the average of an $\mathcal{L}_{+,\times}$-definable predicate over an $\mathcal{L}_{\mathrm{ring}}$-definable set. The archetypical example will be that of a (generalised) character sum. Character sums are ubiquitous in Number Theory. However, we refrain from an extended presentation and refer the reader for example to Chapter 11 of \cite{iwaniec2021analytic}. They will also be prominent in the proof of the main result of this section, so we recall the definition.
\begin{example}
    We call a generalised character sum (over a definable set) an expression of the form
    \[\sum_{\bar{x}\in X(\mathbb{F}_{q})}\Psi(h(\bar{x}))\chi(g(\bar{x}))\]
    where $X$ is an $\mathcal{L}_{\mathrm{ring}}$-definable set and $h,g$ are $\mathcal{L}_{\mathrm{ring}}$-definable functions $\mathbb{F}_{q}^{|\bar{x}|}\rightarrow\mathbb{F}_{q}$.
\end{example}
In \cite{Kowalski2005ExponentialSO} Kowalski gives the following bound on character sums over definable sets over finite fields.
\begin{fact}(Theorem 2 in \cite{Kowalski2005ExponentialSO})There exist constants $K,c>0$ such that for $h,g\in\mathbb{Q}(\bar{X})$ and $X$ an $\mathcal{L}_{\mathrm{ring}}$-definable set whenever for all $b\in\mathbb{F}_{q}$ we have $|X(\mathbb{F}_{q})\cap h^{-1}(b)|<c|X(\mathbb{F}_{q})|$, then
\[\sum_{\bar{x}\in X(\mathbb{F}_{q})}\Psi(h(\bar{x}))\chi(g(\bar{x}))\leq Kq^{-\frac{1}{2}}|X(\mathbb{F}_{q})|.\]    
\end{fact}

Kowalski however does not give bounds for multiplicative character sums over definable sets, that is, sums as above where additive character is trivial. He discusses why they are harder to deal with (Example 5 in \cite{Kowalski2005ExponentialSO}). We, however, need to be able to work with $\chi$ solely as well, as it is a predicate in our language. As we do not aim to give a concrete condition for cancellation in the character sum but are only interested in the fact that there is (a definable) such condition the possible complications are less problematic for us. Also, our $\chi$ is generic in the sense that it is not of finite order in the ultraproduct. Therefore we avoid Example 5 of \cite{Kowalski2005ExponentialSO}. However, note that we could also incorporate multiplicative characters of finite order as they are are already definable in $\mathrm{PF}$.\\
Another difference to \cite{Kowalski2005ExponentialSO} that is worth mentioning is that our proof only use the classical Weil bounds as given by Fact \ref{factperelmutterestimates} and do not employ the cohomological methods due to Grothendieck and Deligne.\\ 
Our result below is a generalisation Theorem 4.1 of \cite{Hrushovski2021AxsTW} to $\mathrm{PF}^{+,\times}$. We adapt the ideas of its proof to our case and provide some of the details left out in \cite{Hrushovski2021AxsTW}. Let us now state the main theorem of this section in its weak version. We will later show in Theorem \ref{corollaryremovinglowerbounds} how we can get rid of the lower bounds $k_{0}$ and $r_{0}$.
\begin{theorem}\label{theoremdefinableintegrationfiniteversion}
    Let $P(\bar{x},\bar{y})$ be an $\mathcal{L}_{+,\times}$-predicate and $B_{\bar{y}}$ an $\mathcal{L}_{\mathrm{ring}}$-definable set. For any finite field $\mathbb{F}_{q}$ where $q=p^{e}$ consider the function $f:\mathbb{F}_{q}^{|\bar{y}|}\rightarrow\mathbb{C}$ given by $f(\bar{y})=0$, if $B_{\bar{y}}=\emptyset$, and otherwise by\[f(\bar{y}):=\frac{\sum_{\bar{x}\in B_{\bar{y}}(\mathbb{F}_{q})}P(\bar{x},\bar{y})}{|B_{\bar{y}}(\mathbb{F}_{q})|}.\] 
Now for any $\epsilon>0$, there exist $k_{0},r_{0}\in\mathbb{N}$ and an $\mathcal{L}_{+,\times}$-predicate $\Tilde{P}(\bar{y})$  such that for any $p\geq k_{0}$, $\Psi$ non-trivial additive character on $\mathbb{F}_{q}$,  $\chi$ multiplicative character of order $\geq r_{0}$ on $\mathbb{F}_{q}$ and any $\bar{a}\in \mathbb{F}_{q}^{|\bar{y}|}$,
\[|\Tilde{P}(\bar{a})-f(\bar{a})|\leq\epsilon.\]
\end{theorem}

We will now proceed with the preparation for the proof of Theorem \ref{theoremdefinableintegrationfiniteversion}. Instead of working on the level of finite fields, similar as in \cite{Hrushovski2021AxsTW}, we will work in a characteristic $0$ ultraproduct thereof as it makes the presentation more readable.

\begin{convention}
From now on, if not stated otherwise, $(F,\Psi,\chi)\models \mathrm{PF}^{+,\times}$
will be an (characteristic $0$) ultraproduct of enriched finite fields, i.e., $(F,\Psi,\chi)=\prod_{\mathcal{U}}(\mathbb{F}_{q},\Psi_{q},\chi_{q})$.
\end{convention}

\begin{definition}\label{definitionintegrals}
We work in the notation of Theorem \ref{theoremdefinableintegrationfiniteversion}. We write the ultralimit of the functions $f:\mathbb{F}_{q}^{|\bar{y}|}\rightarrow\mathbb{C}$ as follows.
\[\int_{B_{\bar{y}}}P(\bar{x},\bar{y})d\mu_{B_{\bar{y}}}^{*}(\bar{x}):=\lim_{\mathcal{U}}f(\bar{y}).\]
%$\int_{B_{\bar{y}}}P(\bar{x},\bar{y})d\mu_{B_{\bar{y}}}^{*}(\bar{x})$ for the function $\lim_{\mathcal{U}}f(\bar{y}):F^{|\bar{y}|}\rightarrow\mathbb{C}$, that is, the ultralimit of the functions $f:\mathbb{F}_{q}^{|\bar{y}|}\rightarrow\mathbb{C}$. 
Moreover, we write \[\int_{D_{\bar{y}}}P(\bar{x},\bar{y})d\mu_{B_{\bar{y}}}^{*}(\bar{x}):=\int_{B_{\bar{y}}}\mathbbm{1}_{D_{\bar{y}}}\cdot P(\bar{x},\bar{y})d\mu_{B_{\bar{y}}}^{*}(\bar{x})\]
%$\int_{D_{\bar{y}}}P(\bar{x},\bar{y})d\mu_{B_{\bar{y}}}^{*}(\bar{x})$ for the expression $\int_{B_{\bar{y}}}\mathbbm{1}_{D_{\bar{y}}}\cdot P(\bar{x},\bar{y})d\mu_{B_{\bar{y}}}^{*}(\bar{x})$
where $D_{\bar{y}}\subseteq B_{\bar{y}}$ is an $\mathcal{L}_{\mathrm{ring}}$-definable subset.
\end{definition}

% \begin{definition}

% We denote by $\mu_{B_{\bar{y}}}^{*}$ the family of \textit{pseudofinite counting measures on} $B_{\bar{y}}$ where for every $\bar{a}\in F^{|\bar{y}|}$ the measure $\mu_{B_{\bar{a}}}^{*}$ assigns to every internal set $D\subseteq B_{\bar{a}}$ the ultralimt (along $\mathcal{U}$) of the quantity $\frac{|D_{q}(\mathbb{F}_{q})|}{|B_{\bar{a}_{q}}(\mathbb{F}_{q})|}$.
    
% \end{definition}

\begin{remark}
    Note that the notation from Definition \ref{definitionintegrals} is justified since the above expressions indeed correspond to integration with respect to (a family of) Keisler measures on the type space. More concretely, for any $\bar{a}\in F^{|\bar{y}|}$ the function\footnote{Here, $\mathcal{C}_{c}(S_{\bar{x}}(F),\mathbb{C})$ is the space of continuous compactly supported functions $g:S_{\bar{x}}(F)\rightarrow\mathbb{C}$, that is, the $\mathbb{C}$-vector space of definable predicates in the variables $\bar{x}$.} $\mathrm{Int}_{B_{\bar{a}}}:\mathcal{C}_{c}(S_{\bar{x}}(F),\mathbb{C})\rightarrow\mathbb{C}$ given by
    \[P(\bar{x},\bar{a})\mapsto\int_{B_{\bar{a}}}P(\bar{x},\bar{a})d\mu_{B_{\bar{a}}}^{*}(\bar{x})\]
    corresponds by the Riesz-Markov-Kakutani theorem to a Keisler measure on $S_{\bar{x}}(F)$. See Section 2.3 of \cite{chavarria2021continuous} for details; in their language $\mathrm{Int}_{B_{\bar{a}}}$ corresponds to a \textit{Keisler functional}.
\end{remark}

From Fact \ref{factzoemasureiscountingmeasure} as well as the fact that the corresponding statement trivially holds in finite fields (as we work with the counting measure) we obtain the following basic properties, which we state here as they will be implicitly used later on.

\begin{fact}\label{factfubini}
For any $\bar{a}\in F^{|\bar{y}|}$ the following holds.
\begin{itemize}
	\item For every $\mathcal{L}_{\mathrm{ring}}$-definable subset $D\subseteq B_{\bar{a}}$ we have \[\mu_{B_{\bar{a}}}^{*}(D)=\int_{B_{\bar{a}}}\mathbbm{1}_{D}(\bar{x})d\mu_{B_{\bar{a}}}^{*}=\mu_{B_{\bar{a}}}(D).\]
	% \item Every $\mathcal{L}_{+,\times}$-definable (with parameters) predicate $P(\bar{x})$ is measurable with respect to $\mu_{B_{\bar{a}}}^{*}$.
	\item The following version of \textit{Fubini} holds. Assume that $\bar{x}=(\bar{x}_{1},\bar{x}_{2})$ and let $\pi_{1}, \pi_{2}$ be the corresponding projections. 
%	Assume that all the fibres $\pi_{1}(\bar{y})$ (for $\bar{y}\in\pi_{1}(B_{\bar{a}})$) have the same dimension and absolute measure (and similarly for $\pi_{2}$). 
For any $\bar{c}\in \pi_{2}(B_{\bar{a}})$, we write $B_{\bar{a}}^{\bar{c}}$ for the set $\pi_{1}(\pi_{2}^{-1}(\bar{c}))$. Assume that $B_{\bar{a}}^{\bar{c}}$ has the same dimension and multiplicity for all $\bar{c}\in \pi_{2}(B_{\bar{a}})$.
Then, for any $\mathcal{L}_{+,\times}$-definable predicate $P(\bar{x})$ the following holds.
	\[ \int_{B_{\bar{a}}}P(\bar{x})d\mu^{*}_{B_{\bar{a}}}(\bar{x})=\int_{\pi_{2}(B_{\bar{a}})} \left(\int_{B_{\bar{a}}^{\bar{x}_{2}}} P(\bar{x}_{1},\bar{x}_{2})d\mu^{*}_{B_{\bar{a}}^{\bar{x}_{2}}}(\bar{x}_{1}) \right)d\mu^{*}_{\pi_{2}(B_{\bar{a}})}(\bar{x}_{2}).
	\]
\end{itemize}
\end{fact}

\begin{convention}
From now on we will simply write $\mu_{B_{\bar{y}}}$ instead of $\mu_{B_{\bar{y}}}^{*}$. When we refer to definability of the measure, we mean definability (in the parameters) of the restriction of $\mu_{B_{\bar{y}}}$ to definable sets as given by Fact \ref{factsizeofdebsets}.
\end{convention}

The following is a generalisation of Proposition 4.1 for $\mathrm{PF}^{+}$ in \cite{Hrushovski2021AxsTW} to the theory $\mathrm{PF}^{+,\times}$. Note that it is the non-standard version of Theorem \ref{theoremdefinableintegrationfiniteversion}. In particular, once we obtain Theorem \ref{theoremdefinableintegrability}, its finite version, Theorem \ref{theoremdefinableintegrationfiniteversion} immediately follows. 

\begin{theorem}\label{theoremdefinableintegrability}
 We work in $F$ (a characteristic $0$ ultraproduct of finite fields). Let $P(\bar{x},\bar{y})$ be an $\mathcal{L}_{+,\times}$-definable predicate and let $B_{\bar{y}}$ be an $\mathcal{L}_{\mathrm{ring}}$-definable set. Then the following expression $\mathrm{Int}_{B_{\bar{y}}}(P)(\bar{y})$ is an $\mathcal{L}_{+,\times}$-definable predicate (in $\bar{y}$).
    \[\mathrm{Int}_{B_{\bar{y}}}(P)(\bar{y})=\int_{B_{\bar{y}}}P(\bar{x},\bar{y})d\mu_{B_{\bar{y}}}(\bar{x}).\]
\end{theorem}

Fact \ref{factdecompintersect} and Fact \ref{factkowlaskireductiontoalgsets} will now allow us to prove the following theorem on (uniform) definability of character sums over definable subsets of algebraic curves that together with our results on quantifier elimination in $\mathcal{L}_{\kappa,+,\times}^{\mathrm{sym}}$ (and reduction to basic predicates) will be the main ingredients to prove Theorem \ref{theoremdefinableintegrability}.
\begin{lemma}\label{theoremdefinablecharactersums}
    Let $C_{\bar{y}}$ be a family of definable sets of dimension $1$ where $\bar{y}$ varies over some definable set $D$. Assume that $B_{\bar{a}}$ is a definable subset of $C_{\bar{a}}$ for any $\bar{a}\in D$. Let $h(\bar{x})$ and $g(\bar{x})$ be an integral linear (resp. multiplicative) transformation. Then the following expression $\mathrm{Int}_{h,g}^{BC}(\bar{y})$ is an $\mathcal{L}_{+,\times}$-definable predicate (in the variables $\bar{y}$) on $D$.
    \[\mathrm{Int}_{h,g}^{BC}(\bar{y}):=\int_{B_{\bar{y}}}\Psi(h(\bar{x}))\chi(g(\bar{x}))d\mu_{C_{\bar{y}}}(\bar{x}).\]
\end{lemma}

\begin{proof}
First we note that
\[\int_{B_{\bar{y}}}\Psi(h(\bar{x}))\chi(g(\bar{x}))d\mu_{C_{\bar{y}}}(\bar{x})=\mu_{C_{\bar{y}}}(B_{\bar{y}})\cdot\int_{B_{\bar{y}}}\Psi(h(\bar{x}))\chi(g(\bar{x}))d\mu_{B_{\bar{y}}}(\bar{x})
\]
We now apply Fact \ref{factkowlaskireductiontoalgsets} with $\beta(\bar{x}):=\Psi(h(\bar{x}))\chi(g(\bar{x}))$ to obtain a partition of $B_{\bar{y}}$ into $\mathcal{L}_{\mathrm{ring}}$-definable sets $B_{\bar{y}}^{1},\dots,B_{\bar{y}}^{l}$ such that for any $1\leq i\leq l$ the following holds. There is a finite tuple $\bar{b}$ in $F^{|\bar{u}|}$ and  algebraic sets $W_{\bar{y}\bar{u}}^{i,j}$ (in the variables $\bar{z}$ with $|\bar{z}|\geq|\bar{x}|$) with $1\leq j\leq e$ such that for all $\bar{a}\in D$ and any $1\leq i\leq l$ we have that $\int_{B_{\bar{a}}^{i}}\Psi(h(\bar{x}))\chi(g(\bar{x}))d\mu_{B_{\bar{y}}}(\bar{x})$ is given by
\[\sum_{1\leq j\leq e}\mu_{B_{\bar{a}}^{i}}(\pi(W_{\bar{a}\bar{b}}^{i,j}))\cdot\left(    \int_{W_{\bar{a}\bar{b}}^{i,j}}\Psi(h(\pi(\bar{z})))\chi(g(\pi(\bar{z})))d\mu_{W_{\bar{a}\bar{b}}^{i,j}}(\bar{z})    \right)
\]where $\pi$ is the projection to the first $|\bar{x}|$-many variables.
It thus suffices to 
prove that an expression of the form
\[\mathrm{Int}_{h,g}^{J}:=\int_{J_{\bar{y}}}\Psi(h(\bar{x}))\chi(g(\bar{x}))d\mu_{J_{\bar{y}}}(\bar{x})\]
is uniformly definable in $\bar{y}$ where $J_{\bar{y}}$ is an algebraic curve (not necessarily irreducible).\\ Moreover, if necessary, we can further partition the parameter space into $\mathcal{L}_{\mathrm{ring}}$-definable sets as to assume that for all $\bar{a},\bar{a}^{\prime}\in D$ we have $\mu(J_{\bar{a}})=\mu(J_{\bar{a}^{\prime}})$ (and $\mathrm{dim}(J_{\bar{a}})=1=\mathrm{dim}(J_{\bar{a}^{\prime}})$) for $\mu$ the multiplicity (which in this case is given by the number of irreducible components of $J_{\bar{a}}$ (resp. $J_{\bar{a}^{\prime}}$) of dimension $1$ defined over $F$.\\
If $J_{\bar{y}}$ is already absolutely irreducible we only need to observe that it is (uniformly) definable in $\bar{y}$ whether $h$ and $g$ are constant on $J_{\bar{y}}(F)$. Whenever it is the case that they are constant for some $\bar{a}$ we have $\mathrm{Int}_{h,g}^{J}(\bar{a})=\Psi(h(\bar{c}))\chi(g(\bar{c}))$ for any $\bar{c}\in J_{\bar{a}}(F)$. Otherwise, we can apply the Weil bounds (Fact \ref{factperelmutterestimates}) and obtain that $\mathrm{Int}_{h,g}^{J}(\bar{a})=0$.\\
To reduce from general curves to absolutely irreducible ones we can use Fact \ref{factdecompintersect} which yields uniform definability (in $\bar{y}$) of the decomposition-intersection procedure applied on $J_{\bar{y}}$ (and we note that lower-dimensional algebraic subsets have measure $0$ and thus do not have to be taken care of, i.e., it suffices to compute the integral for every absolutely irreducible component).
\end{proof}
% \begin{remark}
% Note that in the above Theorem we only prove definability in the parameters of the character sum/integral over a definable set but do not give a concrete description of the corresponding formula. 
% %I.e., it is uniformly definable in $\bar{y}$ whether $h$ and $g$ are constant on the absolutely irreducible curves obtained from $B_{\bar{y}}$ from Fact \ref{factkowlaskireductiontoalgsets} and the decomposition-intersection procedure but do not translate this in an \textit{easy} statement on $h,g$ and $W_{\bar{y}}$.
% See for example the discussion of Example 4 in \cite{Kowalski2005ExponentialSO} to see that this is in general less straightforward in the multiplicative than in the additive case.\footnote{Note however that in the example given there the multiplicative character $\chi$ is of order two, whereas by our assumption the order of $\chi$ tends to $\infty$.}
% \end{remark}
We can now proceed to give the proof of Theorem \ref{theoremdefinableintegrability}. It follows the same ideas as the proof of the purely additive case (Proposition 4.1 in \cite{Hrushovski2021AxsTW}).

\begin{proof}(Proof of Theorem \ref{theoremdefinableintegrability}.)
Write $\bar{x}=(x_{1},\bar{x}_{2})$. Let $\pi_{1}$ be the projection to the first variable (of $\bar{x})$ and $\pi_{2}$ to the others. We partition $B_{\bar{y}}$ into definable subsets $B_{\bar{y},1},\dots,B_{\bar{y},k}$ such that for any fixed $1\leq i\leq k$ the definable sets $B_{\bar{y},i}^{\bar{c}}$ have all the same dimension and multiplicity where $\bar{c}$ varies over $\pi_{2}(B_{\bar{y},i})$. From the equality
\[\int_{B_{\bar{y}}}P(\bar{x},\bar{y})d\mu_{B_{\bar{y}}}(\bar{x})=\sum_{1\leq i\leq k}\mu_{B_{\bar{y}}}(B_{\bar{y},i})\cdot\int_{B_{\bar{y},i}}P(\bar{x},\bar{y})d\mu_{B_{\bar{y},i}}(\bar{x})
\]it follows that it suffices to show that $\int_{B_{\bar{y},i}}P(\bar{x},\bar{y})d\mu_{B_{\bar{y},i}}(\bar{x})$ is definable in $\bar{y}$. We can apply Fubini (Fact \ref{factfubini}) and obtain that $\int_{B_{\bar{y},i}}P(\bar{x},\bar{y})d\mu_{B_{\bar{y},i}}(\bar{x})$ is given by the expression
\[\int_{\pi_{2}(B_{\bar{y},i})}\left( \int_{B_{\bar{y},i}^{\bar{x}_{2}}}P(\bar{x},\bar{y})d\mu_{B_{\bar{y},i}^{\bar{x}_{2}}}(x_{1})  \right)d\mu_{\pi_{2}(B_{\bar{y}})}(\bar{x}_{2})
\]
By induction it then suffices to show that the expression
\[\int_{B_{\bar{y},i}^{\bar{x}_{2}}}P(\bar{x},\bar{y})d\mu_{B_{\bar{y},i}^{\bar{x}_{2}}}(x_{1})\]
is a definable predicate in the variables $(\bar{y},\bar{x}_{2})=(\bar{y},x_{2},\dots,x_{n})$. By construction of $B_{\bar{y},i}$ we have that $B_{\bar{a},i}^{\bar{c}}$ is of the same dimension for all tuples $(\bar{a},\bar{c})$ that satisfy $\bar{a}\in D$ and $\bar{c}\in \pi_{2}(B_{\bar{a},i})$.\\ 
As $B_{\bar{a}}$ has dimension 1, the dimension of all the $B_{\bar{a},i}^{\bar{c}}$ is either $0$ or $1$. If it is $0$, then the sets are finite (and bounded uniformly) and definability of the integral expression (which is then just a finite sum) is clear. Let us thus assume that the dimension is $1$.\\
To reduce the notation we summarise the above as follows: We have reduced the problem to showing that an expression of the form $\int_{B_{\bar{y},i}}P(x,\bar{y})d\mu_{B_{\bar{y},i}}(x)$ is definable in $\bar{y}$ where for all $\bar{a}$ we have that $B_{\bar{a},i}$ is of dimension $1$. For the sake of readability we will now suppress the $i$ in the notation and write $B_{\bar{y}}$ instead.\\
Recall from Definition \ref{defintionlanguageforqepfplustimes} that for some $\mathcal{L}_{\mathrm{ring}}$-formula $\phi(\bar{s},\bar{v})$ and $g,h$ integral linear (resp. multiplicative) transformations we defined $\Theta_{\mathrm{sym}}^{\phi,g,h}(\bar{v})=\sum_{\bar{r}\in X_{\bar{v}}^{\phi}}\Psi(g(\bar{r}))\chi(h(\bar{r}))$ where for $\bar{e}\in F^{|\bar{v}|}$ we denote by $X_{\bar{e}}^{\phi}$ the set defined by $\phi(\bar{s},\bar{e})$. By Theorem \ref{qepfplustimes} and Lemma \ref{lemmaalgebraofbasicpredicatespfplustimes} we can assume that $P(x,\bar{y})$ is of the form $\Theta_{\mathrm{sym}}^{\phi,g,h}(\kappa_{1}(x,\bar{y}),\dots,\kappa_{n}(x,\bar{y}))\cdot \mathbbm{1}_{D_{\bar{y}}}$ as in Definition \ref{defintionlanguageforqepfplustimes} where the $\kappa_{i}(x,\bar{y})$ are $\mathcal{L}_{\kappa,+,\times}^{\mathrm{sym}}$-terms and $D_{\bar{y}}$ is an $\mathcal{L}_{\mathrm{ring}}$-definable set in the variable $x$. Hence, we have to show that the expression
\[\int_{D_{\bar{y}}}\Theta_{\mathrm{sym}}^{\phi,g,h}(\kappa_{1}(x,\bar{y}),\dots,\kappa_{n}(x,\bar{y}))d\mu_{B_{\bar{y}}}(x)\]is uniformly definable in $\bar{y}$. The expression $\kappa_{i}(x,\bar{y})$ is uniformly algebraic\footnote{That is, $e=\kappa_{i}(d,\bar{a})\in acl(d,\bar{a})$, uniformly for all $(d,\bar{a})$.} in $(x,\bar{y})$ and 
similarly the expression $\bar{r}\in X_{\bar{v}}^{\phi}$ is uniformly algebraic in $\bar{v}$. Thus, for any $\bar{a}$ we have that the set
 \[C_{\bar{a}}:=\{(x,\bar{v},\bar{r})\,|\,x\in D_{\bar{a}},\kappa_{1}(x,\bar{a})=v_{1},\dots,\kappa_{n}(x,\bar{a})=v_{n},\bar{r}\in X_{\bar{v}}^{\phi}\}\]
 is a definable set of dimension $1$. We can again partition $B_{\bar{y}}$ into definable subsets $B_{\bar{y},1},\dots,B_{\bar{y},k}$ such that for each such subset $B_{\bar{y},i}$ and all $\bar{a}\in F^{|\bar{y}|}$ all the fibres of the projection to the first variable $\pi_{1}:C_{\bar{a}}\rightarrow B_{\bar{a},i}$ have the same (finite) size. We obtain $\int_{D_{\bar{y}}}\Theta_{\mathrm{sym}}^{\phi,g,h}(\kappa_{1}(x,\bar{y}),\dots,\kappa_{n}(x,\bar{y}))d\mu_{B_{\bar{y}}}(\bar{x})$ is given by 
 \[\sum_{1\leq i\leq k}\mu_{B_{\bar{y}}}(B_{\bar{y},i})\cdot\int_{D_{\bar{y}}}\Theta_{\mathrm{sym}}^{\phi,g,h}(\kappa_{1}(x,\bar{y}),\dots,\kappa_{n}(x,\bar{y}))d\mu_{B_{\bar{y},i}}(x).
\]Moreover we have that $\int_{D_{\bar{y}}}\Theta_{\mathrm{sym}}^{\phi,g,h}(\kappa_{1}(x,\bar{y}),\dots,\kappa_{n}(x,\bar{y}))d\mu_{B_{\bar{y},i}}(x)$ is given by \[\mu_{B_{\bar{y},i}}(D_{\bar{y}})\cdot\int_{C_{\bar{y}}}\Psi(\Tilde{g}(\bar{w}))\chi(\Tilde{h}(\bar{w}))d\mu_{C_{\bar{y}}}(\bar{w})\]where we set $\bar{w}=(x,\bar{v},\bar{r})$ and $\Tilde{g}(\bar{w})=\Tilde{g}(x,\bar{v},\bar{r}):=g(\bar{r})$ (and similarly for $\Tilde{h}$). But the latter expression is uniformly definable in $\bar{y}$ by Lemma \ref{theoremdefinablecharactersums}.
\end{proof}

\begin{remark}\label{remarkKeislermeasure}
In light of Fact \ref{factzoemasureiscountingmeasure} we could ask whether we have a measure-theoretic interpretation of $\mu_{B_{\bar{y}}}$ for general models $F\models \mathrm{PF}^{+,\times}$ as well: By Theorem \ref{theoremdefinableintegrability} $\mathrm{Int}_{B_{\bar{y}}}$ is definable. Then, for any $\bar{a}\in F^{|\bar{y}|}$, the corresponding definable predicate defines a positive linear functional on $\mathcal{C}_{c}(S_{\bar{x}}(F),\mathbb{C})$ and hence, yet again by Riesz-Markov-Kakutani (see Fact 2.13 in \cite{chavarria2021continuous}) $\mu_{B_{\bar{a}}}$ yields a Keisler measure on $S_{\bar{x}}(F)$.
\end{remark}

We will now briefly explain how we can easily strengthen Theorem \ref{theoremdefinableintegrationfiniteversion} and remove the lower bounds $k_{0}$ and $r_{0}$ on the order of $\Psi$ and $\chi$. The strong version of the main theorem of this section then reads as follows.
\begin{theorem}\label{corollaryremovinglowerbounds}
        Let $P(\bar{x},\bar{y})$ be an $\mathcal{L}_{+,\times}$-predicate and $B_{\bar{y}}$ an $\mathcal{L}_{\mathrm{ring}}$-definable set. For any finite field $\mathbb{F}_{q}$ where $q=p^{e}$ consider the function $f:\mathbb{F}_{q}^{|\bar{y}|}\rightarrow\mathbb{C}$ given by $f(\bar{y})=0$, if $B_{\bar{y}}=\emptyset$, and otherwise by\[f(\bar{y}):=\frac{\sum_{\bar{x}\in B_{\bar{y}}(\mathbb{F}_{q})}P(\bar{x},\bar{y})}{|B_{\bar{y}}(\mathbb{F}_{q})|}.\] 
Now for any $\epsilon>0$, there exists an $\mathcal{L}_{+,\times}$-predicate $\Tilde{P}(\bar{y})$  such that for any triple $(\mathbb{F}_{q},\Psi_{q},\chi_{q})$ (consisting of a finite field with additive and multiplicative character) and any $\bar{a}\in \mathbb{F}_{q}^{|\bar{y}|}$,
\[|\Tilde{P}(\bar{a})-f(\bar{a})|\leq\epsilon.\]
\end{theorem}

\begin{proof} We show how we can get rid of the lower bound $r_{0}$, the additive case for $k_{0}$ then follows analogously. We first observe that $\mathrm{ord}(\chi)\geq r_{0}$ is expressible by an $\mathcal{L}_{+,\times}$-sentence. Thus, a case distinction becomes definable. As there are only finitely many possibilities for $\mathrm{ord}(\chi)<r_{0}$, it suffices to treat any of the cases $\mathrm{ord}(\chi)=\Tilde{r}_{0}$ for some given $\Tilde{r}_{0}<r_{0}$ separately. Let $I$ be an enumeration of the finitely many occurrences of $\chi$ in $P(\bar{x},\bar{y})$. We now partition $B_{\bar{y}}$ into $|I|\cdot\Tilde{r}_{0}$ many sets $D_{i}$ such that to each of those subsets corresponds exactly one $|I|$-tuple of $\Tilde{r}_{0}$-th roots of unity $(\lambda_{1},\dots,\lambda_{|I|})$ such that $(\bar{a},\bar{b})\in D_{i}$, if for all $1\leq j\leq |I|$ the $j$-th occurrence of $\chi$ in $P(\bar{a},\bar{b})$ takes value $\lambda_{j}$ (and $(\bar{a},\bar{b})\in B_{\bar{y}}(\mathbb{F}_{q})$). We then associate to $D_{i}$ the $\mathcal{L}_{+}$-formula $P_{i}$ obtained by replacing for any $1\leq j\leq |I|$ the $j$-th occurrence of $\chi$ in $P(\bar{x},\bar{y})$ by $\lambda_{j}$. We then have
\[f(\bar{y})=\frac{\sum_{1\leq i\leq |I|\Tilde{r}_{0}}\sum_{\bar{x}\in D_{i,\bar{y}}(\mathbb{F}_{q})}P_{j}(\bar{x},\bar{y})}{|B_{\bar{y}}(\mathbb{F}_{q})|}.\]
If we can now show that any of the $D_{j}$ is already $\mathcal{L}_{\mathrm{ring}}$-definable and thus $\frac{|D_{j,\bar{y}}(\mathbb{F}_{q})|}{|B_{\bar{y}}(\mathbb{F}_{q})|}$ is uniformly definable (in $\bar{y}$), then we have reduced the problem (for the case $\mathrm{ord}(\chi)=\Tilde{r}_{0}$) to one in $\mathcal{L}_{+}$.\\
But to see that $D_{j}$ is indeed $\mathcal{L}_{\mathrm{ring}}$-definable it suffices to observe that $\chi^{-1}(\lambda_{j})=\{x\in\mathbb{F}_{q}\,|\,\exists y\;y^{\Tilde{r}_{0}}=\frac{x}{a}\}$ where $a\in\mathbb{F}_{q}$ is such that $\chi(a)=\lambda_{j}$. For the additive case we would use the condition $\exists y\;y^{q}-y=x-a$ instead in case that $\Psi_{q}$ is the standard character and otherwise we reduce to this case by replacing above $y^{q}-y$ by $c\cdot y^{q}-y$ for some $c\in\mathbb{F}{q}$ using Fact \ref{remarkdescriptionalladdcharacters}.
\end{proof}

\section{Further remarks}\label{sectionfurtheremarkspfplustimes}
The following three results follow exactly as for $\mathrm{PF}^{+}$ in \cite{Hrushovski2021AxsTW}. For the sake of completeness we outline the respective arguments but omit most of the details.\\

\textbf{Decidability} holds for $\mathrm{PF}^{+,\times}$ in the same sense as it does for $\mathrm{PF}^{+,\times}$. The argument given in 3.14 of \cite{Hrushovski2021AxsTW} transfers directly to our context: Given a sentence $\phi$ and $\epsilon>0$, we can search a quantifier free sentence $\phi^{\prime}$ and a proof for $|\phi-\phi^{\prime}|<\frac{\epsilon}{2}$. Next, $\phi^{\prime}$ is determined by the $\mathcal{L}_{\mathrm{ring}}$-type of a finite Galois extension of $\mathbb{Q}$ together with the values of $\Psi$ and $\chi$ on a finite $\mathbb{Q}$-linearly independent set as well as on a finite multiplicatively independent set. We can search for $\Tilde{\phi}$ that describes the values of those sets. It then follows, that by determining the value of $\Tilde{\phi}$ up to $\frac{\epsilon}{2}$, we can determine to value of $\phi$ up to $\epsilon$.\\

\textbf{Definable sets.} Recall that a definable set in a continuous logic structure is a set for which the distance to it is given by a definable predicate. In our setting this amounts to: $D\subseteq F^{n}$ is $\mathcal{L}_{+,\times}$-definable if there is a definable $\mathcal{L}_{+,\times}$-predicate $\phi(\bar{x})$ such that $\phi(D)=c\neq c^{\prime}=\phi(D^{\mathrm{c}})$. However, there are no new definable sets in $\mathrm{PF}^{+,\times}$ compared to $\mathrm{PF}$. Note that this does by far not hold for type-definable ones, e.g. $\Psi(x)=1$ describes a type-definable set that is not definable in $\mathrm{PF}$ (see Section 2 of \cite{Hrushovski2021AxsTW} for an explanation why the characters are not definable in a model-theoretically tame expansion of $\mathrm{PF}$ in discrete logic).
\begin{theorem}
    Any definable set in $F\models\mathrm{PF}^{+,\times}$ is already definable in its $\mathcal{L}_{\mathrm{ring}}$-reduct (which is a model of $\mathrm{PF}$.
\end{theorem}

\begin{proof}
First, one shows that for an $\mathcal{L}_{\mathrm{ring}}$-definable set $D\subseteq F^{n}$, the set $\Psi^{(n)}\times\chi^{(k)}(D)$ is a finite union of cosets of subtori of $\mathbb{T}^{n+k}$ where the preimages of different cosets are $\mathcal{L}_{\mathrm{ring}}$-definable subsets. This follows exactly as in the proof of Proposition 4.3 in \cite{Hrushovski2021AxsTW}. Next, by Theorem \ref{qepfplustimes} (quantifier elimination) and Lemma \ref{lemmaalgebraofbasicpredicatespfplustimes} (reduction to basic predicates) the problem reduces to consider the image of definable sets under $\Psi^{(n)}\times\chi^{(k)}$. Using the above it is then easy to see that the preimages (under our definable predicate) of connected components (in $\mathbb{C}$) are $\mathcal{L}_{\mathrm{ring}}$-definable.
%the sets $\Psi^{(n)}(D)$ and $\chi^{(n)}(D)$ are finite unions of subtori of $\mathbb{T}^{n}$ {\color{red}THIS IS NOT ENOUGH, CONSIDER UNION}. This follows exactly as in the proof of Proposition 4.3 in \cite{Hrushovski2021AxsTW}. Next, by Theorem \ref{qepfplustimes} (quantifier elimination) and Lemma \ref{lemmaalgebraofbasicpredicatespfplustimes} (reduction to basic predicates) it reduces to {\color{red}passing to finite cover as in \cite{Hrushovski2021AxsTW} is not entirely trivial but as everything remains definable, it should work.}
\end{proof}

\textbf{Simplicity.} Recall that forking-independence in pseudofinite fields coincides with algebraic independence for (algebraically closed) substructures. That is, we say that for (relatively) algebraically closed fields $C\subseteq A,B\subseteq F\models\mathrm{PF}$, $A$ is independent from $B$ over $C$ (written $A\ind_{C} B$) if $A$ is algebraically independent from $B$ over $C$. We will show simplicity for $\mathrm{PF}^{+,\times}$ using essentially the same proof as for Proposition 3.20 in \cite{Hrushovski2021AxsTW}.

\begin{theorem}
    The theory $\mathrm{PF}^{+,\times}$ is simple.
\end{theorem}
\begin{proof}
We will use the characterisation of simple theories by the existence of a good notion of independence which also works in the continuous logic setting (Theorem 1.51 in \cite{simplicityCATs}). We consider the same notion of independence in $\mathrm{PF}^{+,\times}$ as for $\mathrm{PF}$. As it directly inherits all the other necessary properties from $\mathrm{PF}$, it suffices to verify that the Independence theorem, or equivalently 3-amalgamation, holds over models in $\mathrm{PF}^{+,\times}$. We work in a sufficiently saturated model $F\models \mathrm{PF}^{+,\times}$ and have to show the following:
    Let $\mathcal{V}=\{p_{1}(\bar{x}_{1}),p_{2}(\bar{x}_{2}),p_{3}(\bar{x}_{3}),p_{12}(\bar{x}_{12}),p_{13}(\bar{x}_{13}),p_{23}(\bar{x}_{23})\}$ be a system of complete types over some $\mathcal{M}\preceq$ such that
    \begin{itemize}
        \item for any $\bar{a}_{w}\models p_{w}$ the infinite tuple $\bar{a}_{w}$ enumerates its algebraic closure. 
        \item For any $w\subseteq w^{\prime}$ we have $\bar{x}_{w}\subseteq \bar{x}_{w^{\prime}}$ and $p_{w^{\prime}}\restriction_{w}=p_{w}$.
        \item For any $\bar{a}_{ij}\models p_{ij}$ and corresponding subtuples $\bar{a}_{i}$ and $\bar{a}_{j}$, we have $\bar{a}_{i}\ind_{\mathcal{M}}\bar{a}_{j}$.
    \end{itemize}
    Then, there is some $p_{123}(\bar{x}_{123})$ such that the system $(p_{w})_{w\subseteq \{1,2,3\}}$ still satisfies the above properties (in particular for $\bar{a}_{123}\models p_{123}$, we require $\bar{a}_{ij}\ind_{\bar{a}_{j}}\bar{a}_{jk})$ for $\{i,j,k\}=\{1,2,3\}$.
    To prove this, we first consider the 3-amalgamation problem (in $\mathrm{PF}$) given by the $\mathcal{L}_{\mathrm{ring}}$-reducts of our initial system. By simplicity of $\mathrm{PF}$ we find an $\mathcal{L}_{\mathrm{ring}}$-type $\Tilde{p}_{123}$ that is a solution. Given a realisation $\Tilde{a}\models\Tilde{p}_{123}$ we thus have to show that we can extend $\Psi$ and $\chi$ to additive and multiplicative characters consistent with all the types in the system $\mathcal{V}$.\\
    It is sufficient to show that we can do this for $\Psi$ on the $\mathbb{Q}$-vector space  $\langle\Tilde{a}_{12},\Tilde{a}_{13},\Tilde{a}_{23}\rangle_{\mathbb{Q}}$ as well as for $\chi$ on the multiplicative subgroup $\langle\Tilde{a}_{12},\Tilde{a}_{13},\Tilde{a}_{23}\rangle$, 
    %generated by all elements projecting to torsion elements in $\Tilde{a}_{123}^{\times}/\langle{a}_{12},\Tilde{a}_{13},\Tilde{a}_{23}\rangle$,
    because outside of these groups we can simply choose an extension of $\Psi$ (resp. $\chi$) to homomorphisms on $(\Tilde{a}_{123},+)$ and $(\Tilde{a}_{123}^{\times},\cdot)$.
    By compactness this reduces to showing that for any $b_{ij},c_{ij}\in\Tilde{a}_{ij}$ with $b_{12}+b_{13}+b_{23}=0$ and $c_{12}\cdot c_{13}\cdot c_{23}=1$ we have $\sum\Psi(b_{ij})=1$ and $\prod\chi(c_{ij})=1$ where $\Psi(b_{ij})$ and $\chi(c_{ij})$ are as specified by $p_{ij}$.
    Let $p_{ij}^{\mathrm{acl}}$ be the quantifier-free $\mathcal{L}_{ring}$-type (over $M^{\mathrm{acl}}$) of the (field-theoretic) algebraic closure $\Tilde{a}_{ij}^{\mathrm{acl}}$ of $\Tilde{a}_{ij}$.\\
    Then, $\mathrm{tp}^{\mathrm{qf}}_{\mathcal{L}_{\mathrm{ring}}}(\Tilde{a}_{12}^{\mathrm{acl}}/\Tilde{a}_{3}^{\mathrm{acl}})$ is the co-heir of $p_{12}^{\mathrm{acl}}$ and thus there are $\Tilde{b}_{i},\Tilde{c}_{i}\in\Tilde{a}_{i}$ as well as $\Tilde{d}_{i},\Tilde{e}_{i}\in M^{\mathrm{acl}}$ such that $b_{ij}=\Tilde{b}_{i}+\Tilde{d}_{i}+\Tilde{b}_{j}+\Tilde{d}_{j}$ and similarly $c_{ij}=\Tilde{c}_{i}\cdot\Tilde{e}_{i}\cdot \Tilde{c}_{j}\cdot\Tilde{e}_{j}$ (see e.g. 1.9 in \cite{acfa} for a similar argument). Note that $b_{j}=\Tilde{d}_{i}+\Tilde{b}_{j}+\Tilde{d}_{j}\in\Tilde{a}_{j}$ as $b_{ij}$ and $b_{i}:=\Tilde{b}_{i}$ are both in $F$ and similarly for $c_{j}$.
    But this already completes the proof as then $\Psi(b_{ij})$ and $\chi(c_{ij})$ are determined by $\Psi(b_{i})+\Psi(b_{j})$ and $\chi(c_{i})\cdot\chi(c_{j})$ and thus by the types $p_{i}$ and $p_{j}$. As this hold for all pairs $(i,j)$ with $1\leq i\neq j\leq 3$ (after possibly adding or multiplying with elements from $M$ we can assume $b_{i},c_{i}$ to be the same when occurring in different pairs), it follows that $\sum\Psi(b_{ij})=1$ and $\prod\chi(c_{ij})=1$ which we had to show. 
    %But otherwise there is $\tau\in\mathrm{Gal}(\Tilde{a}_{ij}^{\mathrm{acl}}/\Tilde{a}_{ij})$ such that $\tau(\Tilde{b}_{i})-\Tilde{b}_{i}\neq 0$. {\color{red}We need fullness?With existentially closedness for pseudofinite fields equipped with Galois generator this works, because we get a solution of $\tau(x)-x=d_{i}$ in $\mathcal{M}$?With elimination of imaginaries we should be able to code it?}
%    But then $\tau(\Tilde{b}_{i})-\Tilde{b}_{i}=\tau(\Tilde{b}_{j})-\Tilde{b}_{j}$ and by independence of $\Tilde{a}_{i}$ and $\Tilde{a}_{j}$ thus $\tau(\Tilde{b}_{i})-\Tilde{b}_{i}=d_{i}\in M^{\mathrm{acl}}$.
\end{proof}

\begin{remark}
    The above theorem generalises straight-forwardly to higher dimensions. In particular, in $\mathrm{PF}^{+,\times}$ $n$-amalgamation holds for all $n\in\mathbb{N}$ over $\mathrm{acl}^{\mathrm{eq}}$-closed substructures (where $\mathrm{acl}^{\mathrm{eq}}$ is meant as for $\mathrm{PF}$, i.e., a priori without taking into account continuous-logic-imaginaries).
\end{remark}

\textbf{The purely multiplicative case.} All the results on the theory $\mathrm{PF}^{+,\times}$ that we obtained in the last sections hold similarly for the theory $\mathrm{PF}^{\times}$ of pseudofinite fields with a multiplicative character of infinite order using the same proofs. We refrain from restating all the results in this context. However, as some of the difficulties from the mixed context disappear, we give an axiomatisation as well as a language for quantifier elimination that are slightly less technical than the ones that immediately fall of from the $\mathrm{PF}^{+,\times}$. They are closer to the purely additive case in \cite{Hrushovski2021AxsTW}. As the proofs remain essentially the same, instead of giving them here we redirect the reader to Chapter 5 of the author's PhD thesis \cite{ludwig:tel-05236078} where the purely multiplicative case is treated on the way.

\begin{definition}
    Let $F$ be a field. We say that a multiplicative coset defined by an equation of the form $\prod_{1\leq i\leq n} X_{i}^{a_{i}}=b$ for $(a_{1},\dots,a_{n})\in\mathbb{Z}^{n}\backslash\{(0,\dots,0)\}$ and $b\in F^{\times}$ is of height $\leq k$ for $k\in\mathbb{N}$, if it can be defined with $|a_{i}|\leq k$ for all $1\leq i\leq n$.
\end{definition}

\begin{definition}\label{definitionpftimes}
    We define $F$ to be a model of the $\mathcal{L}_{\times}$-theory $\mathrm{PF}^{\times}$ if it satisfies the following axiom scheme: 
    \begin{enumerate}[(1)]
        \item $\mathrm{char}(F)=0$ and $\mathrm{Gal}(F)=\hat{\mathbb{Z}}$.
        \item $\chi: (F^{\times},\cdot)\rightarrow S^{1}$ is a group homomorphism and $\chi(0)=0$.
        \item Let $h\in\mathbb{Q}[Z,Z^{-1},\dots,Z_{n},Z_{n}^{-1}]$ be a finite Laurent polynomial with degree $\leq m$, real valued on $\mathbb{T}^{n}=S^{1}\times\dots\times S^{1}$ with no constant term. For any non zero-degenerate absolutely irreducible curve $C\subset\mathbb{A}^{n}$ defined over $F$ not contained in any multiplicative coset of height at most $m$ the following holds:
        \[\sup\{h(\chi^{n}(\bar{x}))\,|\,\bar{x}\in C^{\prime}(F)\}\geq 0.\]
    \end{enumerate}
\end{definition}

% \begin{lemma}
%     Let $F$ be a field of characteristic $0$ and $\chi:F^{\times}\rightarrow S^{1}$ a homomorphism. Then $F$ satisifies axiom $(3)$ of $\mathrm{PF}^{\times}$ (Definition \ref{definitionpftimes}) if and only if the following holds:\\
%     For any non zero-degenerate, absolutely irreducible curve $C\subset\mathbb{A}^{n}$ over $F$, not contained in any multiplicative coset, $\chi^{n}(C(F))$ is dense in $\mathbb{T}^{n}$.
% \end{lemma}

% \begin{proof}
%     See the proof of the more general Lemma \ref{definitionporpertyequivtoaxiompfplustimes}.
% \end{proof}

% \begin{theorem}
%     Let $\chi_{q}$ denote a multiplicative character on the finite field $\mathbb{F}_{q}$. Let $(F,\chi)$ be any characteristic $0$ ultraproduct of the structures $(\mathbb{F}_{q},\chi_{q})$ where the $\chi_{q}$ are sufficiently general, i.e., $\chi$ is not of finite order.
%   %  goes to infinity for $q\rightarrow\infty$. 
%     Then, it follows that $F\models \mathrm{PF}^{\times}$.
% \end{theorem}

% \begin{proof}
%     See the proof of the more general Theorem \ref{theoremUPoffinitefieldsmodelofpfplustimes}.
% \end{proof}

\begin{definition}(See 3.2 in \cite{Hrushovski2021AxsTW}.)
    Let $\bar{a}=(a_{1},\dots,a_{n})$ be an n-tuple in $F\models \mathrm{PF}^{\times}$. We denote by $Z(\bar{a})$ the unordered tuple (multiset) containing the elements from $\{x\in F\,|\,x^{n}+a_{1}x^{n-1}+\cdots+a_{n}=0\}$ such that every element occurs with its multiplicity.
    We define the terms $\chi_{\mathrm{sym}}^{n}(y_{1},\dots,y_{n})$ as follows:
    \[\chi_{\mathrm{sym}}^{n}(a_{1},\dots,a_{n}):=\sum_{x\in Z(\bar{a})}\chi(x).\]
    We denote by $\mathcal{L}_{\kappa,\times}^{\mathrm{sym}}$ the language consisting of $\mathcal{L}_{\times}$ together with n-ary predicate symbols symbols $\chi_{\mathrm{sym}}^{n}$ for all $n\in \mathbb{N}$ as well as symbols for the $\kappa_{P,Q}$ as in Definition \ref{definitiondefclosedfunctionsforPF}. The $\mathcal{L}_{\kappa,\times}^{\mathrm{sym}}$-theory $\mathrm{PF}^{\times}$ extends the $\mathcal{L}_{\times}$-theory $\mathrm{PF}^{\times}$ by the interpretations of the new symbols as described above.
\end{definition}

\begin{theorem}
    The theory $\mathrm{PF}^{\times}$ has quantifier elimination in the language $\mathcal{L}_{\kappa,\times}^{\mathrm{sym}}$.
\end{theorem}
\begin{remark}
We did not consider the characteristic $p$ case which could be interesting for multiplicative characters. More precisely, one could ask to determine the theory $\mathrm{Th}((\mathbb{F}_{q},\chi_{q})_{q})\cup\{\mathrm{ord}(\chi)=\infty\}\cup\{\mathrm{char}(F)=p\}$. As the Weil bounds can still be applied, most of the results should work using the same arguments as for $\mathrm{PF}^{+,\times}$. One would however have to check if one can apply the Erdős–Turán–Koksma inequality as in Section \ref{sectionlimittheoryplustimes} when working with the multiplicative group of a pseudofinite field in positive characteristic.
Furthermore, it would be interesting to compare this with the work of Chieu Minh Tran in \cite{tran2019tamestructurescharactersums}. He studies the classical first-order theory of the algebraic closure of a finite field $\bar{\mathbb{F}}_{p}$ together with a circular order stemming from the preimage of the circle under a multiplicative character $\chi:\bar{\mathbb{F}}_{q}\rightarrow S^{1}$. The ingredients from number theory are (almost) the same. As in our work, he makes use of character sum estimates as well as Weyl's equidistribution criterion (we use an effective version in Section \ref{sectionlimittheoryplustimes}). Moreover, the geometric axioms obtained for the corresponding model companion look naturally very similar to ours. The work of Tran embeds in the context of interpolative fusion, studied by Tran together with Kruckmann and Walsberg in \cite{interpolativefusion1} and \cite{kruckman2022interpolativefusionsiipreservation}. We do not know how and if our pseudofinite case relates to this. It could also be interesting to compare \cite{tran2019tamestructurescharactersums} where the theory has $\mathrm{TP}_{2}$ to $\mathrm{Th}(\bar{\mathbb{F}}_{p},\chi)$ where the multiplicative character $\chi$ is a predicate in continuous logic. Is the latter tamer in terms of model-theoretic classification?
\end{remark}

\bibliography{bibpfplustimes}
\bibliographystyle{abbrv}
\end{document}